\documentclass[11pt]{article}

\usepackage{amsfonts,amsthm}
\usepackage{amsmath}   
\usepackage{latexsym}
\usepackage{amssymb}
\usepackage{mathrsfs}
\usepackage{pifont}
\usepackage{tikz}
\usepackage[all]{xy}
\usepackage[hidelinks]{hyperref}

\newtheoremstyle{rem}{3pt}{3pt}{}{}% <Space above><Space below><Body font> <Indent amount>
{\bfseries}{.}{.5em}{}% <Theorem head font><Punctuation after theorem head><Space after theorem headi><Theorem head spec>

\newtheorem{theo}{Theorem}[section]

\newtheorem{conc}{Conclusion}
\theoremstyle{rem}
\newtheorem{remax}[theo]{Remark}
\newenvironment{rema}
  {\pushQED{\qed}\remax}
  {\popQED\endremax}
\theoremstyle{definition}
\newtheorem{defi}[theo]{Definition}
\newtheorem{exam}[theo]{Example}

\newtheorem*{term*}{Notation/Terminology}

\newcommand{\bT}{\boldsymbol{t}}

\newcommand{\bC}{\mathbb{C}}
\newcommand{\si}{\sigma}
\newcommand{\mS}{\mathfrak{S}}

\newcommand{\qc}{\mathrm{c}}

\newcommand{\un}[1]{\underline{#1}}

\newcommand{\bc}{\textbf{c}}

\usepackage{geometry}
\begin{document}

\title{Fusion for the Yang-Baxter equation and the braid group}

\author{L. Poulain d'Andecy\footnote{Laboratoire de Math\'ematiques de Reims UMR 9008, Universit\'e de Reims Champagne-Ardenne,
Moulin de la Housse BP 1039, 51100 Reims, France.
loic.poulain-dandecy@univ-reims.fr}}

\date{}
\maketitle
\thispagestyle{empty}

\begin{abstract}
These are the extended notes of a mini-course given at the school \emph{WinterBraids X}. We discuss algebras simultaneously related to: the braid group, the Yang--Baxter equation and the representation theory of quantum groups. The main goal is to explain the idea of the fusion procedure for the Yang--Baxter equation and to show how it leads to new examples of such algebras: the fused Hecke algebras.
\end{abstract}

\vskip 2cm
\begin{center}
\begin{tikzpicture}[scale=0.5]
\node[above] at (0,2) {\scriptsize{$1$}};
\node[above] at (2,2) {\scriptsize{$2$}};
\node[above] at (8,2) {\scriptsize{$3$}};
\draw[thick] (0,2)--(8,-6);
\draw[thick] (2,2)--(2,-6);
\draw[thick] (8,2)--(0,-6);
\draw (1.4,0.6)--(2,0.7);
\node[above] at (1.5,0.65) {\scriptsize{$u_{12}$}};
\draw (3.6,-1.6)--(4.4,-1.6);
\node[above] at (3.9,-1.6) {\scriptsize{$u_{13}$}};
\draw (2,-3.3)--(2.6,-3.4);
\node[above] at (2.5,-3.35) {\scriptsize{$u_{23}$}};
\node at (10,-2) {$=$};
\node[above] at (12,2) {\scriptsize{$1$}};
\node[above] at (18,2) {\scriptsize{$2$}};
\node[above] at (20,2) {\scriptsize{$3$}};
\draw[thick] (12,2)--(20,-6);
\draw[thick] (18,2)--(18,-6);
\draw[thick] (20,2)--(12,-6);
\draw (18.6,0.6)--(18,0.7);
\node[above] at (18.5,0.65) {\scriptsize{$u_{23}$}};
\draw (15.6,-1.6)--(16.4,-1.6);
\node[above] at (15.9,-1.6) {\scriptsize{$u_{13}$}};
\draw (18,-3.3)--(17.4,-3.4);
\node[above] at (17.5,-3.35) {\scriptsize{$u_{12}$}};
\end{tikzpicture}
\end{center}

\vskip 1cm
\begin{center}
\begin{tikzpicture}[scale=0.5]
%\node at (0,3) {\scriptsize{$i$}};
%\node at (1.8,3) {\scriptsize{$i\!+\!1$}};
%\node at (4,3) {\scriptsize{$i\!+\!2$}};
\draw[line width=0.5mm] (4,2)..controls +(0,-4) and +(0,+4) .. (0,-6);
\draw[line width=0.5mm] (2,2)..controls +(0,-3) and +(0,+3) .. (0,-2);
\fill[white] (0.6,-0.2) circle (0.3);
\fill[white] (2,-2) circle (0.3);
\draw[line width=0.5mm] (0,2)..controls +(0,-4) and +(0,+4) .. (4,-6);
\fill[white] (0.6,-3.8) circle (0.3);
\draw[line width=0.5mm] (0,-2)..controls +(0,-3) and +(0,+3) .. (2,-6);
\node at (6,-2) {$=$};
%\node at (8,3) {\scriptsize{$i$}};
%\node at (9.8,3) {\scriptsize{$i\!+\!1$}};
%\node at (12,3) {\scriptsize{$i\!+\!2$}};
\draw[line width=0.5mm] (12,2)..controls +(0,-4) and +(0,+4) .. (8,-6);
\fill[white] (11.4,-0.2) circle (0.3);
\fill[white] (10,-2) circle (0.3);
\draw[line width=0.5mm] (10,2)..controls +(0,-3) and +(0,+3) .. (12,-2);
\draw[line width=0.5mm] (12,-2)..controls +(0,-3) and +(0,+3) .. (10,-6);
\fill[white] (11.4,-3.8) circle (0.3);
\draw[line width=0.5mm] (8,2)..controls +(0,-4) and +(0,+4) .. (12,-6);
\end{tikzpicture}
\end{center}
\vskip 0cm
\begin{center}
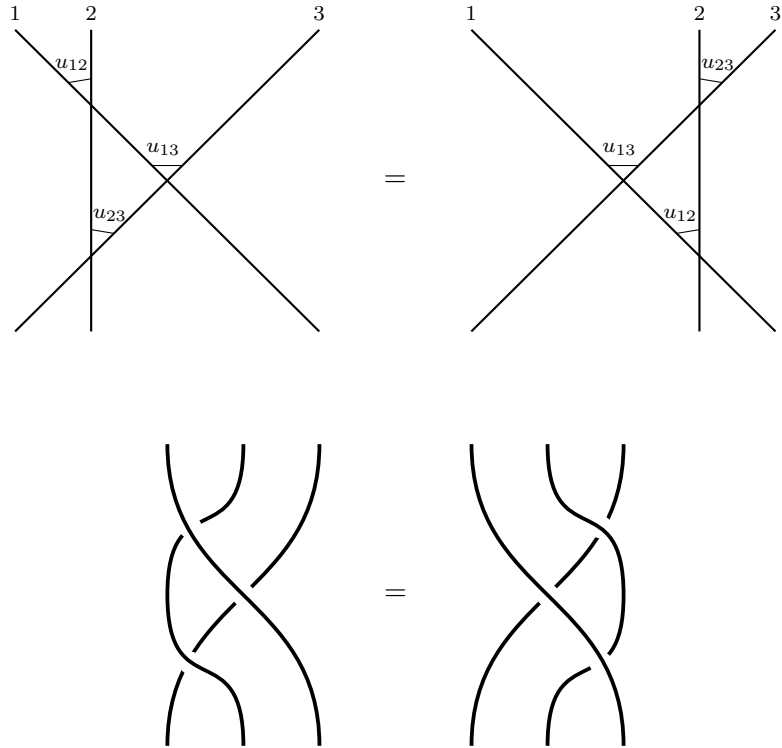
\emph{The Yang--Baxter equation and the braid relation.}
\end{center}

\newpage
\thispagestyle{empty}
\setcounter{tocdepth}{2}
\tableofcontents

\newpage
\section*{Introduction}
\addcontentsline{toc}{section}{Introduction}

The main purpose of these lectures is to discuss algebras with the following features:
\begin{itemize}
\item[$\bullet$] they are quotients of the braid group algebra;
\item[$\bullet$] they contain abstract solutions of the Yang--Baxter equation;
\item[$\bullet$] they admit representations on vector spaces of the form $V^{\otimes n}$.
\end{itemize}
This ``wish list'' of properties is fulfilled in three famous situations, the Temperley--Lieb algebra, the Hecke algebra and the Birman--Murakami--Wenzl algebra. These algebras are well-known quotients of the algebra of the braid group; in fact they are the algebras behind the following invariants of links: the Jones polynomial, the Homfly-pt polynomial and the Kauffman polynomial. 

It is also well-know that they contain solutions of the Yang--Baxter equation (so-called Baxterization formulas). We call them abstract solutions, since they are solutions in some algebras, as opposed to solutions acting on vector spaces (matrix solutions) as genuine solutions of the Yang--Baxter equation should do. To obtain genuine matrix solutions, we need to look for representations of some special form, and these three algebras indeed have such representations on spaces $V^{\otimes n}$. 

The third item on our wish list, the existence of representations on tensor spaces $V^{\otimes n}$, is the key to unifying these examples. In fact, the Temperley--Lieb algebra, the Hecke algebra and the Birman--Murakami--Wenzl algebra (through their actions on tensor spaces) can be seen as centralisers of representations of some interesting algebras: quantum groups. The first goal of these notes is to discuss briefly how the quantum groups (through their centralisers) are nicely designed to produce algebras fulfilling our wish list above.

\vskip .2cm
Once we know that our sought-after algebras can be found among the centralisers of representations of quantum groups, we would like to be able to describe and study these centralisers. The second goal of these notes is to explain how to fit the fusion procedure in this picture. The fusion procedure was designed at the matrix level to produce new solutions of the Yang--Baxter equation. Here we would like to use it in order to understand better the centralisers. This leads us to introduce topological objects that we call fused braids, and to define from them new algebras, called fused Hecke algebras. 

\vskip .2cm
We proceed as follows. The first two sections are introductory to the braid group and the Yang--Baxter equation, and describe in details our favourite examples: the Hecke algebra and the Birman--Murakami--Wenzl algebra. Sections \ref{sec-QG} and \ref{sec-cent} discuss quantum groups and their centralisers. Without full details, they are intended to give an idea of how these objects (quantum groups) are constructed, how they behave and how they are relevant to our subject.

In Section \ref{sec-fusion}, we present the fusion procedure and we revert to a more rigorous presentation for this part. Finally, the last section, presents the fused Hecke algebras. After discussing the braid-like description of these algebras, we go over our wish list, and we check that we indeed produced algebras fulfilling every wishes on the list, and without too much surprise, that these algebras are related to the centralisers of some quantum groups representations.

\vskip .2cm
The lectures were intended to PhD students and young researchers studying braids and their different aspects and applications. In these extended notes, more details are given than during the lectures, but many are still omitted or hidden, and full mathematical precision is not always the priority. Some background in basic algebra and representation theory is assumed.

\paragraph{Aknowledgements.} It is a pleasure to thank the organisers Paolo Bellingeri, Vincent Florens, JB Meilhan, Emmanuel Wagner of the Winterbraids series of schools and the organiser Filippo Callegaro of the session in Pisa, along with all participants of the tenth edition. The work of the author is supported by Agence Nationale de la Recherche Projet AHA ANR-18-CE40-0001.
\newpage

\section{The Braid Group and the Yang--Baxter Equation}\label{sec-def}

\subsection{The Braid group}

For precise definitions and more details on the braid group, we refer to \cite{KT}. Here is an example of a braid (with $5$ strands):
\vskip .4cm
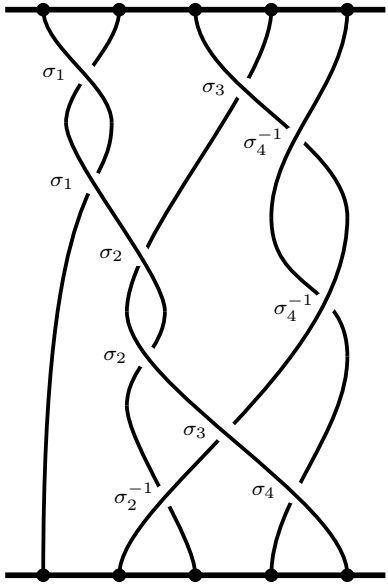
\begin{minipage}{0.3\linewidth}
\begin{center}
\begin{tikzpicture}[scale=0.5]
\draw[line width=0.75mm] (1,0)--(11,0);
\fill[black] (2,0) circle (0.18);\fill[black] (4,0) circle (0.18);
\fill[black] (6,0) circle (0.18);\fill[black] (8,0) circle (0.18);
\fill[black] (10,0) circle (0.18);
\draw[line width=0.5mm] (4,0)..controls +(0,-1) and +(0,+1) .. (2.6,-3);
\node at (2.3,-1.7) {$\scriptstyle{\si_1}$};\fill[white] (3.1,-1.7) circle (0.3);
\draw[line width=0.5mm] (2,0)..controls +(0,-1) and +(0,+1) .. (3.8,-3);
\draw[line width=0.5mm] (3.8,-3)..controls +(0,-2) and +(0,11) .. (2,-15);
\node at (2.5,-4.6) {$\scriptstyle{\si_1}$};\fill[white] (3.3,-4.6) circle (0.3);
\draw[line width=0.5mm] (8,0)..controls +(0,-2) and +(0,+2) .. (4.2,-8);
\node at (3.8,-6.5) {$\scriptstyle{\si_2}$};\fill[white] (4.5,-6.5) circle (0.3);
\draw[line width=0.5mm] (2.6,-3)..controls +(0,-1) and +(0,+1) .. (5.2,-8);
\node at (6.5,-2.1) {$\scriptstyle{\si_3}$};\fill[white] (7.3,-2.1) circle (0.3);
\draw[line width=0.5mm] (6,0)..controls +(0,-2) and +(0,+2) .. (10,-5.5);
\node at (7.8,-3.5) {$\scriptstyle{\si_4^{-1}}$};\fill[white] (8.6,-3.4) circle (0.3);
\draw[line width=0.5mm] (10,0)..controls +(0,-2) and +(0,+2) .. (8,-5.5);
\draw[line width=0.5mm] (5.2,-8)..controls +(0,-1) and +(0,+1) .. (4.2,-10.5);
\node at (3.9,-9.2) {$\scriptstyle{\si_2}$};\fill[white] (4.7,-9.2) circle (0.3);
\draw[line width=0.5mm] (4.2,-10.5)..controls +(0,-1) and +(0,+1) .. (6,-15);
\node at (4.4,-12.9) {$\scriptstyle{\si_2^{-1}}$};\fill[white] (5.2,-12.9) circle (0.3);
\draw[line width=0.5mm] (8,-5.5)..controls +(0,-2) and +(0,+2) .. (10,-9.2);
\node at (8.6,-7.9) {$\scriptstyle{\si_4^{-1}}$};\fill[white] (9.4,-7.8) circle (0.3);
\draw[line width=0.5mm] (10,-9.2)..controls +(0,-2) and +(0,+2) .. (8,-15);
\node at (7.8,-12.8) {$\scriptstyle{\si_4}$};\fill[white] (8.6,-12.8) circle (0.3);
\draw[line width=0.5mm] (10,-5.5)..controls +(0,-4) and +(0,+2) .. (4,-15);
\node at (6,-11.2) {$\scriptstyle{\si_3}$};\fill[white] (6.8,-11.2) circle (0.3);
\draw[line width=0.5mm] (4.2,-8)..controls +(0,-2) and +(0,+2) .. (10,-15);
\draw[line width=0.75mm] (1,-15)--(11,-15);
\fill[black] (2,-15) circle (0.18);\fill[black] (4,-15) circle (0.18);
\fill[black] (6,-15) circle (0.18);\fill[black] (8,-15) circle (0.18);
\fill[black] (10,-15) circle (0.18);
\end{tikzpicture}
\end{center}
\end{minipage}
\hspace{2cm}
\begin{minipage}{0.5\linewidth}
We fix two horizontal parallel lines each having $5$ fixed points, and we connect bijectively each point on the top line to a point on the bottom line by a ``strand''. \\ This is an object in a three-dimensional space, so the strands can pass ``over'' or ``under'' other strands. For now, we may ignore the labels $\sigma_i^{\pm1}$ next to the crossings.\\
Braids are considered up to isotopy, meaning that we can move continuously the strands while leaving their end points fixed, and this is still the same braid.
\end{minipage}

\vskip .4cm
As shown by the example, a braid with $n$ strands is drawn in a rectangular strip with a top line of $n$ fixed dots and a bottom line of $n$ fixed dots. We connect bijectively each top dot to a bottom dot by a strand inside the strip. A strand is a continuous line going from a top point to a bottom point. By convention, we assume that the vertical coordinate is always decreasing along the line (no strand is allowed to go back towards the top). 

At each point of the strip at most two strands are crossing each other, and at each crossing, we indicate which strand passes over the other one. We call a crossing  positive (resp. negative) when the strand coming from the left passes over (resp. under) the strand coming from the right. Such diagram is called a braid with $n$ strands and braids are considered up to isotopy (continuous moves of the strands with fixed end points).

The set of all braids with $n$ strands forms a group, denoted $B_n$, the multiplication being simply the vertical concatenation of diagrams. If $\alpha,\beta\in B_n$, to perform the product $\alpha\beta$, we place the diagram of $\alpha$ on top of the diagram of $\beta$, we identify the bottom dots of $\alpha$ with the top dots of $\beta$, thereby connecting the strands, and we delete these middle dots.

The identity element $1$ of $B_n$ is the braid where the $n$ strands are vertical and parallel. It might not be completely obvious at first sight that the braids form a group. In fact, to find the inverse of a given braid, reflect it through the bottom horizontal line and then invert all the crossings. Up to isotopy, the concatenation with the given braid reduces to the identity element.

The following elements are called elementary braidings:
\begin{center}
 \begin{tikzpicture}[scale=0.3]
\node at (0,0) {$\si_i=$};
\node at (2,3) {$1$};\fill (2,2) circle (0.2);\fill (2,-2) circle (0.2);
\draw[thick] (2,2) -- (2,-2);
\node at (4,0) {$\dots$};
\draw[thick] (6,2) -- (6,-2);
\node at (6,3) {$i-1$};\fill (6,2) circle (0.2);\fill (6,-2) circle (0.2);
\node at (10,3) {$i$};\fill (10,2) circle (0.2);\fill (10,-2) circle (0.2);
\node at (14,3) {$i+1$};\fill (14,2) circle (0.2);\fill (14,-2) circle (0.2);
\draw[thick] (14,2)..controls +(0,-2) and +(0,+2) .. (10,-2);
\fill[white] (12,0) circle (0.4);
\draw[thick] (10,2)..controls +(0,-2) and +(0,+2) .. (14,-2);
\draw[thick] (18,2) -- (18,-2);
\node at (18,3) {$i+2$};\fill (18,2) circle (0.2);\fill (18,-2) circle (0.2);
\node at (20,0) {$\dots$};
\draw[thick] (22,2) -- (22,-2);\fill (22,2) circle (0.2);\fill (22,-2) circle (0.2);
\node at (22,3) {$n$};
\node at (38,0) {$i\in\{1,\dots,n-1\}$;};
\end{tikzpicture}
\end{center}
We can see that sometimes we number the dots $1,...,n$ from left to right to be able to speak of the ``$i$-th dot'' and to help visualize where a crossing is. The previous braid provides an example of positive crossing. As an example of the multiplication rule and of the use of isotopy, one can check that the inverse $\si_i^{-1}$ of $\si_i$ is:
\begin{center}
 \begin{tikzpicture}[scale=0.3]
\node at (-0.5,0) {$\si_i^{-1}=$};
\node at (2,3) {$1$};\fill (2,2) circle (0.2);\fill (2,-2) circle (0.2);
\draw[thick] (2,2) -- (2,-2);
\node at (4,0) {$\dots$};
\draw[thick] (6,2) -- (6,-2);
\node at (6,3) {$i-1$};\fill (6,2) circle (0.2);\fill (6,-2) circle (0.2);
\node at (10,3) {$i$};\fill (10,2) circle (0.2);\fill (10,-2) circle (0.2);
\node at (14,3) {$i+1$};\fill (14,2) circle (0.2);\fill (14,-2) circle (0.2);
\draw[thick] (10,2)..controls +(0,-2) and +(0,+2) .. (14,-2);
\fill[white] (12,0) circle (0.4);
\draw[thick] (14,2)..controls +(0,-2) and +(0,+2) .. (10,-2);
\draw[thick] (18,2) -- (18,-2);
\node at (18,3) {$i+2$};\fill (18,2) circle (0.2);\fill (18,-2) circle (0.2);
\node at (20,0) {$\dots$};
\draw[thick] (22,2) -- (22,-2);\fill (22,2) circle (0.2);\fill (22,-2) circle (0.2);
\node at (22,3) {$n$};
\node at (38,0) {$i\in\{1,\dots,n-1\}$;};
\end{tikzpicture}
\end{center}
It is rather easy to accept that the braid group $B_n$ is generated by the elements $\si_1,\dots,\si_{n-1}$, meaning that any braid can be written as a product of elementary braidings $\si_i$ and their inverses. In the braid pictured above, the elementary braiding are indicated next to each crossing, and we see that this braid is:
\[\si_1\si_3\si_4^{-1}\si_1\si_2\si_4^{-1}\si_2\si_3\si_4\si_2^{-1}\ .\]
As we can see the decomposition of a braid as a product of generators is not unique (for example, we could have ended the word by $\si_2^{-1}\si_4$). This comes from the fact that some relations are satisfied by the generators, reflecting the invariance under isotopy. Namely, in addition to the trivial relations $\si_i\si_i^{-1}=\si_i^{-1}\si_i=1$, it is easy to see that the following relations are satisfied:
\[\si_i\si_j=\si_j\si_i\ \ \ \ \ \text{if $|i-j|>1$,}\] 
since if $|i-j|>1$, then the pairs of strands $i,i+1$ and $j,j+1$ can be manipulated independently. It is also immediate that the following relations are satisfied:
\vskip .1cm
\begin{minipage}{0.3\linewidth}
\begin{center}
\begin{tikzpicture}[scale=0.3]
\node at (0,3) {\scriptsize{$i$}};
\node at (1.8,3) {\scriptsize{$i\!+\!1$}};
\node at (4,3) {\scriptsize{$i\!+\!2$}};
\draw[thick] (4,2)..controls +(0,-4) and +(0,+4) .. (0,-6);
\draw[thick] (2,2)..controls +(0,-3) and +(0,+3) .. (0,-2);
\fill[white] (0.6,-0.2) circle (0.3);
\fill[white] (2,-2) circle (0.3);
\draw[thick] (0,2)..controls +(0,-4) and +(0,+4) .. (4,-6);
\fill[white] (0.6,-3.8) circle (0.3);
\draw[thick] (0,-2)..controls +(0,-3) and +(0,+3) .. (2,-6);
\node at (6,-2) {$=$};
\node at (8,3) {\scriptsize{$i$}};
\node at (9.8,3) {\scriptsize{$i\!+\!1$}};
\node at (12,3) {\scriptsize{$i\!+\!2$}};
\draw[thick] (12,2)..controls +(0,-4) and +(0,+4) .. (8,-6);
\fill[white] (11.4,-0.2) circle (0.3);
\fill[white] (10,-2) circle (0.3);
\draw[thick] (10,2)..controls +(0,-3) and +(0,+3) .. (12,-2);
\draw[thick] (12,-2)..controls +(0,-3) and +(0,+3) .. (10,-6);
\fill[white] (11.4,-3.8) circle (0.3);
\draw[thick] (8,2)..controls +(0,-4) and +(0,+4) .. (12,-6);
\end{tikzpicture}
\end{center}
\end{minipage}
\hspace{0.5cm}
\begin{minipage}{0.6\linewidth}
Indeed, one simply has to move the middle strand. In algebraic terms, we have:
\[\si_i\si_{i+1}\si_i=\si_{i+1}\si_i\si_{i+1}\,,\ \ \ \ \text{for $i=1,\dots,n-2$.}\]
\end{minipage}
\vskip .1cm
\vskip .2cm
It is remarkable that only the two sorts of relations above are enough to completely characterise the braid group. This is the well known following theorem.
\begin{theo}[Artin \cite{Ar1,Ar2,Bo}]\label{theo-artin}
The braid group $B_n$ is generated by $\si_1,\dots,\si_{n-1}$ with defining relations:
\begin{equation}\label{rel-Artin}
\begin{array}{ll}
\si_i\si_{i+1}\si_i=\si_{i+1}\si_i\si_{i+1}\,,\ \ \  & \text{for $i\in\{1,\dots,n-2\}$}\,,\\[0.2em]
\si_i\si_j=\si_j\si_i\,,\ \ \  & \text{for $i,j\in\{1,\dots,n-1\}$ such that $|i-j|>1$}\,.
\end{array}
\end{equation} 
\end{theo}
The phraseology ``defining relations'' means that any other relation involving braids in $B_n$ is implied by these ones (and the trivial ones $\si_i\si_i^{-1}=\si_i^{-1}\si_i=1$). Equivalently, it means that any other group with generators satisfying the relations (\ref{rel-Artin}) must be a quotient of the braid group $B_n$.

\paragraph{Local representations of $B_n$.} 
Let $V$ be a vector space and form the tensor product $V^{\otimes n}$ of $n$ copies of $V$. We take an invertible element $R\in\text{End}(V\otimes V)$ and we define elements $R_1,\dots,R_{n-1}$ like this:
\[\text{Operators $R_1,\dots,R_{n-1}$ on $V^{\otimes n}$:}\qquad\ \ \ \rlap{$\overbrace{\phantom{V\otimes V}}^{R_1}$}V\otimes\overbrace{V\otimes V}^{R_2}\otimes\dots\overbrace{\ldots\otimes V}^{R_{n-1}}\]
More precisely, for $i=1,\dots,n-1$, the operator $R_{i}$ acts on $V^{\otimes n}$ non-trivially only in the copies $i,i+1$ of $V$, where it acts by $R$. That is, it is defined by $R_{i}=\text{Id}_{V^{\otimes i-1}}\otimes R\otimes\text{Id}_{V^{\otimes n-i-1}}$.

Now we say that $R$ defines a local representation of the braid group $B_n$, or that the following map
\[\begin{array}{lcrcl}\rho\ &:&\ B_n & \to & \text{End}(V^{\otimes n})\\[0.5em]
 & & \si_i & \mapsto & R_i \end{array}\]
is a local representation of $B_n$ on the tensor space $V^{\otimes n}$, if the braid relation is satisfied
\begin{equation}\label{braid}
R_1R_{2}R_1=R_{2}R_1R_{2}\ \quad \ \ \ \ \text{on $V\otimes V\otimes V$.} 
\end{equation}
This immediately implies all the other defining relations of $B_n$, thereby making the map $\rho$ to be a representation of the group $B_n$. 
Note that the other relation $R_iR_j=R_jR_i$ if $|i-j|>1$ follows from the locality of the operators $R_i,R_j$: they do not act on the same copies of $V$ if $|i-j|>1$.

A local representation is thus a representation of $B_n$ on a tensor space $V^{\otimes n}$ with the special (local) form explained above for the action of the generators $\si_i$. It depends only on the element $R\in\text{End}(V\otimes V)$ satisfying (\ref{braid}) (and invertible in order to have a representation of a group). 

\begin{rema}
A direct, or numerical, approach to find local representation of the braid group $B_n$ is quite difficult. Indeed if $\dim(V)=n$ then one has to solve cubic equations in $n^4$ variables.
\end{rema}

\begin{exam}
Let $P\in\text{End}(V\otimes V)$ be the permutation operator sending $x\otimes y$ to $y\otimes x$. Then it provides a local representation of the braid group. In fact, it factors through the natural permutation representation of the symmetric group on $V^{\otimes n}$. All this amounts to the simple facts that $P^2=\text{Id}$ and that the following equality is true when one composes transpositions: $(1,2)(2,3)(1,2)=(2,3)(1,2)(2,3)$.
\end{exam}

\subsection{The Yang--Baxter equation}\label{subsec-YB}

Let $V$ be a finite-dimensional vector space. We introduce a bit more notations for operators on $V^{\otimes n}$. If $S\in \text{End}(V\otimes V)$, we use the notation $S_{ij}$ for the operator on $V^{\otimes n}$ acting as $S$ on copies $i$ and $j$ and trivially otherwise. More formally, write $S=\sum_as_a\otimes t_a$\,, where $s_a,t_a\in\text{End}(V)$. Then by definition, we have:
\[S_{ij}=\sum_a \text{Id}_V\otimes\dots\otimes\text{Id}_V\otimes s_a\otimes\text{Id}_V\otimes\dots \otimes\text{Id}_V\otimes t_a \otimes \text{Id}_V\otimes \dots \otimes \text{Id}_V\ ,\]
where $s_a$ is in position $i$ and $t_a$ is in position $j$. We use this notation in particular for the permutation operator $P$ of $V\otimes V$. Explicitly, we have:
\[P_{ij}(v_1\otimes\dots\otimes v_i\otimes\dots \otimes v_j \otimes \dots \otimes v_n)=v_1\otimes\dots\otimes v_j\otimes\dots \otimes v_i \otimes \dots \otimes v_n\ .\]
The basic property about all these notations that we use repeatedly is:
\[S_{ab}P_{ij}=P_{ij}S_{\pi_{i,j}(a)\pi_{i,j}(b)}\,,\ \ \ \ \ \ \ \text{where $\pi_{i,j}$ is the transposition $(i,j)$.}\]
For example, $S_{13}P_{12}=P_{12}S_{23}$ or  $S_{13}=P_{23}S_{12}P_{23}$, etc.

\vskip .2cm
Now we consider a function 
\[\begin{array}{lcrcl}R\ &:&\ \bC^2 & \to & \text{End}(V\otimes V)\\[0.5em]
 & & (u,v) & \mapsto & R(u,v) \end{array}\ .\]
The Yang--Baxter equation (YB equation for short) is a functional equation for such a function $R$:
\begin{equation}\label{YBprime}
R_{12}(u_1,u_2)R_{13}(u_1,u_3)R_{23}(u_2,u_3)=R_{23}(u_2,u_3)R_{13}(u_1,u_3)R_{12}(u_1,u_2)\ \quad \ \ \ \ \text{on $V\otimes V\otimes V$.} 
\end{equation}
In this context, the variables $u_1,u_2,u_3$ are often called \emph{spectral parameters}, or \emph{spectral variables}.

\vskip .2cm
In our perspective, the so-called braided version of this equation is more relevant. From $R$ define another function:
\[\begin{array}{lcrcl}\check R\ &:&\ \bC^2 & \to & \text{End}(V\otimes V)\\[0.5em]
 & & (u,v) & \mapsto & PR(u,v) \end{array}\ ,\] 
 where $P$ is the permutation operator sending $x\otimes y$ to $y\otimes x$. To see better the connections with the braid relation, set $\check R_{1}(u,v):=\check R_{12}(u,v)$ and $\check R_{2}(u,v):=\check R_{23}(u,v)$. Then an easy manipulation shows that the YB equation is equivalent to the following equation only involving the function $\check R$:
\begin{equation}\label{YB}
\check R_{1}(u_1,u_2)\check R_{2}(u_1,u_3)\check R_{1}(u_2,u_3)=\check R_{2}(u_2,u_3)\check R_{1}(u_1,u_3)\check R_{2}(u_1,u_2)\ \quad \ \ \ \ \text{on $V\otimes V\otimes V$.} 
\end{equation}
When precision is needed, we will refer to this version as the braided YB equation.

It happens in many situations that the function $\check R$ (or equivalently $R$) depends only on the ratio $u/v$. In this particular case (set $u=u_1/u_2$ and $v=u_2/u_3$), the braided YB equation becomes:
\begin{equation}\label{YBmult}
\check R_1(u)\check R_{2}(uv)\check R_1(v)=\check R_{2}(v)\check R_1(uv)\check R_{2}(u)\ \quad \ \ \ \ \text{on $V\otimes V\otimes V$.} 
\end{equation}
Instead of a multiplicative version, we sometimes also have the additive one, where the function $\check R$ depends only on the difference $u-v$:
\begin{equation}\label{YBadd}
\check R_1(u)\check R_{2}(u+v)\check R_1(v)=\check R_{2}(v)\check R_1(u+v)\check R_{2}(u)\ \quad \ \ \ \ \text{on $V\otimes V\otimes V$.} 
\end{equation}

\begin{exam}[\textbf{constant solution, Yang solution}]\label{ex-solYB}$\ $\\
$\bullet$ For a first encounter with the YB equation, one may ask for constant solutions. The equation (\ref{YB}) for a constant function is simply the braid relation. Thus, a constant solution of the braided YB equation is equivalent to a local representation of the braid group, as defined previously.

\noindent $\bullet$ One can check by hand that the following function satisfies (\ref{YBprime}):
\[R(u,v)=\text{Id}_{V\otimes V}+\frac{P}{u-v}\,,\]
or equaivalently that the braided YB equation (\ref{YB}) is satisfied by:
\[\check R(u,v)=PR(u,v)=P+\frac{\text{Id}_{V\otimes V}}{u-v}\,.\] 
Note that in this verification, one uses that $P$ satisfies the braid relation (this is an example of \emph{Baxterization}, as we will see later). This solution is an example depending on the spectral parameters only through their difference. Namely we have that $\check R(u)=P+\frac{\text{Id}_{V\otimes V}}{u}$ satisfies (\ref{YBadd}).
\end{exam}

\subsubsection{The Yang--Baxter equation in Physics}

The Yang--Baxter equation is certainly one of the fundamental equations in theoretical and mathematical physics. Its history is complex and naturally most of it requires some background in physics that we are not ready to discuss here. So let us only give some hints and vague indications which might be enough to get a feeling of its importance and to feel a little motivated for its study. Some references among the huge literature on this subject are \cite{AYP,Bax,Fad,FST,GRS,Ji-int,Ji-ed,JM,KBI}. Though the YB equation belongs to the theoretical and mathematical side of Physics, it is interesting to note that connections with more experimental physics can be investigated today, see for example \cite{Bat,BF,GBL}.

\paragraph{Factorization of interaction processes.} Probably the most mathematicians-friendly interpretation of the YB equation in physics is as a factorisation property. Let us depict the interaction of two particles called $1$ and $2$ like this: 
\begin{equation}\label{YB-traj}
\begin{tikzpicture}[scale=0.3]
\node[above] at (0,0) {\scriptsize{$1$}};
\node[above] at (2,0) {\scriptsize{$2$}};
\draw[thick] (0,0)--(2,-2);
\draw[thick] (2,0)--(0,-2);
\draw (0.6,-0.6)--(1.4,-0.6);
\node[above] at (0.95,-0.6) {\scriptsize{$u$}};
\end{tikzpicture}\ \ \ \ :\ \ \ \ \text{space-time trajectories for the interactions of two particles,}
\end{equation}
where we think of the vertical direction as time and of the horizontal one as space. The parameter $u$ is the difference of rapidity of the two particles\footnote{In relativity, the rapidity corresponding to the speed $v$ is $\text{tanh}^{-1}(\frac{v}{c})$ and is additive for one-dimensional motion ($c$ is the speed of light).}. If we think of the particles as being described by a vector space $V$ (their internal state space), the interaction is then controlled by an operator, say $S(u)$, in $\text{End}(V\otimes V)$.

Now consider $n$ particles on the line (numbered from $1$ to $n$) and they all interact with each other. We imagine our $n$ particles  coming closer and closer to each other, interacting when they cross each other (controlled by the scattering operator $S(u)$), and after having all interacted once with each other, their order is completely reversed. As an example, if we consider three particles, then we find two different possibilities for ordering these interactions. The space-time trajectories are shown below:
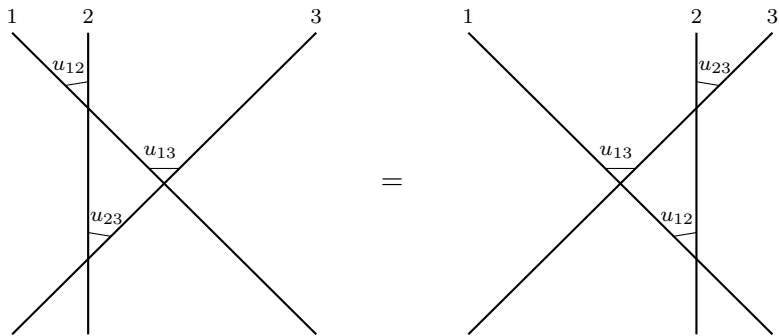
\begin{center}
\begin{tikzpicture}[scale=0.5]
\node[above] at (0,2) {\scriptsize{$1$}};
\node[above] at (2,2) {\scriptsize{$2$}};
\node[above] at (8,2) {\scriptsize{$3$}};
\draw[thick] (0,2)--(8,-6);
\draw[thick] (2,2)--(2,-6);
\draw[thick] (8,2)--(0,-6);
\draw (1.4,0.6)--(2,0.7);
\node[above] at (1.5,0.65) {\scriptsize{$u_{12}$}};
\draw (3.6,-1.6)--(4.4,-1.6);
\node[above] at (3.9,-1.6) {\scriptsize{$u_{13}$}};
\draw (2,-3.3)--(2.6,-3.4);
\node[above] at (2.5,-3.35) {\scriptsize{$u_{23}$}};
\node at (10,-2) {$=$};
\node[above] at (12,2) {\scriptsize{$1$}};
\node[above] at (18,2) {\scriptsize{$2$}};
\node[above] at (20,2) {\scriptsize{$3$}};
\draw[thick] (12,2)--(20,-6);
\draw[thick] (18,2)--(18,-6);
\draw[thick] (20,2)--(12,-6);
\draw (18.6,0.6)--(18,0.7);
\node[above] at (18.5,0.65) {\scriptsize{$u_{23}$}};
\draw (15.6,-1.6)--(16.4,-1.6);
\node[above] at (15.9,-1.6) {\scriptsize{$u_{13}$}};
\draw (18,-3.3)--(17.4,-3.4);
\node[above] at (17.5,-3.35) {\scriptsize{$u_{12}$}};
\end{tikzpicture}
\end{center}
where $u_{ij}=u_i-u_j$ are the differences of the rapidities. The Yang--Baxter equation is simply the hypothesis that these two interacting processes lead to the same result. Namely, in operator notation, we find:
\[S_{12}(u_{12})S_{13}(u_{13})S_{23}(u_{23})=S_{23}(u_{23})S_{13}(u_{13})S_{12}(u_{12})\ ,\]
which is the Yang--Baxter equation for the function $S\ :\ (u,v)\mapsto S(u-v)$.

The fundamental fact about this is the following. Once we assume the YB equation for the interaction of three particles, then it follows that the interaction of $n$ particles, decomposed as a sequence of $2$-particles interactions, is independent of the chosen sequence. In other words, the $n$ particles interaction is unambiguously given by the two particles interaction.  The YB equation is a compatibility condition for the factorization of the $n$-body interaction in terms of the two-body interaction. This originates in the papers \cite{McG,Yang1,Yang2,ZZ} and accounts for the first half of the name of the equation.

\begin{rema}
The above interpretation of the YB equation is very similar to the following classical fact about the braid group and the symmetric group. The longest element of the symmetric group (reversing the order of $1,\dots,n$) can be factorised in different ways as a minimal-length product of simple transpositions, and all these different factorisations are seen to be equal just by assuming the braid relation. The YB equation above plays the same role as the braid equation here.
\end{rema}

\paragraph{2-dimensional statistical physics.} The YB equation has a long history in 2-dimensional statistical models, and we shall not attempt to describe it. Suffices it to say that it seems to have made its first apparition (somewhat hidden) in the solution of the Ising model by Onsager in the 40's \cite{Ons} and its importance was gradually recognized, culminating in the work of Baxter in the 70's \cite{Bax}. It is instructing to search for explicit apparitions of the YB equation in the classical book by Baxter \cite{Bax}. The first two are located near the end of the discussions on the Ising model (\S 7.13) and on the ice model (\S 9.7). Remarkably, the next one (\S10.4) appears at the beginning of the discussion of the 8-vertex model. The meaning is clear: the decisive step of promoting the YB equation as the key to solvability of these models was done. At that time, the YB equation was called the star-triangle relation, or the parametrized star-triangle relation to emphasize the presence of the spectral parameter.

\vskip .2cm
In these statistical models, there is a partition function which is the thing one would like to calculate. The partition function can be expressed as the trace of some powers of an object called the transfer matrix $t(u)$ depending on a parameter. So of course, one natural approach is to try to diagonalize this matrix. Then the key property that one wishes for this transfer matrix is its commutation property:
\[[t(u),t(v)]=0\ \ \ \ \ \forall u,v\,,\]
since it provides an infinite number of commuting operators. The YB equation (for which operator? we shall see below) was identified as a condition ensuring the commutativity of the transfer matrices.

After this very brief summary, let us at least explain how an operator in $\text{End}(V\otimes V)$ arises in this setting. A vertex model in 2-dimensional statistical physics starts with a lattice in the plane, say a square lattice for definiteness. Each edge can be in a certain state $s$ taken in a given finite set $S$. For example $S$ can be of cardinal $2$, corresponding to the two orientations of a spin, or two states ``empty'' and ``occupied''. The interactions happen at the vertices (hence the name ``vertex models'') and at each vertex, there is a Boltzmann weight\footnote{In statistical physics, the Boltzmann weight of a state is equal to $e^{-E/k_BT}$ where $E$ is the energy, $T$ is the temperature and $k_B$ is the Boltzmann constant. The Boltzmann weight expresses the probability for a system to be in a state with energy $E$.} depending on the states of the edges connected to this vertex. To define the partition function, take a configuration of the lattice (fix a state for each edge) and take the product of the Boltzmann weights of all vertices. Then sum over all the  possible configurations of the lattice.

Thus to define the model, it suffices to give the Boltzmann weight of a vertex for each configuration of the edges. If we organise these weights in an array as follows:
\begin{center}
\begin{tikzpicture}[scale=0.25]
\draw[thick] (5,0)--(15,0);
\draw[thick] (10,5)--(10,-5);
\node at (6,1) {$\scriptstyle{i}$};
\node at (14,1) {$\scriptstyle{j}$};
\node at (11,4) {$\scriptstyle{l}$};
\node at (11,-4) {$\scriptstyle{k}$};
\node at (25,0) {$\leadsto$};
\node at (35,0) {$\bigl(R_{ik}^{jl}\bigr)_{i,j,k,l\in S}$};
\end{tikzpicture}
\end{center}
we realise that it is an array of numbers with $4$ indices in $S$. Thus this can be seen as an element of $\text{End}(V\otimes V)$, where $V$ is a vector space with basis indexed by $S$. The appearance of the spectral parameter $u$ in this context is not easily explained, it can be seen as a clever way of parametrizing the Boltzmann weights (one can start from the occurrences of the YB equation in \cite{Bax} previously indicated, and read backwards to see where the spectral parameter came from). Some symmetries can be required, restricting the values of the various weights. For example if $\dim(V)=2$, for the 6-vertex model (respectively, 8-vertex model), only 6 entries (respectively, 8 entries) of the $R$-matrix are non-zero.

If it happens that this operator compiling the various Boltzmann weights satisfies the YB equation, then the commuting property of the transfer matrix is ensured.

\paragraph{Quantum spin chains.} A quantum spin chain is a quantum-mechanical model of interacting particles on a line. The Hilbert space of states is the tensor product $V^{\otimes n}$ ($n$ particles, each having $V$ as its Hilbert space of states), and there is a Hamiltonian describing the interactions. For example, $V$ can be $\mathbb{C}^2$ and the Hamiltonian can be:
\[\sum_{i=1}^n(J_x\sigma_i^x\sigma_{i+1}^x+J_y\sigma_i^y\sigma_{i+1}^y+J_z\sigma_i^z\sigma_{i+1}^z)\,,\]
where $\sigma^x,\sigma^y,\sigma^z$ are the famous Pauli matrices and the indices $i,i+1$ indicate the positions in the tensor product $V^{\otimes n}$ ($n+1$ is understood as $1$). This is called the XYZ model. It reduces to the XXZ model if we take $J_x=J_y$ and further to the $XXX$ model (or Heisenberg spin chain) if $J_x=J_y=J_z$.

In quantum mechanics, the goal is clear: we want the eigenvalues and eigenvectors of the Hamiltonian. The connections with the seemingly unrelated statistical models of the previous paragraph is as follows. The transfer matrix of a vertex model is also an operator on some tensor product $V^{\otimes n}$. And it was found that for the ice model, the eigenvectors of the transfer matrix are the same as for the $XXZ$ Hamiltonian. In fact, it turned out that the Hamiltonian commutes with the transfer matrix $t(u)$, thus explaining this coincidence. Even better, it was finally realized that the Hamiltonian is included in a sense in the transfer matrix (the Hamiltonian can be recovered from $t(u)$), and this already for the $XYZ$ model (the corresponding vertex model is the 8-vertex model). So in fact quantum spin chains also fall in the realm of models governed by solutions of the YB equation.

\paragraph{Algebraic Bethe Ansatz.} In the end of the 70s, a general method to construct and study integrable systems was developed under the name of Algebraic Bethe Ansatz (see \cite{Fad,RTF} and references therein). An essential ingredient taken as a starting point is an $R$-matrix, namely a solution of the YB equation. Long story short, an $L$ matrix, a monodromy matrix and a transfer matrix are successively constructed from it, and ultimately, the commutation relation for the transfer matrix follows from the YB equation.

This approach put forward the algebraic relations between the various operators involved (including of course the YB equation), and it naturally led to wonder about abstract and general algebraic structures behind all this. This was the birth of quantum groups, and the moment in time where the YB equation started to diffuse into mathematics. Quantum groups are one of our next subjects of discussions. It is remarkable that these structures (quantum groups) originating from considerations in mathematical physics turn out to be very important in modern representation theory, for example for the symmetric group in positive characteristic. We will not talk about that.

\newpage
\section{The Guiding Examples: Hecke and Birman--Murakami--Wenzl}\label{sec-ex}

\subsection{The Hecke algebra}\label{subsec-Hec}

There are several possible equivalent definitions for the Hecke algebra, in quite different contexts. The most natural here is as a quotient of the algebra of the braid group. The idea is as follows. The braid group is a fairly complicated group algebraically since we can make an arbitrary number of crossing between strands. There is no way to reduce this number of crossings (except in the trivial situation when a positive crossing meets a negative crossing on the same two strands, that is, when we can use the relation $\si_i\si_i^{-1}=1$). So the naive idea is to add an algebraic relation allowing to deal with all these crossings. One of the simplest way to do this results in the Hecke algebra.

More precisely, we add to the definition of the braid group the following relation:
\begin{equation}\label{rel-skeinH}
 \begin{tikzpicture}[scale=0.25]
\draw[thick] (0,2)..controls +(0,-2) and +(0,+2) .. (4,-2);
\fill[white] (2,0) circle (0.4);
\draw[thick] (4,2)..controls +(0,-2) and +(0,+2) .. (0,-2);
\node at (6,0) {$=$};
\draw[thick] (12,2)..controls +(0,-2) and +(0,+2) .. (8,-2);
\fill[white] (10,0) circle (0.4);
\draw[thick] (8,2)..controls +(0,-2) and +(0,+2) .. (12,-2);
\node at (17,0) {$-\,(q-q^{-1})$};
\draw[thick] (21,2) -- (21,-2);\draw[thick] (25,2) -- (25,-2);
\end{tikzpicture}
\end{equation}
where $q$ is a parameter, which can be seen as an indeterminate or as a non-zero complex number. This relation is a local relation, meaning that for any braid and any of its crossing, the braid is equal to a sum of two terms: the braid obtained by replacing the crossing by its opposite and $\pm(q-q^{-1})$ (depending on the sign of the original crossing) times the braid obtained by replacing the crossing by two pieces of vertical strands. In particular this allows one to transform all the negative crossings into positive ones. 

Algebraically, the new relation is equivalent to imposing $\sigma_i^{-1}=\sigma_i-(q-q^{-1})$ for all generators $\si_i$. Of course, if $q^2\neq1$ we leave the realm of groups with such a relation, and we end up with an algebra. So the algebraic definition of the Hecke algebra $H_n(q)$ goes by defining it as the algebra generated by elements $\si_1,\dots,\si_{m-1}$ with defining relations:
\begin{equation}\label{rel-H}
\begin{array}{ll}
\si_i\si_{i+1}\si_i=\si_{i+1}\si_i\si_{i+1}\,,\ \ \  & \text{for $i\in\{1,\dots,n-2\}$}\,,\\[0.2em]
\si_i\si_j=\si_j\si_i\,,\ \ \  & \text{for $i,j\in\{1,\dots,n-1\}$ such that $|i-j|>1$}\,,\\[0.2em]
\si_i^2=1+(q-q^{-1})\si_i\,,\ \ \  & \text{for $i\in\{1,\dots,n-1\}$}\,.
\end{array}
\end{equation} 
It is an algebra over $\bC$ if $q$ is a complex number, and can be defined also as an algebra over $\bC[q,q^{-1}]$ if $q$ is an indeterminate. Comparing with the algebraic presentation of the braid group $B_n$ in Theorem \ref{theo-artin}, we see that the Hecke algebra $H_n(q)$ is a quotient of the group algebra of the braid group $B_n$ (the algebra $\bC B_n$ consisting of linear combinations of elements of $B_n$). The characteristic equation of degree 2 for the generators is:
\[(\si_i-q)(\si_i+q^{-1})=0\ .\]
If $q^2=1$ the relations (\ref{rel-H}) are defining relations for the symmetric group $S_n$, and the Hecke algebra $H_n(1)$ is thus the group algebra $\bC S_n$. In this case, the generators $\si_i$ correspond to the transpositions $(i,i+1)$. It should feel natural since if $q^2=1$, the local Hecke relation says that we can forget about the signs of the crossings, and so all topological information is lost, and what remains of a braid is simply the permutation of the $n$ points induced by the strands. 

Thus the Hecke algebra is a deformation of the group algebra of the symmetric group $S_n$. The symmetric group can be seen as a ``simplification'' of the braid group, but it is in some sense ``too simple''. The Hecke algebra retains more information (in particular, topological) because of the additional freedom given by the deformation parameter $q$. 

\begin{rema}[Jones polynomial]
A breakthrough in the theory of invariants of links was the discovery of the Jones polynomial (and its generalisation, the so-called HOMFLY-PT polynomial) \cite{HOMFLY,Jo-pol,PT}. In algebraic terms, this invariant is obtained from a Markov trace on the chain of Hecke algebras (see for example \cite[\S4.5]{GP}). 
\end{rema}

\begin{rema}
There are several related algebraic structures called Hecke algebras:\\
$\bullet$ One purely algebraic definition is as deformations of Coxeter groups. Our example corresponds to a finite Coxeter group of type $A$, \emph{a.k.a.} the symmetric group.\\
$\bullet$ There is also a very general definition of Hecke algebras (explaining the name ``Hecke'') as endomorphism algebras of induced representations. We know of several coincidences of the two definitions. Let us only indicate that if we take $G$ the group of invertible matrices over a finite field with $q$ elements, and induce to $G$ the trivial representation of the subgroup consisting of upper-diagonal matrices, then the endomorphism algebra turns out to be the Hecke algebra $H_n(q)$.\\
$\bullet$ In our case of $H_n(q)$, there is another possible definition as the centraliser of the action of a quantum group in a tensor product of representations. We will come back to this later.
\end{rema}

\paragraph{The YB equation.} For any $i\in\{1,\dots,n-1\}$, define the following function taking values in the algebra $H_n(q)$:
\begin{equation}\label{Baxt-H}
\si_i(u)=\si_i+(q-q^{-1})\frac{1}{u-1}\ .
\end{equation}
Then a straightforward calculation using the relations in $H_n(q)$ shows that:
\[\si_i(u)\si_{i+1}(uv)\si_i(v)=\si_{i+1}(v)\si_i(uv)\si_{i+1}(u)\ ,\]
namely, the braided YB equation is satisfied in $H_n(q)$ by these functions.

So we have a solution of the braided YB equation inside the Hecke algebra $H_n(q)$, given by a rather simple formula. However this is not a \emph{genuine} solution of the YB equation since it is not yet an operator on a tensor product of vector spaces. The ingredient one has to add is a little bit of representation theory, namely a representation of $H_n(q)$ on a tensor space $V^{\otimes n}$. More precisely, we want a local representation of $H_n(q)$, that is, a representation of the form:
\[\sigma_i\mapsto \check R_{i}=\text{Id}_{V^{\otimes i-1}}\otimes \check R\otimes\text{Id}_{V^{\otimes n-i-1}}\,,\] 
for some operator $\check R\in\text{End}(V\otimes V)$. Indeed, assume that we have such a representation and set:
\[\check R(u)=\check R+(q-q^{-1})\frac{\text{Id}_{V\otimes V}}{u-1}\ .\]
Then it follows from (\ref{Baxt-H}) that $\check R(u)$ satisfies the braided YB equation.

\vskip .2cm
This is all very nice, only if we can find local representations of the Hecke algebra $H_n(q)$. It turns out that there are some. In fact, take any vector space $V$, then we can construct a local representation of $H_n(q)$ on $V^{\otimes n}$. The rough idea is that for the symmetric group $S_n$, there is one local representation which is simply by permuting the $n$ components in the tensor product. As the Hecke algebra $H_n(q)$ is a deformation of $S_n$, one may try to deform the permutation representation. It works.

More precisely, fix a basis $(e_1,\dots,e_N)$ of $V$ and define a linear operator on $V\otimes V$ by:
\begin{equation}\label{rep-Hn}\check R(e_a\otimes e_b):=\left\{\begin{array}{ll}
q\,e_a\otimes e_b\ \  & \text{if $a=b$,}\\[0.8em]
e_b\otimes e_a+(q-q^{-1})\,e_a\otimes e_b & \text{if $a<b$,}\\[0.4em]
e_b\otimes e_a & \text{if $a>b$.}
\end{array}\right.\ \ \ \text{where $a,b=1,\dots,N$.}
\end{equation}
One can check by hand that this provides a local representation of the braid group: $\check R_i\check R_{i+1}\check R_i=\check R_{i+1}\check R_i\check R_{i+1}$; and moreover, that it factors through the Hecke algebra since the Hecke relation $\check R^2=\text{Id}_{V\otimes V}+(q-q^{-1})\check R$ is satisfied. So finally we get a solution of the braided YB equation on any vector space $V$ by the above procedure.
\begin{exam}\label{ex-solYBq}
For example, if $\dim(V)=2$, in the lexicographic ordering of the basis of $V\otimes V$, we find:
\[\check R=\left(\begin{array}{cccc}
q & \cdot & \cdot & \cdot \\
\cdot  & q-q^{-1} & 1& \cdot \\
\cdot  & 1 & 0& \cdot \\
\cdot & \cdot & \cdot & q
\end{array}\right)\ ,\ \ \ \ \ \ 
\check R(u)=\left(\begin{array}{cccc}
\displaystyle\frac{qu-q^{-1}}{u-1} & \cdot & \cdot & \cdot \\
\cdot   & \displaystyle\frac{(q-q^{-1})u}{u-1}& 1& \cdot \\
\cdot  & 1 & \displaystyle\frac{q-q^{-1}}{u-1} & \cdot \\
\cdot & \cdot & \cdot & \displaystyle\frac{qu-q^{-1}}{u-1}
\end{array}\right)\ .\]
The solution $\check R(u)$ coming from the Hecke algebra is a deformation of the Yang solution $\check R(\alpha)=P+\frac{\text{Id}}{\alpha}$ from Example \ref{ex-solYB} in the following sense: set $u=q^{2\alpha}$ and take the limit $q\to 1$.
\end{exam}

\subsection{The Birman--Murakami--Wenzl algebra}

The Hecke algebra is more or less the most general quotient of the algebra of the braid group by a quadratic characteristic equation for the generators. So following the idea of simplifying the braid group, it is natural to turn our attention to quotients by a cubic characteristic equation.

\paragraph{On quotients of the braid group.} A naive idea for defining quotients of the braid group algebra $\bC B_n$ is to add the following relations:
\[\si_i^m=a_0+a_1\si_i+\dots+a_{m-1}\si_i^{m-1}\ \ \ \text{for $i=1,\dots,n-1$.}\]
In this context, we would like the generators $\si_i$ to stay invertible so that, up to a renormalisation of the generators $\si_i$ which does not change the braid relations, we can assume that $a_0=1$. Note that all generators $\si_i$ are conjugated in the braid group so they must satisfy the same characteristic equation (hence the same relation for all $i=1,\dots,n-1$). If $m=2$, up to some conventions, the resulting algebra is the Hecke algebra $H_n(q)$.

A first natural step is to look at the group we obtain if we consider the relation $\si_i^m=1$. Then we could try to see the algebra resulting from the relations above as deformations of these groups. This would certainly be nice if the groups (respectively, the algebras) were finite (respectively, finite-dimensional). This question has been settled by Coxeter \cite{Co} for the groups resulting in the following remarkable result: if we make the quotient of the braid group $B_n$ on $n$ strands by the relation $\si_i^m=1$ then the resulting group is finite if and only if $\frac{1}{n}+\frac{1}{m}>\frac{1}{2}$. This gives the following possibilities:
\[(n,m)\in\bigl\{\ (2,m)\,,\ \ (n,2)\,,\ \ (3,3)\,,\ \ (3,4)\,,\ \ (3,5)\,\ \ (4,3)\,,\ \ (5,3)\ \bigr\}\ .\]
The first two families correspond to cyclic groups and symmetric groups. Aside from these groups, there are 5 other situations giving a finite group.

The second step is to try to deform the above finite groups by considering an arbitrary characteristic equation as above. It turns out that there is such a deformation theory for these groups, since they all belong to the family of finite complex reflection groups. The resulting algebra are called cyclotomic Hecke algebras, and they are flat deformations in the sense that their dimensions remains equal to the order of the corresponding group. We refer to \cite{Ma1,Ch}.

\paragraph{An example of a cubic quotient: the BMW algebra.} As we have discussed just above, a generic cubic quotient of the braid group algebra does not have to be finite-dimensional (as soon as $n>5$). So an idea is to again add relations to make it finite-dimensional. Complete study of such quotients turns out to be quite involved, and it seems to be still an open question what would be the most generic finite-dimensional quotient of the braid group algebra including a cubic characteristic equation for the generators, see for example \cite{Ma2}. The Birman--Murakami--Wenzl (BMW) algebra \cite{BW,Mu1} is the most well-known particular case of such a quotient. 

The topological definition of the BMW algebra uses the notion of tangles, generalising braids. We keep the rectangular strip with the two lines of $n$ dots, and we still connect each dot to another dot by a strand inside the strip. However now we are allowed to connect a top dot to another top dot (and similarly for bottom dots). An example is:
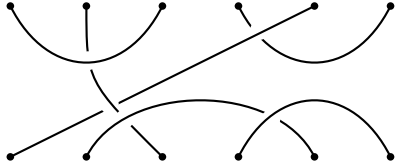
\begin{center}
\begin{tikzpicture}[scale=0.5]
\fill[black] (0,0) circle (0.1);
\fill[black] (2,0) circle (0.1);
\fill[black] (4,0) circle (0.1);
\fill[black] (6,0) circle (0.1);
\fill[black] (8,0) circle (0.1);
\fill[black] (10,0) circle (0.1);
\draw[thick] (6,0)..controls (7,-2) and (9,-2) .. (10,0);
\fill[white] (6.55,-0.65) circle (0.25);
\draw[thick] (8,0)--(0,-4);
\fill[white] (2.6,-2.6) circle (0.25);
\draw[thick] (2,0)..controls (2,-2) .. (4,-4);
\fill[white] (2.05,-1.45) circle (0.25);
\fill[white] (3,-3) circle (0.25);
\draw[thick] (0,0)..controls (1,-2) and (3,-2) .. (4,0);
\draw[thick] (2,-4)..controls (3,-2) and (7,-2) .. (8,-4);
\fill[white] (6.85,-3) circle (0.25);
\draw[thick] (6,-4)..controls (7,-2) and (9,-2) .. (10,-4);
\fill[black] (0,-4) circle (0.1);
\fill[black] (2,-4) circle (0.1);
\fill[black] (4,-4) circle (0.1);
\fill[black] (6,-4) circle (0.1);
\fill[black] (8,-4) circle (0.1);
\fill[black] (10,-4) circle (0.1);
\end{tikzpicture}
\end{center}
We can multiply such objects by vertical concatenation as before. The BMW algebra is generated by such objects, with an adequate notion of isotopy and adding some local relations. Note that closed links (possibly intertwined with the strands) can live inside the rectangular strip, due to the desired stability by concatenation. One important local relation is the Kauffman skein relation:
\begin{equation}\label{rel-skeinBMW}
 \begin{tikzpicture}[scale=0.25]
\draw[thick] (0,2)..controls +(0,-2) and +(0,+2) .. (4,-2);
\fill[white] (2,0) circle (0.4);
\draw[thick] (4,2)..controls +(0,-2) and +(0,+2) .. (0,-2);
\node at (6,0) {$-$};
\draw[thick] (12,2)..controls +(0,-2) and +(0,+2) .. (8,-2);
\fill[white] (10,0) circle (0.4);
\draw[thick] (8,2)..controls +(0,-2) and +(0,+2) .. (12,-2);
\node at (17,0) {$=-(q-q^{-1})$};
\node at (22,0) {$\Bigl($};
\draw[thick] (23.5,2) -- (23.5,-2);\draw[thick] (27.5,2) -- (27.5,-2);
\node at (29.5,0) {$-$};
\draw[thick] (31.5,2)..controls (32.5,0) and (34.5,0).. (35.5,2);
\draw[thick] (31.5,-2)..controls (32.5,0) and (34.5,0).. (35.5,-2);
\node at (37,0) {$\Bigr)$};
\end{tikzpicture}
\end{equation}
allowing in some sense to ``resolve'' the crossings, in a way similar to the Hecke algebra situation, though more involved. We will not give the precise topological definition of the BMW algebra (see for example \cite{Mor}), and we shall be happy with an algebraic description by generators and relations. 

The BMW algebra $BMW_n(a,q)$ is the quotient of the braid group algebra $\bC B_n$ by the relations:
\[
\begin{array}{ll}
e_i \si_i =  a e_i & \text{for $i=1,\dots, n-1$,}\\[0.4em]
e_i \si_{i+1}^{\pm1} e_i = a^{\mp1} e_{i} & \text{for $i=1,\dots, n-2$,}
\end{array}\]
where we have set $e_i = 1-\frac{\si_i -\si_i^{-1}}{q - q^{-1}}$. Here $q$ and $a$ are two non-zero complex numbers and $q^2\neq 1$. In the realisation as an algebra of tangles, the element $e_i$ is 
\begin{center}
 \begin{tikzpicture}[scale=0.3]
\node at (-0.5,0) {$e_i=$};
\node at (2,3) {$1$};\fill (2,2) circle (0.2);\fill (2,-2) circle (0.2);
\draw[thick] (2,2) -- (2,-2);
\node at (4,0) {$\dots$};
\draw[thick] (6,2) -- (6,-2);
\node at (6,3) {$i-1$};\fill (6,2) circle (0.2);\fill (6,-2) circle (0.2);
\node at (10,3) {$i$};\fill (10,2) circle (0.2);\fill (10,-2) circle (0.2);
\node at (14,3) {$i+1$};\fill (14,2) circle (0.2);\fill (14,-2) circle (0.2);
\draw[thick] (10,2)..controls (11,0) and (13,0) .. (14,2);
\fill[white] (12,0) circle (0.4);
\draw[thick] (10,-2)..controls (11,0) and (13,0) .. (14,-2);
\draw[thick] (18,2) -- (18,-2);
\node at (18,3) {$i+2$};\fill (18,2) circle (0.2);\fill (18,-2) circle (0.2);
\node at (20,0) {$\dots$};
\draw[thick] (22,2) -- (22,-2);\fill (22,2) circle (0.2);\fill (22,-2) circle (0.2);
\node at (22,3) {$n$};
\node at (38,0) {$i\in\{1,\dots,n-1\}$,};
\end{tikzpicture}
\end{center}
as can be seen from the local Kauffman relation. The first relation, when written only in terms of $\si_i$, is:
\[(\si_i-q)(\si_i+q^{-1})(\si_i-a)=0\ ,\]
so that we have a cubic characteristic equation for the generators. The other relation involves three different strands. It is clear from the algebraic presentation that if we set $e_i=0$ we recover the Hecke algebra $H_n(q)$.

\begin{rema}[Kauffman polynomial]
The chain of BMW algebras supports a Markov trace, different from the one on the chain of Hecke algebras, resulting in an invariant of links called the Kauffman polynomial \cite{Kau}.
\end{rema}

\paragraph{The YB equation.} For any $i\in\{1,\dots,n-1\}$, define the following function taking values in the algebra $BMW_n(q,a)$:
\begin{equation}\label{Baxt-BMW}
\si_i(u)=\si_i+(q-q^{-1})\frac{1}{u-1}+(q-q^{-1})\frac{1}{a^{-1}qu+1}e_i\ .
\end{equation}
It turns out that the braided YB equation is satisfied in the BMW algebra by these functions \cite{Jo89}:
\[\si_i(u)\si_{i+1}(uv)\si_i(v)=\si_{i+1}(v)\si_i(uv)\si_{i+1}(u)\ .\]

Following the same steps as for the Hecke algebra, it remains to discuss whether we can find some local representation of the BMW algebra, in order to obtain genuine matrix solutions of the YB equation. As for the Hecke algebra, it turns out to be possible. We do not give the details  and refer to \cite{Mu2,Tu} for explicit formulas. In both cases (Hecke and BMW), the existence of a local representation has a far-reaching significance in the context of Schur--Weyl dualities. We will discuss this later in Section \ref{sec-cent}.

\subsection{Baxterization}

The two preceding subsections provide examples of Baxterization formulas (inside an algebra). We will quickly review this notion, very natural in our discussion. The terminology is due to V. Jones \cite{Jo-bax}. Some references for Baxterization are \cite{BM,CGX,CFRV,CPdA2,IO,KMN,ZGB}.

Start with a local representation of the braid group given by an operator $\check R\in\text{End}(V\otimes V)$. Baxterizing this solution is, roughly speaking, to find a way to add the spectral parameters. More precisely, one looks for an expression $\check R(u)$ depending on a parameter $u$, and built out of the matrix $\check R$, such that the braided YB equation is satisfied:
\[\check R_1(u)\check R_{2}(uv)\check R_1(u)=\check R_{2}(v)\check R_1(uv)\check R_{2}(u)\ \ \ \ \ \text{on $V\otimes V\otimes V$}.\]
The expression for $\check R(u)$ is seen as a function of $u$ taking values in $\text{End}(V\otimes V)$. Typically it takes values in the subalgebra of $\text{End}(V\otimes V)$ generated by $\check R$. For example, if $\check R$ is diagonalisable, one can look at an expression of the form:
\[\check R(u)=\sum_{i} f_i(u)P_i\ ,\]
where $P_i$ are the projectors on eigenspaces of $\check R$, and one tries to find expressions for the functions $f_i(u)$ such that the YB equation is satisfied. usually, there is a way ro recover the operator $\check R$ by taking some limiting values of the parameters $u$.

\vskip .2cm
At the level of an algebra $A$, a Baxterization formula often refers to the following construction. Assume for simplicity that $A$ is generated by elements $\sigma_1,\dots,\sigma_{n-1}$ satisfying the braid relations. Then a Baxterization formula in the algebra $A$ is an explicit formula for a function $\sigma_i(u)$ with values in the algebra $A$ (or more precisely, 
in the subalgebra generated by $\sigma_i$), such that the braided YB equation is satisfied inside $A$:
\[\sigma_i(u)\sigma_{i+1}(uv)\sigma_i(v)=\sigma_{i+1}(v)\sigma_i(uv)\sigma_{i+1}(u)\ .\]
The expression of $\sigma_i(u)$ in terms of $\sigma_i$ should be the same for every $i$. 

If we have a Baxterization for an algebra $A$, then we should look for local representations of $A$, that is, representations of 
$A$ in $\text{End}(V^{\otimes n})$ for some vector space $V$ given in the form:
\[\sigma_i\mapsto \text{Id}^{\otimes i-1}\otimes \check R\otimes \text{Id}^{\otimes n-i-1}\ ,\]
where $\check R\in\text{End}(V\otimes V)$. Applying the representation, the function $\sigma_i(u)$ is sent to a function with values in $\text{End}(V^{\otimes n})$ satisfying the YB equation for matrices.

We have already seen three Baxterizations formulas, one for the Hecke algebra (\ref{Baxt-H}), one for the BMW algebra (\ref{Baxt-BMW}) and one for the symmetric group in Example \ref{ex-solYB}. The Baxterization formula for the symmetric group (more accurately, its group algebra) is:
\[\check R(u)=\si_i+\frac{1}{u}\,,\ \ \ \ \ \text{(where $\si_i=(i,i+1)$)}.\]

\section{Quantum Groups}\label{sec-QG}

We seek for a common ground to interpret both the Hecke algebra and the BMW algebra, and if possible explaining that both of them admit a Baxterization formula. To do so, we make a little detour through quantum groups. 

Quantum groups were introduced in the mid 80s by Drinfeld and Jimbo \cite{Dr-QG,Ji-QG}. Some constructions in special cases were already found in \cite{KR,Skl}. Some general references on quantum groups are \cite{CP,Dri,Kas,KRT,KS}.

\subsection{Some properties of quantum groups}

Let $\mathfrak{g}$ be a complex simple Lie algebra (for example, $\mathfrak{g}=sl_N$). To keep these notes reasonable in size, we will not give the definitions of quasi-triangular Hopf algebras of which the quantum groups are famous examples (we refer the interested reader to \cite[\S 4]{CP}, \cite[\S 3\&8]{Kas}, \cite[\S 2]{KRT} and \cite[\S 1]{KS}). We will only pick up from this whole theory the properties especially relevant for us and briefly discuss them. This section and the following are only meant to roughly indicate, not going into details, that there exists interesting algebraic objects (quantum groups) which provide us with representations of the braid group and even more, with solutions of the Yang--Baxter equation.

To get a first feeling of these objects, let us mention that the quantum group that we denote $U_q(\mathfrak{g})$ is an associative algebra\footnote{Regarding the terminology ``quantum group'', it may be satisfying enough to note that the quantum groups $U_q(\mathfrak{g})$ are deformations of the classical objects $U(\mathfrak{g})$ associated to Lie groups and moreover, during the deformation, something which was commutative (the coproduct) becomes non-commutative. More convincing explanations can be found in \cite{Dri}.}, which can be defined explicitly by generators and relations. As a vector space, it looks like the universal enveloping algebra $U(\mathfrak{g})$ of the Lie algebra $\mathfrak{g}$. However, the way of multiplying elements is different. It is a deformation of the multiplication in $U(\mathfrak{g})$, in the sense that there is a certain way to send the parameter $q$ to $1$ which recovers the algebra $U(\mathfrak{g})$. The example of $\mathfrak{g}=sl_2$ will be treated after the general discussion.

\vskip .2cm
To summarize, the main properties that we are going to discuss here are\footnote{The existence of a trivial representation and of contragredient (or dual) representations are important properties of quantum groups that we omit in our discussion.}:
\begin{itemize}
\item the fact that we can make tensor products of representations of $U_q(\mathfrak{g})$;
\item the fact that $V\otimes W$ and $W\otimes V$ are isomorphic as representations of $U_q(\mathfrak{g})$, and that the isomorphisms provide local representations of the braid group;
\item and finally, the fact that we can upgrade this picture, using affine quantum groups, to get solutions of the Yang--Baxter equation.
\end{itemize}
\begin{rema}
From general deformation theory, it is known that there is no non-trivial deformation of the algebra $U(\mathfrak{g})$ for a simple Lie algebra $\mathfrak{g}$ (see \cite{Kas}). It means that in a sense $U_q(\mathfrak{g})$ and $U(\mathfrak{g})$ are isomorphic as algebras (even though the explicit isomorphism is not obvious at all). At first, this remark can be somewhat disturbing for we realise that the new multiplication of $U_q(\mathfrak{g})$ is more or less equivalent to the usual one of $U(\mathfrak{g})$. However, one has to keep in mind that $U_q(\mathfrak{g})$ and $U(\mathfrak{g})$ are not isomorphic as Hopf algebras, meaning that the coproduct (see below) of $U_q(\mathfrak{g})$ is really different from the usual one of $U(\mathfrak{g})$. Roughly speaking, we could say that going from $U(\mathfrak{g})$ to $U_q(\mathfrak{g})$ does not change much the algebra structure, and in particular the representation theory, but it changes non-trivially the way to perform tensor products of representations.
\end{rema}

Let $V_1$ and $V_2$ be two representations of $U_q(\mathfrak{g})$. The vector space $V_1\otimes V_2$ can be made a representation of $U_q(\mathfrak{g})$. This is made possible by the existence of a so-called coproduct $\Delta$, which is a morphism of algebras:
\[\Delta\ :\ \ U_q(\mathfrak{g})\ \longrightarrow\ U_q(\mathfrak{g})\otimes U_q(\mathfrak{g})\ .\]
Note that the vector space $V_1\otimes V_2$ naturally carries a representation of $U_q(\mathfrak{g})\otimes U_q(\mathfrak{g})$. Then the way to construct a representation of $U_q(\mathfrak{g})$ on $V_1\otimes V_2$ is simply by precomposing with the coproduct.

More explicitly, denote $\rho_1$ the morphism $U_q(\mathfrak{g})\to \text{End}(V_1)$ and similarly for $\rho_2$. Then, explicitly, the representation of $U_q(\mathfrak{g})$ on $V_1\otimes V_2$ is defined as:
\[\begin{array}{ccccc}
U_q(\mathfrak{g}) & \stackrel{\Delta}{\longrightarrow} & U_q(\mathfrak{g})\otimes U_q(\mathfrak{g}) & \stackrel{\rho_1\otimes\rho_2}{\longrightarrow} & \text{End}(V_1\otimes V_2)\\[0.5em]
a & \mapsto & \Delta(a)=\sum a'\otimes a'' & \mapsto & \sum \rho_1(a')\otimes\rho_2(a'')
\end{array}\]
Here and in some places below, we use the useful Sweedler notation $\Delta(a)=\sum a'\otimes a''$, where the summation index is omitted. In rigorous mathematical notations, it would have been $\Delta(a)=\sum_{i=1}^N a'_i\otimes a''_i$ for some $N>0$ and some $a'_i,a''_i\in U_q(\mathfrak{g})$.
\begin{exam}\label{exam-coproduct}
For a group algebra $\bC G$, the map $\Delta$ from $\bC G$ to $\bC G\otimes \bC G$ defined by $\Delta(g)=g\otimes g$ for any $g\in G$ extends to a morphism of algebras and leads to the standard way of performing tensor products of representations of a group. 

For a Lie algebra $\mathfrak{g}$, the map $\Delta$ from $U(\mathfrak{g})$ to $U(\mathfrak{g})\otimes U(\mathfrak{g})$ defined by $\Delta(x)=x\otimes 1+1\otimes x$ for any $x\in \mathfrak{g}$ extends to a morphism of algebras and leads to the standard way of performing tensor products of representations of a Lie algebra.
\end{exam}
Not any morphism of algebras $U_q(\mathfrak{g})\to U_q(\mathfrak{g})\otimes U_q(\mathfrak{g})$ will work as a nice coproduct. One property which is usually required is the coassociativity property. In formulas, it reads 
$$(\Delta\otimes \text{Id})\circ\Delta=(\text{Id}\otimes\Delta )\circ\Delta\ .$$
This is an equality for applications from $U_q(\mathfrak{g})$ to $U_q(\mathfrak{g})^{\otimes 3}$. Its meaning is the following. At this point, there are two ways of performing the tensor product of three representations: $(V_1\otimes V_2)\otimes V_3$ and $V_1\otimes (V_2\otimes V_3)$ and the resulting representations on $V_1\otimes V_2\otimes V_3$ do not have to be isomorphic. The coassociativity property ensures that the two representations in fact simply coincide. Thanks to the coassociativity, we do not need to worry about putting parentheses in a tensor product, since every possible ways of performing it lead to the same result.

It might be worth emphasizing that, in general, for an arbitrary associative algebra, there is no natural coproduct and no natural way of performing tensor product of representations (for example, for the Hecke algebra $H_n(q)$). So the existence of a coproduct is a first remarkable property of quantum groups.

\vskip .2cm
Now, for our discussion, one crucial property of tensor products of representations of $U_q(\mathfrak{g})$ is the following: for any two representations $V,W$ of $U_q(\mathfrak{g})$, the representations $V\otimes W$ and $W\otimes V$ are isomorphic, that is, we have an invertible linear operator $\check R_{V,W}$ between $V\otimes W$ and $W\otimes V$ which commutes with the action of $U_q(\mathfrak{g})$:
\[ \check R_{V,W}\ :\ \ V\otimes W\ \stackrel{\sim}{\to}\ W\otimes V\ ,\ \ \ \ \ \ \ \text{isomorphism of $U_q(\mathfrak{g})$-representations.}\]
\begin{rema}
The permutation operator $\text{P}_{V,W}$ sending $v\otimes w$ to $w\otimes v$, is the natural isomorphism of vector spaces between $V\otimes W$ and $W\otimes V$. For $U(\mathfrak{g})$, the permutation provides the isomorphism of representations. This is not so for $U_q(\mathfrak{g})$ (that is, $\check R_{V,W}\neq \text{P}_{V,W}$) and this is one of the important points about quantum groups.
\end{rema}
The existence of the isomorphisms $\check R_{V,W}$ is not enough for our purpose, we need moreover the fact that these isomorphisms satisfy a compatibility condition when considering a tensor product of three representations:
\begin{center}
 \begin{tikzpicture}[scale=0.3]
\node at (0,0) {$V_1\otimes V_2\otimes V_3$};
\draw[thick][->] (4,0.1) -- (7,2.1);
\draw[thick][->] (4,-0.1) -- (7,-2.1);
\node at (11,2.1) {$V_2\otimes V_1\otimes V_3$};
\node at (11,-2.1) {$V_1\otimes V_3\otimes V_2$};
\draw[thick][->] (15,2.1) -- (18,2.1);
\draw[thick][->] (15,-2.1) -- (18,-2.1);
\node at (22,2.1) {$V_2\otimes V_3\otimes V_1$};
\node at (22,-2.1) {$V_3\otimes V_1\otimes V_2$};
\draw[thick][->] (26,2.1) -- (29,0.1);
\draw[thick][->] (26,-2.1) -- (29,-0.1);
\node at (33,0) {$V_3\otimes V_2\otimes V_1$};
\end{tikzpicture}
\end{center}
This diagram shows the two possible paths from $V_1\otimes V_2\otimes V_3$ to $V_3\otimes V_2\otimes V_1$ by applying the isomorphisms $\check R_{V_i,V_j}$. The compatibility condition is that these two paths coincide. This condition writes as an equality of linear operator from $V_1\otimes V_2\otimes V_3$ to $V_3\otimes V_2\otimes V_1$:
\begin{equation}\label{comp}
\bigl(\check R_{V_1,V_2}\otimes \text{Id}_{V_3}\bigr)\circ\bigl(\text{Id}_{V_2}\otimes \check R_{V_1,V_3}\bigr)\circ\bigl(\check R_{V_2,V_3}\otimes \text{Id}_{V_1}\bigr)=\bigl(\text{Id}_{V_1}\otimes \check R_{V_2,V_3}\bigr)\circ\bigl(\check R_{V_1,V_3}\otimes \text{Id}_{V_2}\bigr)\circ\bigl(\text{Id}_{V_3}\otimes \check R_{V_1,V_2}\bigr)\ .
\end{equation}
As we will see just below, this is at this step of the story that the braid group makes its appearance.

\subsection{Local representations of the braid group}\label{subsec-QGbr}

If we consider the situation where $V_1=V_2=V_3=V$ and we denote simply by $\check R$ the isomorphism $\check R_{V,V}$ which is an invertible element in $\text{End}(V\otimes V)$, then the compatibility condition above becomes simply the braid relation:
\[
\check R_1\check R_{2}\check R_1=\check R_{2}\check R_1\check R_{2}\ \quad \ \ \ \ \text{on $V\otimes V\otimes V$.} 
\]
The operator $\check R$ on $V\otimes V$ is called the $R$-matrix associated to the representation $V$ of the quantum group $U_q(\mathfrak{g})$.
\begin{conc}\label{conc-QGbr}
For any quantum group $U_q(\mathfrak{g})$ and any representation $V$ of $U_q(\mathfrak{g})$, we have a local representation of the braid group on $V^{\otimes n}$.
\end{conc}
Explicitly, that is to say that we have elements $\check R_1,\dots,\check R_{n-1}$ in $\text{End}(V^{\otimes n})$, which are constructed from the element $\check R\in\text{End}(V\otimes V)$ like this:
\[\rlap{$\overbrace{\phantom{V\otimes V}}^{\check R_1}$}V\otimes\overbrace{V\otimes V}^{\check R_2}\otimes\dots\overbrace{\ldots\otimes V}^{\check R_{n-1}}\]
and which satisfy the braid relations.

\paragraph{Universal $R$-matrix.} For an algebra $A$ admitting a coproduct $\Delta$, we denote $\Delta^{op}$ the composition of $\Delta$ with the permutation of $A\otimes A$ ($a\otimes b\mapsto b\otimes a$). We say that the coproduct is cocommutative when the images $\Delta(a)$ are invariant under the transposition of the two components, that is, when $\Delta=\Delta^{op}$. Note that the two classical examples of coproducts for groups and for Lie algebras (Example \ref{exam-coproduct}) are cocommutative. For an algebra with a cocommutative coproduct, it is easy to see that the permutation operator from $V\otimes W$ to $W\otimes V$ provides an isomorphism of representations. So in this case, the obtained local representation of the braid group is always the trivial one given by the permutation operator.

One of the main interest of quantum groups is that their coproduct is not cocommutative, and thus the obtained local representations of the braid group are not trivial. So regarding the construction of local representations of the braid group, the quantum groups $U_q(\mathfrak{g})$ are a considerable improvement compared to their classical analogue $U(\mathfrak{g})$.

\vskip .2cm
For a non-cocommutative coproduct on an algebra $A$, there is a natural condition  weakening the cocommutativity that ensures the isomorphisms of representations $V\otimes W$ and $W\otimes V$. This is the existence of an invertible element $\mathcal{R}$ in $A\otimes A$ such that:
\begin{equation}\label{univ-R}
\mathcal{R}\cdot\Delta(a)=\Delta^{op}(a)\cdot\mathcal{R}\ \ \ \ \forall a\in A\ .
\end{equation}
It follows immediately from this that the following element provides an isomorphism $\check R_{V,W}$ of representation between $V\otimes W$ and $W\otimes V$:
\[\check R_{V,W}=P_{V,W}\cdot \rho_{V\otimes W}(\mathcal{R})\,,\]
where $\rho_{V\otimes W}(\mathcal{R})$ is the image of the element $\mathcal{R}$ of $A\otimes A$ in the representation on $V\otimes W$; it belongs to $\text{End}(V\otimes W)$. And $P_{V,W}$ is the permutation operator from $V\otimes W$ to $W\otimes V$.

At this point, there is no reason for the compatibility condition (\ref{comp}) to be satisfied. A way to get it is obtained through the so-called quasitriangularity conditions, which read as:
\begin{equation}\label{quasiT}
(\Delta\otimes\text{Id}_A)(\mathcal{R})=\mathcal{R}_{13}\mathcal{R}_{23}\ \ \ \text{and}\ \ \ (\text{Id}_A\otimes\Delta)(\mathcal{R})=\mathcal{R}_{13}\mathcal{R}_{12}\ .
\end{equation}
These two conditions, which are equations in $A\otimes A\otimes A$, imply in particular the constant YB equation $\mathcal{R}_{12}\mathcal{R}_{13}\mathcal{R}_{23}=\mathcal{R}_{23}\mathcal{R}_{13}\mathcal{R}_{12}$ directly in $A\otimes A\otimes A$. From this, the compatibility condition for all representations as in (\ref{comp}) is ensured, and in particular the local representations of the braid group as in Conclusion \ref{conc-QGbr} are obtained.

\vskip .2cm
Such an element $\mathcal{R}$ satisfying (\ref{univ-R}) and (\ref{quasiT}) is called a universal $R$-matrix. It turns out that the quantum group $U_q(\mathfrak{g})$ admits a universal $R$-matrix satisfying all the required properties. We note that an appropriate definition of $U_q(\mathfrak{g})$ and $U_q(\mathfrak{g})\otimes U_q(\mathfrak{g})$ is needed, involving formal power series and suitable completions (see the example of $sl_2$ below).

\subsection{Solutions of the YB equation}\label{subsec-QGYB}

The quantum groups were originally designed for being useful in the study of the YB equation, but so far we have only discussed how they are related to local representations of the braid group. In short we must explain how to add the spectral parameters. This requires increasing by one (rather big) step the technical difficulties by going to the so-called ``affine'' Lie algebras. We are thus going to be even more sketchy than before. We note that the use of affine quantum groups for the YB equation was present since the origin of the theory of quantum groups \cite{Ji86}. References for this section are \cite{CP1,CP,EFK,FR,GRS}. The YB equation is now connected with many more mathematical structures than only quantum groups, see for example \cite{He} and references therein.

\vskip .2cm
The property we would like to discuss is the existence of solutions of the YB equation on representations of $U_q(\mathfrak{g})$. We will only try to convey the idea that this remarkable property comes from the existence of affine quantum groups with good properties. We organise the discussion in a series of steps towards the main conclusion.

\vskip .2cm
$\bullet$ The affine Lie algebra $\hat{\mathfrak{g}}$ is a central extension of the Lie algebra $\mathfrak{g}\otimes \bC[t,t^{-1}]$ of Laurent polynomials with coefficients in $\mathfrak{g}$ (we will forget the central extension in our discussion). If $V$ is a representation of $\mathfrak{g}$ then for any non-zero $a\in\bC$, there is a representation $V(a)$ of $\hat{\mathfrak{g}}$ called evaluation representation. As a representation of $\mathfrak{g}$, $V(a)$ is simply $V$, while the generator $t$ is evaluated to the number $a$. The new parameter $a$ in the fixed representation $V$ of $\mathfrak{g}$ will play the role of the spectral parameter in the YB equation. Only problem is we need to do the same thing for representations of $U_q(\mathfrak{g})$. This is where a quantum version of $\hat{\mathfrak{g}}$ is required.

\vskip .2cm
$\bullet$ The affine Lie algebra $\hat{\mathfrak{g}}$ is not a simple Lie algebra. However, it belongs to the family of Kac--Moody Lie algebra. In a few words, if we cast the algebraic definition of $\mathfrak{g}$ such that it depends only on its Cartan matrix, then we realise that we can use a similar definition for generalised Cartan matrices. This results in the Kac--Moody Lie algebras, including the affine Lie algebras. It turns out that the definition of the quantum group $U_q(\mathfrak{g})$ that we have avoided to give also only depends on the Cartan matrix. And thus it is perfectly generalisable to the generalised Cartan matrices leading to quantum groups associated to Kac--Moody Lie algebras. The point is that there is a definition of affine quantum groups $U_q(\hat{\mathfrak{g}})$ similar to the definition of quantum groups $U_q(\mathfrak{g})$ (the resulting algebras are quite more complicated though). 

As in the preceding discussion for $U_q(\mathfrak{g})$, there is an analogous story for tensor products of representations of $U_q(\hat{\mathfrak{g}})$. The algebra $U_q(\hat{\mathfrak{g}})$ still admits a coproduct, so we can make tensor products of representations, and there is in some completion a universal $R$-matrix allowing, roughly speaking, to reproduce the discussion before Formula (\ref{comp}). 

\vskip .2cm
$\bullet$ Generalising the natural embedding of $\mathfrak{g}$ in $\hat{\mathfrak{g}}$, one may regard $U_q(\mathfrak{g})$ as a subalgebra of $U_q(\hat{\mathfrak{g}})$. This allows to restrict representations of $U_q(\hat{\mathfrak{g}})$ to representations of $U_q(\mathfrak{g})$. Let us assume for now that for a representation $V$ of $U_q(\mathfrak{g})$, there exists an analogue of evaluation representations. Namely that we have representations $V(a)$ of $U_q(\hat{\mathfrak{g}})$, depending on a parameter $a$, such that the restriction to $U_q(\mathfrak{g})$ is the representation $V$ we started with.

For two ``evaluation'' representations $V(a)$ and $V(b)$ of $U_q(\hat{\mathfrak{g}})$, we can find for generic values of $a$ and $b$ an operator intertwining the representations $V(a)\otimes V(b)$ and $V(b)\otimes V(a)$. This operator, denoted $\check R(a,b)$, thus lives in $\text{End}(V\otimes V)$ since the vector space underlying both representations $V(a)$ and $V(b)$ is $V$. In fact, it turns out that it can be seen as a rational function in the parameters $a,b$:
\[\begin{array}{lcrcl}\check R\ &:&\ \bC^2 & \to & \text{End}(V\otimes V)\\[0.5em]
 & & (a,b) & \mapsto & \check R(a,b) \end{array}\ ,\]
The compatibility condition applied to $V(a)\otimes V(b)\otimes V(c)$ looks like this:
\begin{center}
 \begin{tikzpicture}[scale=0.3]
\node at (0,0) {$V(a)\otimes V(b)\otimes V(c)$};
\draw[thick][->] (6,0.1) -- (8,2.1);
\draw[thick][->] (6,-0.1) -- (8,-2.1);
\node at (14,2.1) {$V(b)\otimes V(a)\otimes V(c)$};
\node at (14,-2.1) {$V(a)\otimes V(c)\otimes V(b)$};
\draw[thick][->] (20,2.1) -- (22,2.1);
\draw[thick][->] (20,-2.1) -- (22,-2.1);
\node at (28,2.1) {$V(b)\otimes V(c)\otimes V(a)$};
\node at (28,-2.1) {$V(c)\otimes V(a)\otimes V(b)$};
\draw[thick][->] (34,2.1) -- (36,0.1);
\draw[thick][->] (34,-2.1) -- (36,-0.1);
\node at (42,0) {$V(c)\otimes V(b)\otimes V(a)$};
\end{tikzpicture}
\end{center}
For the function $\check R(a,b)$, this reads:
\[
\check R_{12}(b,c)\check R_{23}(a,c)\check R_{12}(a,b)=\check R_{23}(a,b)\check R_{12}(a,c)\check R_{23}(b,c)\ \quad \ \ \ \ \text{on $V\otimes V\otimes V$.} 
\]
This is the braided YB equation. Note that in the end the spectral parameters came from a family of representations of $U_q(\hat{\mathfrak{g}})$ extending the given representation of $U_q(\mathfrak{g})$. So it would not be too far from truth to say that the YB equation is the braid relation applied to some representations of affine quantum groups. 

At last we can draw our second conclusion regarding the uses of quantum groups in our context. 
\begin{conc}\label{conc-QGYB}
For a quantum group $U_q(\mathfrak{g})$ and some of its representations $V$, we have a solution of the YB equation on $V\otimes V$.
\end{conc}
Let us admit for the sake of simplicity that the function $\check R(a,b)$ depends only on the ratio of the two parameters, and therefore is equivalent to a function of one parameter $\check R(u)$. The conclusion says that we have functions $\check R_1(u),\dots,\check R_{n-1}(u)$ taking values in $\text{End}(V^{\otimes n})$, which are constructed from the function $\check R(u)$ like this:
\[\rlap{$\overbrace{\phantom{V\otimes V}}^{\check R_1(u)}$}V\otimes\overbrace{V\otimes V}^{...}\otimes\dots\overbrace{\ldots\otimes V}^{\check R_{n-1}(u)}\]
and which satisfy the braided YB equation:
\[\check R_i(u)\check R_{i+1}(uv)\check R_i(v)=\check R_{i+1}(v)\check R_i(uv)\check R_{i+1}(u)\ .\]

\paragraph{On the validity of Conclusion \ref{conc-QGYB}.} In the above statement we had to reduce the generality by writing ``some of its representations''. This is because we have assumed that for a representation $V$ of $U_q(\mathfrak{g})$, we had a family of ``evaluation'' representations $V(a)$ of $U_q(\hat{\mathfrak{g}})$. This assumption is not always satisfied and we will discuss briefly its validity.

First, for the simplest situation $\mathfrak{g}=sl_N$, this assumption is valid for any representation $V$ of $U_q(\mathfrak{g})$. Indeed, for $\mathfrak{g}=sl_N$, there is a morphism from $U_q(\hat{\mathfrak{g}})$ to $U_q(\mathfrak{g})$ depending on a parameter $a$, which is the analogue of the evaluation morphism from $\hat{\mathfrak{g}}$ to $\mathfrak{g}$. This allows immediately to construct the ``evaluation'' representations $V(a)$ of $U_q(\hat{\mathfrak{g}})$ from a representation $V$ of $U_q(\mathfrak{g})$.

For other simple Lie algebras $\mathfrak{g}$, such an analogue of evaluation morphism does not exist. However, what is true in general is that there is a family of automorphisms of $U_q(\hat{\mathfrak{g}})$ depending on a parameter $a$. So assume that a representation $V$ of $U_q(\mathfrak{g})$ can be extended to a representation of $U_q(\hat{\mathfrak{g}})$. Then by twisting this representation by the automorphisms, one get a family of representations $V(a)$ of $U_q(\hat{\mathfrak{g}})$ extending $V$ (they do not come anymore from a morphism $U_q(\hat{\mathfrak{g}})\to U_q(\mathfrak{g})$ but they still serve for finding solutions of YB equation on $V$).

This leaves us with the following conclusion: we have a solution of the YB equation on representations $V$ of $U_q(\mathfrak{g})$ which can be extended to $U_q(\hat{\mathfrak{g}})$. To insist once more, for $\mathfrak{g}=sl_N$, all representations can be extended but for arbitrary $\mathfrak{g}$ this is not true. This leads to the natural concept of ``minimal affinization'' of a representation \cite{CP2}. Note that reversing the point of view, one could also say: we have a solution of the YB equation on representations $V$ of $U_q(\mathfrak{g})$ which can be obtained from restrictions of representations of $U_q(\hat{\mathfrak{g}})$.

\begin{rema}
Example \ref{ex-solYBq} is associated to $U_q(sl_2)$ and the fundamental representation. However, before that, we had the example \ref{ex-solYB}, the Yang solution, which is a certain limit of the one from $U_q(sl_2)$. In fact, this example also comes from an algebra, sharing many properties with $U_q(\hat{\mathfrak{g}})$, which is called the Yangian $Y(\mathfrak{g})$ (see \cite{CP}). This algebra, with no parameter $q$, can be seen a certain limit of $U_q(\hat{\mathfrak{g}})$ and also provides solutions of the YB equation along the same lines as outlined above. 
\end{rema}

\subsection{Example of $U_q(sl_2)$}\label{subsec-sl2}

This is about time for an explicit example, so we will take the simplest situation $\mathfrak{g}=sl_2$, see \cite{Kas}. As a vector space, $sl_2$ is generated by:
\[h=\left(\begin{array}{cc} 1 & 0 \\ 0 & -1 \end{array}\right)\,,\ \ \ x=\left(\begin{array}{cc} 0 & 1 \\ 0 & 0 \end{array}\right)\,,\ \ \ y=\left(\begin{array}{cc} 0 & 0 \\ 1 & 0 \end{array}\right)\,.\]
Its Lie algebra structure is given by:
\begin{equation}\label{rel-sl2}
[h,x]=2x\,,\ \ [h,y]=-2y\,,\ \ [x,y]=h\ .
\end{equation}
Thus, as an algebra, $U(sl_2)$ is generated by elements $h,x,y$ with relations (\ref{rel-sl2}) as defining relations. Note that we abuse notations by keeping the same names $x,y,h$ for the generators of $U(sl_2)$. In $U(sl_2)$, the bracket $[a,b]$ simply means $ab-ba$.

\paragraph{The algebra $U_q(sl_2)$.} 
The quantum group $U_q(sl_2)$ is often defined as the algebra generated by $X,Y,K$ with $K$ required to be invertible, and with defining relations:
\begin{equation}\label{def-sl2-q1}
KXK^{-1}=q^2 X\,,\ \ \ KYK^{-1}=q^{-2}Y\,,\ \ \ XY-YX=\frac{K-K^{-1}}{q-q^{-1}}\ .
\end{equation}
In this version, it is not so clear how to relate it to $U(sl_2)$ (at this point, the limit $q\to 1$ does not really make sense) and it would also require some completions to be able to find a universal $R$-matrix.

\vskip .2cm
The alternative standard definition of $U_q(sl_2)$, the one which we are going to consider from now is the following one. As a vector space, we define $U_q(sl_2)$ to be
\[U(sl_2)[[\alpha]]=\{c_0+c_1\alpha+\dots\,,\ c_i\in U(sl_2)\}\,,\]
the vector space of formal power series in $\alpha$ with coefficients in $U(sl_2)$. We abuse again notation and keep the names $x,y,h$ for the generators of $U(sl_2)$ but we insist that we mean here only the vector space $U(sl_2)$; the multiplication will be different. The multiplication of elements of $U_q(sl_2)$ is the usual multiplication of formal series, together with the following (new) defining relations between the elements $x,y,h$:
\begin{equation}\label{def-sl2-q2}
hx-xh=2x\,,\ \ hy-yh=-2y\,,\ \ \ xy-yx=\frac{e^{\alpha h}-e^{-\alpha
h}}{e^{\alpha}-e^{-\alpha}}\ .
\end{equation}
One has to notice that the right hand side of the last relation is indeed a power series in $\alpha$. We should insist that as a vector space, $U(sl_2)[[\alpha]]$ is not the same as $U(sl_2)$ with
coefficients in $\mathbb{C}[[\alpha]]$. This is because $U(sl_2)$ is not finite dimensional. For
example, $e^{\alpha h}$ is in $U(sl_2)[[\alpha]]$.

One can see the connections with the previous definition of $U_q(sl_2)$ as follows: if one sets $q=e^{\alpha}$ and $K=e^{\alpha h}$ then Relations (\ref{def-sl2-q1}) are satisfied. One advantage of the definition of $U_{q}(sl_2)$ involving power series in $\alpha$ is that the connection with $U(sl_2)$ is quite transparent. The limit $\alpha=0$ is well-defined and one recovers the algebra $U(sl_2)$ as can be seen immediately in the defining relations.

\paragraph{Coproduct and $R$-matrix for $U_q(sl_2)$.} 
We start by defining the tensor product $U_q(sl_2)\hat{\otimes} U_q(sl_2)$ as follows (and similarly for $U_{\alpha}(sl_2)^{\hat{\otimes} L}$ for any $L$). As a vector space, this is:
\[\bigl(U(sl_2)\otimes U(sl_2)\bigr)[[\alpha]]\ .\]
The power series are multiplied as usual, and the elements $x,y,h$ still satisfy the relations (\ref{def-sl2-q2}) in
each factor of the tensor product; moreover elements in different factors commute.

Then, the coproduct is defined on the generators by:
\begin{equation}\label{coprod-sl2}\Delta(x)=x\otimes e^{-\alpha h/2}+e^{\alpha h/2}\otimes x\,,\ \ \ \Delta(y)=y\otimes
e^{-\alpha h/2}+e^{\alpha h/2}\otimes y\,,\ \ \ \Delta(h)=h\otimes 1+1\otimes h\,,
\end{equation}
and extended to $U_{q}(sl_2)$ by $\Delta(\sum c_i\alpha^i)=\sum \Delta(c_i)\alpha^i$. One can check that $\Delta$ extends to an algebra homomorphism from $U_{q}(sl_2)$ to $U_q(sl_2)\hat{\otimes} U_q(sl_2)$. It is immediate to see that the limit $\alpha=0$ gives back the usual coproduct on $U(sl_2)$.

Note that $U_q(sl_2)\hat{\otimes} U_q(sl_2)$ is not the same vector space as $U_q(sl_2)\otimes U_q(sl_2)$. For example $e^{\alpha (h\otimes h)}$ is in the former but not in the latter. At this point, it may seem superfluous to consider the completed tensor product $\hat{\otimes}$, since the coproduct takes values in the usual tensor product. However, the completed tensor product is relevant because it turns out that there is an universal $R$-matrix $\mathcal{R}$ with the required properties in $U_q(sl_2)\hat{\otimes} U_q(sl_2)$. An explicit formula is:
\begin{equation}\label{R-sl2}\mathcal{R}=e^{\alpha(h\otimes h)/2}\sum_{n\geq 0}\frac{(q-q^{-1})^n}{[n]_q!}q^{n(n-1)/2}\Bigl(e^{-\alpha h/2}y\otimes x e^{\alpha h/2}\Bigr)^n\,,
\end{equation}
where we have set $q=e^{\alpha h}$, $[n]_q!=[2]_q\cdot [3]_q\dots [n]_q$ where $[k]_q=\frac{q^k-q^{-k}}{q-q^{-1}}$. The term in the sum being a multiple of $\alpha^n$, it follows that the element $\mathcal{R}$ is a well-defined element of $U_q(sl_2)\hat{\otimes} U_q(sl_2)$.

\paragraph{Example of representations of $U_q(sl_2)$.} By checking directly the defining relations, we find that:
\[h\mapsto\left(\begin{array}{cc} 1 & 0 \\ 0 & -1 \end{array}\right)\,,\ \ \ x\mapsto\left(\begin{array}{cc} 0 & 1 \\ 0 & 0 \end{array}\right)\,,\ \ \ y\mapsto\left(\begin{array}{cc} 0 & 0 \\ 1 & 0 \end{array}\right)\,,\]
defines a representation of $U_q(sl_2)$ on a vector space $V$ of dimension 2. It is called the vector representation (it looks exactly the same as the natural representation of $sl_2$, but look at the next example). To calculate the $R$-matrix for this representation, one notes that the images of $x$ and $y$ are nilpotent of index 2, and thus it is enough to take only $n=0,1$ in the sum in (\ref{R-sl2}). Say $(v_1,v_2)$ was the basis of $V$ used to give the matrices above, then in the basis $(v_i\otimes v_j)_{i,j=1,2}$ ordered lexicographically, we easily find the image of $\mathcal{R}$ and in turn the $R$-matrix:
\[\rho_{V\otimes V}(\mathcal{R})=\left(\begin{array}{cccc}
e^{\alpha/2} & \cdot & \cdot & \cdot \\
\cdot  & e^{-\alpha/2} & \cdot & \cdot \\
\cdot  & e^{-\alpha/2}(q-q^{-1}) & e^{-\alpha/2} & \cdot \\
\cdot & \cdot & \cdot & e^{\alpha/2}
\end{array}\right)\ \ \ \ \Rightarrow\ \ \ \check R_{V,V}=e^{-\alpha/2}\left(\begin{array}{cccc}
q & \cdot & \cdot & \cdot \\
\cdot  & q-q^{-1} & 1& \cdot \\
\cdot  & 1 & 0& \cdot \\
\cdot & \cdot & \cdot & q
\end{array}\right)\,,\]
where we recall that the matrix $\check R$ is obtained by multiplying by the permutation operator. Up to a global factor, we recover the solution in Example \ref{ex-solYBq} coming from the Hecke algebra for a vector space of dimension 2.

\vskip .2cm
As a first instance of a procedure whose generalisation we are going to discuss later, we note that the $R$-matrix $\check R_{V,V}$ commutes with the action of $U_q(sl_2)$ on $V\otimes V$. Namely, it commutes with the image of $\Delta(h),\Delta(x),\Delta(y)$ in $\text{End}(V\otimes V)$. This can be easily checked here with a direct calculation, but as we will see, it follows by construction of $R$-matrices. 

Therefore, the action of $U_q(sl_2)$ leaves invariant the eigenspaces of $\check R_{V,V}$. A basis of the eigenspace for the eigenvalue $q$ is $(v_1\otimes v_1,\,e^{-\alpha/2}v_1\otimes v_2+e^{\alpha/2}v_2\otimes v_1,\,v_2\otimes v_2)$. This is the deformation of the symmetric square of $V$. On this subspace, we find the following representation of $U_q(sl_2)$:
\[h\mapsto\left(\begin{array}{ccc} 2 & 0 & 0 \\ 0 & 0 & 0\\0 & 0 & -2 \end{array}\right)\,,\ \ \ x\mapsto\left(\begin{array}{ccc} 0 & q+q^{-1} & 0 \\ 0 & 0 & 1\\0 & 0 & 0 \end{array}\right)\,,\ \ \ y\mapsto\left(\begin{array}{ccc} 0 & 0 & 0 \\ 1 & 0 & 0\\0 & q+q^{-1} & 0 \end{array}\right)\,.\]
This is the deformation for $U_q(sl_2)$ of the three-dimensional representation of $sl_2$. The $R$-matrix associated to this representation can also be calculated directly from the formula (\ref{R-sl2}), we leave the details for the reader.

We note that the $R$-matrix for this representation does not come anymore from the Hecke algebra, and one goal of these notes is to discuss the algebras controlling this solution and the analogues for higher symmetric powers.

\newpage
\section{Centralisers of Tensor Representations of Quantum Groups}\label{sec-cent}

\subsection{Properties of centralisers}

\paragraph{Generalities.} Let $\rho\ :\ A\to\text{End}(E)$ be a representation of an algebra $A$ on a vector space $E$. The centraliser of the representation is a classical object in linear algebra. In words, it is the commutant in $\text{End}(E)$ of the family of endomorphisms $\{\rho(a)\}_{a\in A}$\,. More formally, we define:
\[\text{End}_A(E):=\{x\in\text{End}(E)\ |\ \rho(a)x=x\rho(a)\,,\ \forall a\in A\}\ .\]
The meaning of the centraliser in representation theory is quickly understood when one considers a semisimple representation $E$, that is a representation $E$ which decomposes as a direct sum of irreducible representations:
\[E=\bigoplus_i E_i^{\oplus m_i}\,,\]
where $m_i\in\mathbb{Z}_{>0}$ are the multiplicities of the irreducible representations $E_i$ in $E$. Using Schur Lemma, the centraliser $\text{End}_A(E)$ is described as follows. First it leaves invariant the subspaces $E_i^{\oplus m_i}$. Besides, for a given $i$,  choose a basis for the representation $E_i$ and denote by $\rho_i(a)$ the matrices of elements of $A$ acting in this basis of $E_i$. Form a basis of $E_i^{\oplus m_i}$ by concatenating $m_i$ times this same basis. Then on $E_i^{\oplus m_i}$, the operators of $A$ and of the centralisers look as follows:
\[A\ni a\mapsto\left(\begin{array}{cccc}\rho_i(a) & 0 & \dots & 0\\
0 & \rho_i(a) & \ddots & \vdots\\
\vdots & \ddots & \ddots & 0\\
0 & \ldots & 0 & \rho_i(a)\\\end{array}\right)\,,\ \ \ \ \ \ \text{End}_A(E)\ni x\mapsto\left(\begin{array}{cccc}x_{11}\text{Id} & x_{12}\text{Id} & \dots & x_{1m_i}\text{Id}\\
x_{21}\text{Id} & x_{22}\text{Id} & \ddots & \vdots\\
\vdots & \ddots & \ddots & \vdots\\
x_{m_i1}\text{Id} & \ldots & \ldots & x_{m_im_i}\text{Id}\\\end{array}\right)\,,\]
as block-matrices with $m_i\times m_i$ blocks. Thus we see that the centraliser $\text{End}_A(E)$ is isomorphic to a direct sum of matrix algebras:
\[\text{End}_A(E)\cong \bigoplus_i \text{Mat}_{m_i}(\bC)\ ,\]
the sizes of which correspond to the multiplicities.

One can go further and note that $E$ can also be seen as a representation of the centraliser $\text{End}_A(E)$, and that for this representation, the centraliser is the image of $A$ in $\text{End}(E)$. This is an instance of a double centralising theorem (note that the assumption that the representation of $A$ we started with is semisimple is important). Finally, as the images commute, $E$ can also be seen as a representation of the algebra $A\otimes \text{End}_A(E)$ (the multiplication is performed independently in each factor).

To summarise, we have the following decompositions of $E$ as a representation, respectively, of $A$, of $\text{End}_A(E)$ and finally of $A\otimes \text{End}_A(E)$:
\begin{equation}\label{dec-gen}
E=\bigoplus_i E_i^{\oplus m_i}\,,\quad\ \ \ E=\bigoplus_i \bigl(\bC^{m_i}\bigr)^{\oplus \dim(E_i)}\,,\quad\ \ \ E=\bigoplus_i E_i\otimes \bC^{m_i}\,.
\end{equation}

\paragraph{Centralisers of $U_q(\mathfrak{g})$.} Now, we are mainly interested here in the case of $A=U_q(\mathfrak{g})$ a quantum group as in Section \ref{sec-QG} and $E$ is a tensor product of representations. More precisely, we take a representation $V$ of $U_q(\mathfrak{g})$ and we form the tensor product $V^{\otimes n}$. We have explained the fundamental property of quantum groups that this tensor product $V^{\otimes n}$ is also a representation of $U_q(\mathfrak{g})$. The objects of main interest to us are the centralisers of such representations:
\[\text{End}_{U_q(\mathfrak{g})}(V^{\otimes n})\ .\]

Now it is time to draw on our previous discussion on quantum groups and fulfill the objective set up in the introduction. Recall from Section \ref{subsec-QGbr} (Conclusion \ref{conc-QGbr}) that we have elements $\check R_1,\dots,\check R_{n-1}$ of $\text{End}(V^{\otimes n})$, satisfying the braid relations. 

Moreover recall from Section \ref{subsec-QGYB} that, assuming the representation $V$ is such that Conclusion \ref{conc-QGYB} is valid,  we have functions $\check R_1(u),\dots,\check R_{n-1}(u)$ taking values $\text{End}(V^{\otimes n})$ satisfying the braided YB equation. The quantum groups are so cleverly designed that the following holds.
\begin{conc}\label{conc-finQG}$\ $
\begin{enumerate}
\item The elements $\check R_1,\dots,\check R_{n-1}$ belong to the centraliser $\text{End}_{U_q(\mathfrak{g})}(V^{\otimes n})$.

\item The functions $\check R_1(u),\dots,\check R_{n-1}(u)$ take values in the centraliser $\text{End}_{U_q(\mathfrak{g})}(V^{\otimes n})$.
\end{enumerate}
\end{conc}
\begin{proof}[Proof of Conclusion \ref{conc-finQG}] The facts stated above being clearly fundamental in our discussion, this is a good point to give one proof in these notes. It may help to see a little bit better how the subtle machinery of quantum groups works. So let us prove item 1.

Denote $\rho\ :\ U_q(\mathfrak{g})\to\text{End}(V)$ the representation. Let $a$ be an arbitrary element of $U_q(\mathfrak{g})$ and denote $\Delta(a)=\sum a'\otimes a''$ its coproduct. 

We discussed in Section \ref{sec-QG} the coassociativity property of the coproduct and how it allows to make tensor products without worrying about parentheses. Let us see more explicitly how this happens. Recall that the representation $\rho^{(2)}\ :\ U_q(\mathfrak{g})\to\text{End}(V\otimes V)$ is constructed using the coproduct and is given by 
$$\rho^{(2)}(a)=(\rho\otimes \rho)\bigl(\Delta(a)\bigr)=\sum\rho(a')\otimes \rho(a'')\ .$$
To construct a representation on $V^{\otimes n}$, one should first put parentheses in $V^{\otimes n}$ and then perform a sequence of tensor products of two spaces. For $n=3$, we have two possibilities:
\[V\otimes (V\otimes V)\ \ \ \ \text{or}\ \ \ \ (V\otimes V)\otimes V\ .\]
These two possibilities correspond to representations on $V^{\otimes 3}$ given by:
\[\text{applying $\rho\otimes \rho\otimes \rho$ on:}\quad\ \ \bigl((\text{Id}\otimes \Delta)\circ\Delta\bigr)(a)\quad\ \ \ \text{or}\quad\ \ \ \bigl((\Delta\otimes\text{Id})\circ\Delta\bigr)(a)\ .\]
The coassociativity condition $(\text{Id}\otimes \Delta)\circ\Delta=(\Delta\otimes\text{Id})\circ\Delta$ ensures that these two possibilities give the same representation. Now, similarly, any way of putting parentheses on $V^{\otimes n}$ also corresponds to a certain composition of coproducts. For example, if $n=5$, here is one possibility:
\[V\ \ \leadsto\ \ V\otimes V\ \ \leadsto\ \ (V\otimes V)\otimes V\ \ \leadsto\ \ (V\otimes (V\otimes V))\otimes V\ \ \leadsto\ \ (V\otimes (V\otimes V))\otimes (V\otimes V)\,,\]
corresponding to the following composition: $(\text{Id}\otimes\text{Id}\otimes\text{Id}\otimes\Delta)\circ(\text{Id}\otimes\Delta\otimes \text{Id})\circ(\Delta\otimes \text{Id})\circ\Delta$. So to check that we really obtain always the same representation on $V^{\otimes n}$, we first consider one way of doing it, for example by defining $\Delta^{(2)}=\Delta$ and for $n>2$:
\[\Delta^{(n)}:=(\Delta\otimes \text{Id}\otimes\dots\otimes\text{Id})\circ\Delta^{(n-1)}\,,\]
which corresponds to the parenthesized product $\Bigl(\!\!...\!\bigl((V\otimes V)\otimes V\bigr)\otimes\dots\otimes V\Bigr)$. And then it is easy to check by induction on $n$ that it coincides with any other possibilities:
\[\Delta^{(n)}=(\text{Id}\otimes \dots\otimes\text{Id}\otimes\Delta\otimes \text{Id}\otimes \dots\otimes\text{Id})\circ\Delta^{(n-1)}\,,\ \ \ \ \ \text{for any position of $\Delta$ among the Id's}\ .\]
So finally there is only one representation on $\rho^{(n)}:U_q(\mathfrak{g})\to\text{End}(V^{\otimes n})$, and it is given by the following recursive formula for any choice of $i$:
\[\rho^{(n)}(a)=\sum\rho^{(n-i)}(a')\otimes \rho^{(i)}(a'')\ \ \ \ \ \text{for $i\in\{1,\dots,n-1\}$.}\]
\underline{$\bullet$ $n=2$:} After all this preparation, we are ready for the proof at last. Consider the element $\mathcal{R}$ in $U_q(\mathfrak{g})\otimes U_q(\mathfrak{g})$ satisfying:
\[\mathcal{R}\Delta(a)=\Delta^{op}(a)\mathcal{R}\,,\]
and recall that the element $\check R\in\text{End}(V\otimes V)$ comes from $\mathcal{R}$ in the sense: $\check R=P \rho\otimes\rho(\mathcal{R})$ (where $P$ is the permutation operator of $V\otimes V$). So applying the representation $\rho\otimes\rho$ to the equality above, we find:
\[P\check R\,\Bigl(\sum\rho(a')\otimes\rho(a'')\Bigr)=\Bigl(\sum\rho(a'')\otimes\rho(a')\Bigr)\,P\check R\ .\]
Now it is clear that, moving the permutation operator to the left and removing it, we have:
\[\check R\rho^{(2)}(a)=\rho^{(2)}(a)\check R\ ,\]
which is the statement that $\check R$ belongs to $\text{End}_{U_q(\mathfrak{g})}(V\otimes V)$.

\noindent \underline{$\bullet$ $n>2$:} By induction on $n$, using respectively the formulas:
\[\rho^{(n)}(a)=\sum\rho^{(n-1)}(a')\otimes \rho(a'')\ \ \ \ \text{and}\ \ \ \ \ \rho^{(n)}(a)=\sum\rho(a')\otimes \rho^{(n-1)}(a'')\,,\]
we obtain first that $\check R_1,\dots,\check R_{n-2}$ commutes with $\rho^{(n)}(a)$ and second that $\check R_{n-1}$ commutes with $\rho^{(n)}(a)$ as well.

For item 2, we can follow the same reasoning for $\check R(u)$ in the suitable tensor product of ``evaluation'' representations $V(a)$ of the affine quantum group $U_q(\hat{\mathfrak{g}})$. Restricted to $U_q(\mathfrak{g})$, the dependance on the ``evaluation'' parameter disappears and the tensor product is simply $V^{\otimes n}$. We skip the details here.
\end{proof}

\begin{exam}\label{ex-slN-cl}
Let us discuss a classical example of centralisers. Consider the universal enveloping algebra $U(sl_N)$ (so this is the classical limit of the preceding situation). The Lie algebra $sl_N$ is the Lie algebra of $N\times N$ matrices with zero trace. As such there is a natural representation of dimension $N$ (send a matrix to itself). Call it $V$ and construct the representation $V^{\otimes n}$. Recall that for a Lie algebra, the action of an element $g$ on $V^{\otimes n}$ is the sum of the action on each factor:
\[g\otimes\text{Id}\otimes \dots\otimes \text{Id}+\text{Id}\otimes g\otimes \text{Id}\otimes\dots\otimes \text{Id}+\dots+\text{Id}\otimes\dots \otimes \text{Id}\otimes g\ .\]
Now it is clear that any permutation of the factors in $V^{\otimes n}$ commutes with the action of $U(sl_N)$. In particular, take $P_i$ to be the transposition of positions $i$ and $i+1$, then we have that $P_1,\dots,P_{n-1}$ belong to $\text{End}_{U(sl_N)}(V^{\otimes n})$. These operators satisfy the braid relations, so we just described the analogue of item 1 of Conclusion \ref{conc-finQG} in this case. Concerning item 2, the solutions of the braided YB equation (with additive spectral parameters) are given by the Yang formula, and they obviously take values in $\text{End}_{U(sl_N)}(V^{\otimes n})$:
\[\check R_i=P_i+\frac{Id_{V^{\otimes n}}}{u}\ .\]
Moreover, in this case, we can obtain a complete description of the centraliser $\text{End}_{U(sl_N)}(V^{\otimes n})$ since we can show that it is actually generated by the permutations (so it is a quotient of the algebra of the symmetric group). This is the classical Schur--Weyl duality, which we will describe for $U_q(sl_N)$ in the next subsection. 
\end{exam}

\subsection{Schur--Weyl duality}

At this point, we should be convinced that the centralisers $\text{End}_{U_q(\mathfrak{g})}(V^{\otimes n})$ are the algebras we were looking for, in the sense that they satisfy the properties discussed in the introduction: they contain a quotient of the braid group algebra and they contain solutions of the YB equation. About a representation on a tensor space $V^{\otimes n}$, well, they are defined as subalgebras of $\text{End}(V^{\otimes n})$, so they admit one by definition.

However, we are far from an explicit algebraic (or any other) description of $\text{End}_{U_q(\mathfrak{g})}(V^{\otimes n})$ for any $\mathfrak{g}$ and $V$. So the conclusion above is more a general statement indicating that this should be very interesting for us to study these algebras. Moreover, the claimed objective that we were going to ``explain'' why the Hecke (and BMW) algebras fit so nicely in our picture is still not achieved. This is what we are going to do now.

\paragraph{Quantum group $U_q(sl_N)$.} So we take $\mathfrak{g}=sl_N$ for some $N>0$ and, for the representation $V$, we take the analogue for $U_q(sl_N)$ of the natural representation of $sl_N$. From the general representation theory of quantum groups, we know that any (finite-dimensional irreducible) representation of $\mathfrak{g}$ can be deformed to a representation of $U_q(\mathfrak{g})$. However, let us be more explicit than that.

First, the quantum group $U_q(sl_N)$ should be described a bit more. We build on our previous description of $U_q(sl_2)$ in Section \ref{subsec-sl2} involving power series in $\alpha$ and so on, so we do not repeat our discussion ($q$ is still $e^{\alpha}$). Natural generators of $sl_N$ are the elements $h_i,x_i,y_i$, with $i=1,\dots,N-1$. The triplet of elements $h_i,x_i,y_i$ corresponds to the following matrices of $sl_N$: take the matrices $\left(\begin{array}{cc} 1 & 0 \\ 0 & -1 \end{array}\right)$, $\left(\begin{array}{cc} 0 & 1 \\ 0 & 0 \end{array}\right)$, $\left(\begin{array}{cc} 0 & 0 \\ 1 & 0 \end{array}\right)$ as for $sl_2$ and plug them in lines and columns $i$ and $i+1$.

The multiplication is deformed as follows. First each triplet $h_i,x_i,y_i$ satisfies the relations of $U_q(sl_2)$ in (\ref{def-sl2-q2}), replacing $x$ to $x_i$, $y$ to $y_i$ and $h$ to $h_i$. Then we keep the following usual relations of matrices:
\[[h_i,x_{i\pm1}]=-x_{i\pm1}\,,\ \ \ \ \ [h_i,y_{i\pm1}]=y_{i\pm1}\,,\ \ \ \ \ [h_i,x_{j}]=[h_i,y_{j}]=0\ \ (|i-j|>1)\,,\]
\[[h_a,h_b]=0\ \ (\forall a,b)\,,\ \ \ \ \ \ [x_i,y_j]=0\ \ (i\neq j)\ .\]
And finally we add the deformation of the so-called Serre relations (which is really the new ingredient compared to $sl_2$):
\[x_jx_i^2-(q+q^{-1})x_ix_jx_i+x_i^2x_{j}=y_jy_i^2-(q+q^{-1})y_iy_jy_i+y_i^2y_{j}=0\ \ \ \ \text{ if $|i-j|=1$\,.}\]
For each triplet $h_i,x_i,y_i$, the coproduct is given by the same formulas as for $U_q(sl_2)$ in (\ref{coprod-sl2}). That is it, we have defined $U_q(sl_N)$ and how to make tensor products of its representations.

As for $sl_2$, the standard representation $V$ of $U_q(sl_N)$ is simply given by assigning to each $h_i,x_i,y_i$ its natural $N\times N$ matrix, as if we were dealing with $sl_N$. One can be surprised maybe, but all defining relations of $U_q(sl_N)$ are indeed satisfied. So again, and especially in the standard representation, we are not changing much the algebra structure, but the ``quantum'' novelty of $U_q(sl_N)$ compared to classical $U(sl_N)$ is that we have a different way to construct a representation on $V^{\otimes n}$. Comparing to Example \ref{ex-slN-cl}, now the permutations do not belong to the centraliser $\text{End}_{U_q(sl_N)}(V^{\otimes n})$. This is where we meet again, at last, with the Hecke algebra.

\paragraph{Hecke algebra $H_n(q)$ and Jimbo--Schur--Weyl duality.} Informations on the classical Schur--Weyl duality can be found in \cite{FH,GW,We}, and on the quantum version in \cite{CP,Ji86,KS}. Recall that we have a representation of the Hecke algebra $H_n(q)$ on the tensor product $V^{\otimes n}$. This was given explicitly in (\ref{rep-Hn}). So our situation is that we have two algebras represented on the same vector space, and we can picture it like this:
\[U_q(sl_N)\ \ \stackrel{\rho^{(n)}}{\longrightarrow}\ \ \ \ \text{End}(V^{\otimes n})\ \ \ \ \stackrel{\pi}{\longleftarrow}\ \ H_n(q)\ \]
Now we can state the (first part of) the Schur--Weyl duality.
\begin{theo}[Schur--Weyl I]\label{SW1}
The centraliser $\text{End}_{U_q(sl_N)}(V^{\otimes n})$ is the image of the Hecke algebra $H_n(q)$:
\[\text{End}_{U_q(sl_N)}(V^{\otimes n})=\pi\bigl(H_n(q)\bigr)\ .\]
\end{theo}
Since we are in the semisimple situation, from the general consideration on centralisers sketched above, we have actually that $\pi\bigl(H_n(q)\bigr)$ and $\rho^{(n)}\bigl(U_q(sl_N)\bigr)$ are the mutual centralisers of each other. This is why we call it a duality.

\vskip .2cm
From the generalities discussed above, see around Formulas (\ref{dec-gen}), we know that the centraliser is related to various decompositions of the representation $V^{\otimes n}$ into irreducible summands. To describe this in our particular situation, we need some notations about irreducible representations of $U_q(sl_N)$ and of $H_n(q)$.

A partition $\lambda$ of $n$ is a family of integers $\lambda=(\lambda_1,\dots,\lambda_l)$ such that $\lambda_1\geq\lambda_2\geq\dots\geq\lambda_l\geq 0$ and $\lambda_1+\dots+\lambda_l=n$. We note $\lambda\vdash n$ and we say that $\lambda$ is a partition {\em of size} $n$. 
The number $\ell(\lambda)$ of non-zero parts is called the {\em length} of $\lambda$. Then we have:
\begin{itemize}
\item The quantum group $U_q(sl_N)$ has irreducible representations indexed by highest weights, which are here identified with partitions $\lambda$ such that $\ell(\lambda)\leq N$. Let us denote them by $L_{\lambda}^N$. Note that here the size of $\lambda$ can be arbitrary. We put a $N$ in the notation because the same partition can index a representation of $U_q(sl_N)$ for various $N$. The standard representation $V$ corresponds to $\lambda=(1)$. 
\item The Hecke algebra $H_n(q)$ has irreducible representations indexed by partitions $\lambda$ of size $n$. Let us denote them by $S_{\lambda}$ (our convention is such that $S_{(n)}$ is the one-dimensional representation $\sigma_i\mapsto q$).
\end{itemize}
Now the decomposition of $V^{\otimes n}$, respectively, as a representation of $U_q(sl_N)$, as a representation of $H_n(q)$ and finally as a representation of $U_q(sl_N)\otimes H_n(q)$ is:
\begin{equation}\label{dec-SW}
V^{\otimes n}=\bigoplus_{\substack{\lambda\vdash n\\[0.2em] \ell(\lambda)\leq N}} (L_{\lambda}^N)^{\oplus \dim(S_{\lambda})}\,,\quad\ \ \ V^{\otimes n}=\bigoplus_{\substack{\lambda\vdash n\\[0.2em] \ell(\lambda)\leq N}} S_{\lambda}^{\oplus \dim(L_{\lambda}^N)}\,,\quad\ \ \ V^{\otimes n}=\bigoplus_{\substack{\lambda\vdash n\\[0.2em] \ell(\lambda)\leq N}} L_{\lambda}^N\otimes S_{\lambda}\,.
\end{equation}
The duality here is apparent. Note that the knowledge of the representations $S_{\lambda}$ of $H_n(q)$ can be used to construct the representations $L_{\lambda}^N$ of $U_q(sl_N)$ (this is the natural point of view in these notes), or the other way around.

\paragraph{Kernel in the Schur--Weyl duality.}
The first step of the Schur--Weyl duality is a great step in the understanding of $\text{End}_{U_q(sl_N)}(V^{\otimes n})$. Indeed Theorem \ref{SW1} asserts that it is generated by the elements $\check R_1,\dots,\check R_{n-1}$, and that these elements, in addition to satisfying the braid relations, also satisfy the Hecke relations. However, we can not say that the centraliser \emph{is} the Hecke algebra $H_n(q)$. It only says that the centraliser is the image of the map $\pi$, and this map may have a kernel. So the second part of the Schur--Weyl duality shall be the description of this kernel. It turns out that there is a quite simple description of the kernel. The kernel was actually implicitly described in (\ref{dec-SW}) but let us not worry about that and just give directly its algebraic description.

\vskip .2cm
We need to go a little bit into the algebraic structure of the Hecke algebra $H_n(q)$. For any element $w$ of the symmetric group $\mS_n$, let $w=s_{a_1}\dots s_{a_k}$ be a reduced (\emph{i.e.} minimal length) expression for $w$ in terms of the generators $s_i=(i,i+1)$, and denote $\ell(w)=k$. Then define $\si_w:=\si_{a_1}\dots \si_{a_k}\in H_n(q)$. This definition does not depend on the reduced expression for $w$ and the set $\{\si_w\}_{w\in\mS_n}$ forms a basis of $H_n(q)$.

For $L\in\mathbb{Z}_{\geq 0}$, we define the $q$-numbers as follows:
\begin{equation}\label{quantum-numbers}
[L]_q:=\frac{q^L-q^{-L}}{q-q^{-1}}=q^{L-1}+q^{L-3}+\dots+q^{-(L-1)}\ \ \ \ \ \text{and}\ \ \ \ \ \ \ [L]_q!:=[1]_q[2]_q\dots[L]_q\ .
\end{equation}
Then the \textbf{$q$-symmetriser} in $H_n(q)$ is the following element:
\begin{equation}\label{def-P}
P_n=\frac{\sum_{w\in\mS_n}q^{\ell(w)}\si_w}{\sum_{w\in\mS_n}q^{2\ell(w)}}=\frac{q^{-n(n-1)/2}}{[n]_q!}\sum_{w\in\mS_n}q^{\ell(w)}\si_w\ .
\end{equation}
and the \textbf{$q$-antisymmetriser} in $H_n(q)$ is the following element:
\begin{equation}\label{def-P'}
P'_n=\frac{\sum_{w\in\mS_n}(-q^{-1})^{\ell(w)}\si_w}{\sum_{w\in\mS_n}q^{-2\ell(w)}}=\frac{q^{n(n-1)/2}}{[n]_q!}\sum_{w\in\mS_n}(-q^{-1})^{\ell(w)}\si_w\ .
\end{equation}
We take this opportunity to emphasize the nice formula $\sum_{w\in\mS_n}q^{2\ell(w)}=q^{n(n-1)/2}[n]_q!$.
\begin{rema}
For $q=1$, the $q$-symmetriser becomes the usual symmetriser in the group algebra $\bC\mS_n$: the sum of all elements of $\mS_n$ divided by $n!$. Similarly, the $q$-antisymmetriser becomes the usual antisymmetriser: the sum of all elements of $\mS_n$ multiplied by their signature and divided by $n!$.
\end{rema}

We remark that up to now the description of the centraliser does not explicitly depend on the dimension $N$ of $V$, since the Hecke algebra does not depend on $N$. Of course, this dependence is hidden in the map $\pi$ and in its kernel. So it is only natural that it appears explicitly now.
\begin{theo}[Schur--Weyl II]\label{SW2}
\begin{itemize}
\item If $n\leq N$ then the kernel of the map $\pi$ is $\{(0)\}$  
\item If $n>N$ then the kernel of the map $\pi$ is generated in $H_n(q)$ by the element $P'_{N+1}$.
\end{itemize}  
\end{theo}
Note that we see the element $P'_{N+1}$ as an element of $H_n(q)$ by the natural inclusion of $H_{N+1}(q)$ in $H_n(q)$ (if $n>N$). In words, to obtain the complete description of the centraliser, we start from the Hecke algebra $H_n(q)$ and, if $n>N$, we cancel the $q$-antisymmetriser on $N+1$ letters.

Explicitly, the centraliser $\text{End}_{U_q(sl_N)}(V^{\otimes n})$ is isomorphic to the algebra generated by $\si_1,\dots,\si_{n-1}$ with defining relations:
\[\begin{array}{ll}
\si_i\si_{i+1}\si_i=\si_{i+1}\si_i\si_{i+1}\,,\ \ \  & \text{for $i\in\{1,\dots,n-2\}$}\,,\\[0.2em]
\si_i\si_j=\si_j\si_i\,,\ \ \  & \text{for $i,j\in\{1,\dots,n-1\}$ such that $|i-j|>1$}\,,\\[0.2em]
\si_i^2=1+(q-q^{-1})\si_i\,,\ \ \  & \text{for $i\in\{1,\dots,n-1\}$}\,,\\[0.2em]
P'_{N+1}=0 & \text{if $n>N$.} 
\end{array}
\]

\begin{exam}[Temperley--Lieb algebra]\label{exa-TL}
Let $N=2$. In this case, the centraliser $\text{End}_{U_q(sl_2)}(V^{\otimes n})$ is called the Temperley--Lieb algebra. The additional relation $P'_3=0$ reads:
\[1-q^{-1}(\si_1+\si_2)+q^{-2}(\si_1\si_2+\si_2\si_1)-q^{-3}\si_1\si_2\si_1=0\ .\]
One can show that it implies the same relation with indices $i,i+1$ for all $i=1,\dots,n-2$. Then setting $\tau_i:=\si_i-q$, one recovers the other standard presentation of the Temperley--Lieb algebra:
\[\tau_i^2=-(q+q^{-1})\tau_i\,,\ \ \ \ \tau_i\tau_{i+1}\tau_i=\tau_i\,,\ \ \ \ \tau_{i+1}\tau_{i}\tau_{i+1}=\tau_{i+1}\ \ \ \text{and}\ \ \ \tau_i\tau_j=\tau_j\tau_i\ \text{if $|i-j|>1$.}\]
\end{exam}

\paragraph{And the BMW algebra?}
Now we explained, at length, how the Hecke algebra fits in the story of centralisers of quantum groups representations. This beautiful story has its counterpart for the BMW algebra. Indead instead of $sl_N$ we can consider the other classical Lie algebras $so_N$ or $sp_N$. We can form the quantum groups $U_q(so_N)$ or $U_q(sp_N)$, and we can consider the tensor product $V^{\otimes n}$, where $V$ is the (analogue of the) vector representation for $U_q(so_N)$ or $U_q(sp_N)$. Then the centralisers of $V^{\otimes n}$ are described in a way similar to above, where the Hecke algebra is replaced by another algebra. This other algebra is the BMW algebra (with the parameter $a$ specialised to a power of $q$ depending whether we have $so_N$ or $sp_N$). For details about that along the same lines as above, we refer for example to \cite[\S 8.6]{KS}. This concludes our discussion on the interpretation of Hecke and BMW algebras as centralisers of representations of quantum groups.

\newpage
\section{Fusion Procedure for the Yang--Baxter Equation}\label{sec-fusion}

The fusion procedure was designed very early in the history of the YB equation \cite{Ji86,Ji-int,KRS}. It consists in a general procedure to construct new solutions starting from a known one. This procedure can be split in two steps and presented schematically like in the following table:

\medskip
\begin{center}
\begin{tabular}{r|c|l}

\textbf{Vector spaces} & \textbf{Matrices} & \textbf{Algebras} \\
\hline & & \\

$V$ \ \ & basic solution on $V\otimes V$ & $\hookleftarrow\ $ Hecke algebra $H_n(q)$\\[0.5em]

(generic fusion) $\downarrow$ \ \ \   & $\downarrow$ & \\[0.5em]

$V^{\otimes k}$ & solution on $V^{\otimes k}\otimes V^{\otimes k}$ & $\hookleftarrow\ $ (bigger) Hecke algebra $H_{kn}(q)$\\[0.5em]

(projection) $\downarrow$ \ \ \  & $\downarrow$ &\\[0.5em]

$Pr(V^{\otimes k})$ & \ \ solution on $Pr(V^{\otimes k})\otimes Pr(V^{\otimes k})$\ \  & $\hookleftarrow\ $ \emph{?? fused Hecke algebra ??}

\end{tabular}
\end{center}
\medskip

To give a rough idea of what is going on, the first two columns describe the fusion procedure that was originally designed, for matrix solutions of the YB equation. Starting from a known solution on $V$, an explicit formula builds a new solution on the much larger vector space $V^{\otimes k}$. This is the first step. The second step identifies invariant subspaces in $V^{\otimes k}$ for this new solution. This is quite general and valid for any solution $R(u)$ of the YB equation.

It is important to indicate that we have some free parameters in the fusion procedure, and it turns out that for some specific choices of these free parameters, the invariant subspaces can be especially interesting. To give at once the explicit example we will be interested in, say the solution on $V$ comes from the Hecke algebra and look at the tensor product $V^{\otimes n}$ as a representation of $sl_N$ (or $U_q(sl_N)$), where $N=\dim(V)$. Then it decomposes as:
\[V^{\otimes k}=S^k(V)\oplus \dots\,,\]
where $S^k(V)$ is the $k$-th symmetric power of $V$ (or its analogue for $U_q(sl_N)$). If we choose wisely and precisely the parameters in the fusion procedure, then the big solution $R^{(k)}(u)$ restricts to a solution on the subspace $S^k(V)$. This is for the matrix side of the picture.

\vskip .2cm
The third column intends to indicate what should be the ``algebra'' counterpart of this. Well it turns out that the first step of the fusion procedure can be carried out directly in the algebra. So if we started in $H_n(q)$, we end up in the larger Hecke algebra $H_{kn}(q)$. The number of strands has been multiplied by $k$, this is called cabling from the point of view of braids.

It also turns out that the projector, denoted here $Pr$, on the subspace $S^k(V)$, can be seen directly in the Hecke algebra $H_k(q)$. In fact this is part of the statement of the Schur--Weyl duality, and it turns out that in this case $Pr$ is simply the (image of the) $q$-symmetriser of $H_k(q)$. With a little algebraic thinking on how to project a representation, we realise that the algebra controlling the new fused solution on $S^k(V)$ should be:
\[\text{``fused Hecke algebra''}\ =\ P_{k,n}.H_{kn}(q).P_{k,n}\ ,\]
where the idempotent $P_{k,n}$ in $H_{kn}(q)$ is made up by plugging the $q$-symmetriser on the first $k$ strands, and also on the next $k$ strands and so on up to the last $k$ strands. If we want to conclude this section now, we can sum up the algebraic procedure like this:
\begin{conc}
The fused Hecke algebra is obtained with the following two steps:
\[H_{n}(q)\ \ \stackrel{\text{``cabling''}}{\longrightarrow}\ \ H_{kn}(q)\ \ \stackrel{\text{``projecting''}}{\longrightarrow}\ \ P_{k,n}.H_{kn}(q).P_{k,n}\ \ \text{($=$ fused Hecke algebra)}\]
\end{conc}

If one is satisfied enough with this brief description, one can proceed directly to Section \ref{sec-fusedH} for a description of the fused Hecke algebra (and should not be surprised to find a solution of the YB equation in there). Otherwise, in the present section, we shall present more precisely the machinery of the fusion procedure. More details can be found in \cite{PdA1}.

\subsection{Fusion procedure for matrix solutions}

In what follows, we make the choice to work with multiplicative spectral parameters in the YB equation. This is more adapted to the Hecke algebra situation. Of course, everything is also valid in the additive case with suitable modifications.

\paragraph{First step.} We start with an arbitrary solution of the YB equation on a vector space $V$. Here we use the following form of the YB equation:
\[
R_{12}(u)R_{13}(uv)R_{23}(v)=R_{23}(v)R_{13}(uv)R_{12}(u)\ \quad \ \ \ \ \text{on $V\otimes V\otimes V$.} 
\]
Let $\bc:=(c_1,\dots,c_k)$ be a $k$-tuple of non-zero complex parameters. We consider the space $V^{\otimes k}\otimes V^{\otimes k}$ and we find it convenient to label the copies of $V$ by $1,\dots,k,\un{1},\dots,\un{k}$ (from left to right). For example, the operator $R_{a,\un{b}}(u)$  stands for the operator $R_{a,k+b}(u)$ with the notation used before. 

Then we define the following operator in $\text{End}(V^{\otimes k}\otimes V^{\otimes k})$.
\begin{equation}\label{def-fus-R}
R^{(\bc)}(u):=\prod_{i=1,\dots,k}^{\leftarrow}R_{1,\un{i}}(u\frac{c_1}{c_i})R_{2,\un{i}}(u\frac{c_2}{c_i})\dots\dots R_{k,\un{i}}(u\frac{c_k}{c_i})\ ,
\end{equation}
where the arrow means that the product is taken from left to right in decreasing order in the index $i$ (the first factor is $R_{1,\un{n}}(u\frac{c_1}{c_k})$ and the last is $R_{n,\un{1}}(u\frac{c_k}{c_1})$). 

By a direct calculation using repeatedly the YB equation for $R(u)$, one can check that this new operator also satisfies the YB equation:
\[
R^{(\bc)}_{12}(u)R^{(\bc)}_{13}(uv)R^{(\bc)}_{23}(v)=R^{(\bc)}_{23}(v)R^{(\bc)}_{13}(uv)R^{(\bc)}_{12}(u)\ \quad \ \ \ \text{on $V^{\otimes k}\otimes V^{\otimes k}\otimes V^{\otimes k}$,} 
\]
where now the indices $1,2,3$ refer to the new vector space $V^{\otimes k}$.

\begin{exam}
Within the ``space-time trajectories'' interpretation of the YB equation in (\ref{YB-traj}), this first step of the fusion procedure is easy to understand. The new solution, given by Formula (\ref{def-fus-R}), corresponds to interactions of multiplets of $k$ particles. We have a ``beam'' of $k$ particles crossing another beam, and the operator $R^{(k)}(u)$ describes this interaction, factorised as a product of $2$-particles interaction. For example, for $k=2$, the operator $R^{(2)}(u)$ corresponds to:
\begin{center}
 \begin{tikzpicture}[scale=0.2]
\node at (0,11) {$\scriptstyle{1}$};
\node at (5,11) {$\scriptstyle{2}$};
\node at (15,11) {$\scriptstyle{\un{1}}$};
\node at (20,11) {$\scriptstyle{\un{2}}$};
\draw[thick] (0,10) -- (15,-10);
\draw[thick] (5,10) -- (20,-10);
\draw[thick] (15,10) -- (0,-10);
\draw[thick] (20,10) -- (5,-10);
\node at (10,5.2) {$\scriptstyle{u_{2\un{1}}}$};
\draw (9.2,4.2) -- (10.8,4.2);
\node at (7.5,1.9) {$\scriptstyle{u_{1\un{1}}}$};
\draw (6.7,0.9) -- (8.3,0.9);
\node at (12.5,1.9) {$\scriptstyle{u_{2\un{2}}}$};
\draw (11.7,0.9) -- (13.3,0.9);
\node at (10,-1.4) {$\scriptstyle{u_{1\un{2}}}$};
\draw (9.2,-2.4) -- (10.8,-2.4);
\node at (25,0){$:$};
\node at (50,0){$\ \ \ R^{(2)}(u)=R_{1\un{2}}(u_{1\un{2}})R_{2\un{2}}(u_{2\un{2}})R_{1\un{1}}(u_{1\un{1}})R_{2\un{1}}(u_{2\un{1}})\ ,$};
\end{tikzpicture}
\end{center}
where $u_{i\un{j}}:=u\displaystyle\frac{c_i}{c_j}$. The parameters $\bc=(c_1,\dots,c_k)$ accounts for the different rapidities of the particles. Note that we should have taken in general a different set $\un{\bc}=(c_{\un{1}},\dots,c_{\un{k}})$ of parameters for the particles $\un{1},\dots,\un{k}$. For simplicity, and as it is enough for what follows, we took the same set $\bc$ twice.

So this first step amounts to the remarkable property that if the one-particle interaction satisfies the YB equation then so does the $k$-particles interaction.
\end{exam}

\begin{rema}
The extension from one particle to $k$ particles is reminiscent in braid theory of a particular case of what is called \emph{cabling}. Start from a braid, and duplicate all strands, replacing a strand by $k$ new strands. This gives a new braid. For example, on the elementary braiding for $k=2$:
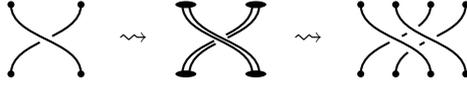
\begin{center}
 \begin{tikzpicture}[scale=0.23]
\fill[black] (0,2) circle (0.2);
\fill[black] (0,-2) circle (0.2);
\fill[black] (4,2) circle (0.2);
\fill[black] (4,-2) circle (0.2);

\draw[thick] (4,2)..controls +(0,-2) and +(0,+2) .. (0,-2);
\fill[white] (2,0) circle (0.4);
\draw[thick] (0,2)..controls +(0,-2) and +(0,+2) .. (4,-2);
 
\node at (7,0) {$\rightsquigarrow$}; 
 
\fill (10,2) ellipse (0.6cm and 0.2cm);\fill (10,-2) ellipse (0.6cm and 0.2cm);
\fill (14,2) ellipse (0.6cm and 0.2cm);\fill (14,-2) ellipse (0.6cm and 0.2cm);

\draw[thick] (13.8,2)..controls +(0,-2) and +(0,+2) .. (9.8,-2);
\draw[thick] (14.2,2)..controls +(0,-2) and +(0,+2) .. (10.2,-2);
\fill[white] (12,0) circle (0.4);
\draw[thick] (10.2,2)..controls +(0,-2) and +(0,+2) .. (14.2,-2);
\draw[thick] (9.8,2)..controls +(0,-2) and +(0,+2) .. (13.8,-2);

\node at (17,0) {$\rightsquigarrow$}; 

\fill[black] (20,2) circle (0.2);
\fill[black] (20,-2) circle (0.2);
\fill[black] (22,2) circle (0.2);
\fill[black] (22,-2) circle (0.2);
\fill[black] (24,2) circle (0.2);
\fill[black] (24,-2) circle (0.2);
\fill[black] (26,2) circle (0.2);
\fill[black] (26,-2) circle (0.2);

\draw[thick] (24,2)..controls +(0,-2) and +(0,+2) .. (20,-2);
\draw[thick] (26,2)..controls +(0,-2) and +(0,+2) .. (22,-2);
\fill[white] (22,0) circle (0.4);
\fill[white] (24,0) circle (0.4);
\fill[white] (23,0.5) circle (0.4);
\fill[white] (23,-0.5) circle (0.4);
\draw[thick] (20,2)..controls +(0,-2) and +(0,+2) .. (24,-2);
\draw[thick] (22,2)..controls +(0,-2) and +(0,+2) .. (26,-2);
\end{tikzpicture}
\end{center}
It is immediate that the new braids obtained from the elementary braidings also satisfy the braid relations. This is the analogue of the first step of the fusion procedure.
\end{rema}

\paragraph{Second step.} Now we consider the following operator on $V^{\otimes k}$:
\begin{equation}\label{F}
F(\bc):=\prod_{1\leq i<j\leq k}^{\rightarrow}R_{j,i}(\frac{c_j}{c_i})\ ,
\end{equation}
where the arrow indicates that the product is taken from left to right in increasing order for the lexicographic order: $F(\bc)=R_{21}(\frac{c_2}{c_1})R_{31}(\frac{c_3}{c_1})\dots R_{k1}(\frac{c_k}{c_1})\cdot R_{32}(\frac{c_3}{c_2})......R_{k,k-1}(\frac{c_k}{c_{k-1}})$\ . The subspace we are interested in is the image of this operator $F(\bc)$:
\begin{equation}\label{W}
W_{\bc}:=\textrm{Im}(F(\bc))\subset V^{\otimes k}\ .
\end{equation}
We note that $F(\bc)$ may not be defined for any value of the parameters $\bc$, since the function $R(u)$ can have singularities (for example, the solution coming from the Hecke algebra can not be evaluated at $u=1$). Nevertheless, we assume in our discussion that the set of parameters $\bc$ is such that the operator $F(\bc)$ is a well-defined element of $\text{End}(V^{\otimes k})$.

So now we have a subspace $W_{\bc}\otimes W_{\bc}$ of $V^{\otimes k}\otimes V^{\otimes k}$. The main fact about it is that it is an invariant subspace for the operator $R^{(\bc)}(u)$. In other words, the operator $R^{(\bc)}(u)$ restricts to $W_{\bc}\otimes W_{\bc}$ and provides a solution of the YB equation on this subspace.
\begin{conc}\label{conc-fus}
For any set of parameters $\bc=(c_1,\dots,c_k)$ such that the subspace $W_{\bc}\subset V^{\otimes n}$ is defined, we have a solution $R^{fus}(u)$ of the YB equation on $W_{\bc}\otimes W_{\bc}$.
\end{conc}
We denoted by $R^{fus}(u)$ the restriction to $W_{\bc}\otimes W_{\bc}$ of $R^{(\bc)}(u)$, and we call this element $R^{fus}(u)$ a fused solution.

\vskip .2cm
We should note that for parameters $\bc$ such that the operator $F(\bc)$ is invertible then the subspace $W_{\bc}$ is in fact $V^{\otimes k}$ and the second step has given us nothing new. Of course, if we went into all this trouble, it is because we hope to find more interesting subspaces than that. We will show convincing examples below, in the Hecke algebra situation.

\subsection{Fusion procedure, the braided version}
Recall from Section \ref{subsec-YB} that the YB equation for $R(u)$ is equivalent to the braided YB equation for $\check R(u)$, where both are related by the permutation operator $P$ of $V\otimes V$:
\[R(u)=P\check R(u)\ .\]
In our setting, it is the braided version $\check R(u)$ which can be seen as coming from a Baxterization formula of an algebra. While $R(u)$ a priori can not be expressed directly in the algebra, because the permutation operator $P$ may have no meaning in the algebra. So our goal should be first to express the fusion procedure using only the ``braided'' operator $\check R(u)$. It turns out to be possible.

\vskip .2cm
In what follows, we use indices $1,\dots,2k$ to label the factors in $V^{\otimes k}\otimes V^{\otimes k}=V^{\otimes 2k}$.
For the first step of the fusion procedure, a direct calculation shows that
\[R^{(\bc)}(u)=P_{V^{\otimes k},V^{\otimes k}}\cdot \prod_{i=1,\dots,k}^{\leftarrow} \check R_i(u\frac{c_1}{c_i})\check R_{i+1}(u\frac{c_2}{c_i})\dots\dots \check R_{i+k-1}(u\frac{c_k}{c_i})\ ,\]
where $P_{V^{\otimes k},V^{\otimes k}}$ is the permutation operator on $V^{\otimes k}\otimes V^{\otimes k}$. But of course, multiplication on the left by $P_{V^{\otimes k},V^{\otimes k}}$ is just a way to go the braided YB equation. So we conclude that we could just have started with the following definition\footnote{However, the nice interpretation with interacting particles would have been less apparent.}:
\begin{equation}\label{fus-op}
\check R^{(\bc)}(u)=\prod_{i=1,\dots,k}^{\leftarrow} \check R_i(u\frac{c_1}{c_i})\check R_{i+1}(u\frac{c_2}{c_i})\dots\dots \check R_{i+k-1}(u\frac{c_k}{c_i})\ .
\end{equation}
And the first step of the fusion procedure asserts that the operator $\check R^{(\bc)}(u)$ satisfies the braided YB equation:
\[\check R^{(\bc)}_1(u)\check R^{(\bc)}_2(uv)\check R^{(\bc)}_1(v)=\check R^{(\bc)}_2(v)\check R^{(\bc)}_1(uv)\check R^{(\bc)}_2(u)\quad\ \ \ \ \ \text{on $V^{\otimes k}\otimes V^{\otimes k}\otimes V^{\otimes k}$.}\]

Now for the second step, we should try to express the operator $F(\bc)$ also in terms of the braided element $\check R(u)$. This is indeed possible, a direct calculation showing that:
\[F(\bc)=\Bigl(\prod_{i=1,\dots,k-1}^{\rightarrow} \check R_i(\frac{c_{i+1}}{c_1})\check R_{i-1}(\frac{c_{i+1}}{c_{2}})\dots\dots\check R_1(\frac{c_{i+1}}{c_i})\Bigr)\cdot P_{w_0}\ \ \ \ \quad\ \text{on $V^{\otimes k}$,}\]
where $P_{w_0}$ is the operator on $V^{\otimes k}$ permuting the factors according to the permutation completely reversing the indices: $(1,\dots,k)$ goes to $(k,\dots,1)$. The map $P_{w_0}$ is invertible so it sends $V^{\otimes k}$ ot itself, and therefore the image of $F(\bc)$ is the same as the image of the element in the big parenthesis just above.
Thus we conclude that we could just have started with the following definition:
\begin{equation}\label{fus-F}
\check F(\bc)=\prod_{i=1,\dots,k-1}^{\rightarrow} \check R_i(\frac{c_{i+1}}{c_1})\check R_{i-1}(\frac{c_{i+1}}{c_{2}})\dots\dots\check R_1(\frac{c_{i+1}}{c_i})\ ,
\end{equation}
and define
\begin{equation}\label{W2}
W_{\bc}:=\textrm{Im}(\check F(\bc))\subset V^{\otimes k}\ .
\end{equation}
We have explained that this definition of $W_{\bc}$ coincides with the original one in (\ref{W}). The second step of the fusion procedure asserts that the operator $\check R^{(\bc)}(u)$ leaves the subspace $W_{\bc}\otimes W_{\bc}$ invariant. In fact, we can show that:
\[\check R^{(\bc)}(u)\ \text{commutes with } \check F(\bc)\otimes \check F(\bc)\,. \]
We denote by $\check R^{fus}(u)$ the restriction to $W_{\bc}\otimes W_{\bc}$ of $\check R^{(\bc)}(u)$, and we also call this element $\check R^{fus}(u)$ a fused solution. At the end of the day, of course we have:
\[\check R^{fus}(u)=P_{W_{\bc}\otimes W_{\bc}}R^{fus}(u)\ \ \ \ \ \ \ \ \text{on $W_{\bc}\otimes W_{\bc}$\ ,}\]
but our point in this subsection was mostly to express $\check R^{fus}(u)$ and $\check F(\bc)$ in terms of $\check R(u)$.

\subsection{Fusion procedure in the Hecke algebra}

We will now describe the procedure directly in the Hecke algebra. We deal with the example of the Hecke algebra to simplify the discussion, but part of it, especially the beginning, can be formulated for an arbitrary algebra with a Baxterization formula and a local representation.

So we assume now that the solution $\check R(u)$ comes from the Hecke algebra. It means that $\check R_i(u)$ is the image of $\si_i(u)=\si_i+(q-q^{-1})\frac{1}{u-1}$ in the local representation of the Hecke algebra $H_n(q)$ on $V^{\otimes n}$, see Section \ref{subsec-Hec}.

\vskip .2cm
The first step of the fusion procedure is easy to formulate in the algebra. From its definition in (\ref{fus-op}), it is immediate that the new solution $\check R^{(\bc)}(u)$ in $\text{End}(V^{\otimes k}\otimes V^{\otimes k})$ is the image of the following element of $H_{2k}(q)$:
\begin{equation}\label{fus-op-alg}
\si_1^{(\bc)}(u)=\prod_{i=1,\dots,k}^{\leftarrow} \si_i(u\frac{c_1}{c_i})\si_{i+1}(u\frac{c_2}{c_i})\dots\dots\si_{i+k-1}(u\frac{c_k}{c_i})\ .
\end{equation}
By shifting the indices by $k$, we can easily define $\si_2^{(\bc)}(u)$, and in turn $\si_3^{(\bc)}(u)$ and so on. The elements $\si_1^{(\bc)}(u),\dots,\si_{n-1}^{(\bc)}(u)$ satisfy the braided Yang--Baxter equation
\[\sigma^{(\bc)}_i(u)\sigma^{(\bc)}_{i+1}(uv)\sigma^{(\bc)}_i(v)=\sigma^{(\bc)}_{i+1}(v)\sigma^{(\bc)}_i(uv)\sigma^{(\bc)}_{i+1}(u)\ ,\]
directly in the algebra, which is here $H_{kn}(q)$ (the number of strands has simply been multiplied by $k$). The local representation of the algebra $H_{kn}(q)$ is naturally on $(V^{\otimes k})^{\otimes n}$, so we have interpreted the solution $\check R^{(\bc)}(u)$ as coming from an algebra with a local representation.

\vskip .2cm
Then we note that from its definition (\ref{fus-F}), the operator $\check F(\bc)$ in $\text{End}(V^{\otimes k})$ is the image of the following element of $H_{k}(q)$:
\begin{equation}\label{fus-F-alg}
\Phi(\bc)=\prod_{i=1,\dots,k-1}^{\rightarrow} \si_i(\frac{c_{i+1}}{c_1})\si_{i-1}(\frac{c_{i+1}}{c_{2}})\dots\dots\si_1(\frac{c_{i+1}}{c_i})\ .
\end{equation}
What does it mean in the Hecke algebra that $\check R^{(\bc)}(u)$ leaves invariant the subspace $W_{\bc}\otimes W_{\bc}$? Well, first we have to note that $H_k(q)\otimes H_k(q)$ is included in $H_{2k}(q)$. The first factor is generated by $\si_1,\dots,\si_{k-1}$, the second factor is generated by $\si_{k+1},\dots,\si_{2k-1}$. So we have an element $\Phi(\bc)\otimes \Phi(\bc)$ in $H_{2k}(q)$. Then the second step of the fusion procedure, directly inside the Hecke algebra is that:
\begin{equation}\label{secstep-alg}
[\si_1^{(\bc)}(u),\Phi(\bc)\otimes \Phi(\bc)]=0\ \ \ \ \ \text{in $H_{2k}(q)$.}
\end{equation}
This formula can be obtained directly from a straightforward calculation using repeatedly the Yang--Baxter equation for $\sigma_i(u)$.

To conclude and summarize the description of the fusion procedure inside the algebra:
\begin{center}
\begin{tabular}{c|c|c}
& Algebra & Representation\\[0.5em]
\hline & & \\[-0.2em]
1st step & $\sigma_i^{(\bc)}(u)$: satisfies YB in $H_{kn}(q)$ & $\check R_i^{(\bc)}(u)$: satisfies YB on $(V^{\otimes k})^{\otimes n}$\\[0.8em]
2nd step & $\sigma_1^{(\bc)}(u)$ commutes with $\Phi(\bc)\otimes \Phi(\bc)$ &  $\check R^{(\bc)}(u)$ commutes with $\check F(\bc)\otimes \check F(\bc)$\\[-0.2em]
 && $\Downarrow$\\[-0.2em]
 && $\check R_i^{(\bc)}(u)$ leaves invariant $(W_{\bc})^{\otimes n}$
\end{tabular}
\end{center}

\paragraph{Interesting subspaces $W_{\bc}$.} Finally, we are ready to study explicit examples of interesting subspaces $W_{\bc}$. To do so, we have to study the function $\Phi(\bc)$ and to find particular specialisations of the parameters $\bc$ which makes it remarkable.

\vskip .2cm
First, a tiny bit of combinatorics notations. A pair $(x,y)\in\mathbb{Z}^2$ is called a {\em node}. The $q$-content of the node $\theta=(x,y)$ is $\qc(\theta)=q^{2(y-x)}$. 
 
We have already defined the notion of a partition $\lambda\vdash k$. The Young diagram of $\lambda=(\lambda_1,\dots,\lambda_l)$ is the set of nodes $(x,y)$ such that $x\in\{1,\dots,l\}$ and $y\in\{1,\dots,\lambda_x\}$. The Young diagram of $\lambda$ will be seen as a left-justified array of $l$ rows such that the $j$-th row contains $\lambda_j$ nodes for all $j=1,\dots,l$ (a node will be pictured by an empty box). 

A standard Young tableau of shape $\lambda\vdash k$ is an assignment of a number from $\{1,\dots,k\}$ to each node of the Young diagram of $\lambda$, such that the numbers are strictly ascending along rows and down columns of the Young diagram. It is represented by filling the nodes of the Young diagram of $\lambda$ with the numbers from $\{1,\dots,k\}$. Here is an example
\[\bT=\begin{array}{ccc}
\fbox{\scriptsize{$1$}} & \hspace{-0.35cm}\fbox{\scriptsize{$2$}} & \hspace{-0.35cm}\fbox{\scriptsize{$4$}} \\[-0.2em]
\fbox{\scriptsize{$3$}} & &
\end{array}\ : \qquad \qc_1(\bT)=1\,,\ \ \ \qc_2(\bT)=q^2\,,\ \ \ \qc_3(\bT)=q^{-2}\,,\ \ \ \qc_4(\bT)=q^4\,.\]
The partition is $\lambda=(3,1)$. Next to the Young tableau $\bT$ is its sequence of $q$-contents, where $\qc_i(\bT)$ is the $q$-content of the node of $\bT$ with number $i$.

Here is the result providing nice specialisations of $\Phi(\bc)$. Below, $\si_{w_0}^{-1}$ is the element of the standard basis of the Hecke algebra corresponding to the longest element $w_0$ of the symmetrix group $\mS_k$. Its presence here is harmless to us since, being an invertible element, in a representation it will not change the subspace obtained as the image of $\Phi(\bc)$. Classical versions for the symmetric group can be found in \cite{Ch2,Mol,Na}.
\begin{theo}[\cite{IMOs}]\label{thm-idem-Hn}
Let $\bT$ be a standard Young tableau of shape $\lambda\vdash k$. The element obtained by the following consecutive evaluations
\begin{equation}\label{idem-Hn}
\Phi(\bc)\si_{w_0}^{-1}\Bigr\vert_{c_1=\qc_1(\bT)}\Bigr\vert_{c_2=\qc_2(\bT)}\dots \Bigr\vert_{c_k=\qc_k(\bT)}
\end{equation}
is non-zero and proportional to an idempotent $E_{\bT}$ of $H_k(q)$.
\end{theo}
In fact, we can say more about $E_{\bT}$. The idempotent $E_{\bT}$ is a primitive idempotent corresponding to the irreducible representation $S_{\lambda}$ of $H_k(q)$. That is equivalent to say that the action of $E_{\bT}$ in an irreducible representation of $H_k(q)$ is as follows:
\[E_{\bT}(S_{\lambda})\ \text{is a one-dimensional subspace}\quad \ \ \ \text{and}\quad\ \ \ \ E_{\bT}(S_{\mu})=\{0\}\ \ \text{if $\mu\neq\lambda$}.\]
So finally, take the parameters $\bc$ as the sequence of contents of a standard Young tableau of some shape $\lambda\vdash k$. Recall from the Schur--Weyl duality that the tensor product $V^{\otimes k}$ decomposes as follows:
\[ V^{\otimes k}=\bigoplus_{\substack{\mu\vdash k\\[0.2em] \ell(\mu)\leq N}} L_{\mu}^N\otimes S_{\mu}\,,\]
as a representation of $U_q(sl_N)\otimes H_k(q)$. To determine the subspace $W_{\bc}$, we need to calcuate the image in this representation of $\Phi(\bc)$, which is the same as the image of the idempotent $E_{\bT}$. From what we have said above, we can conclude:
\[E_{\bT}(V^{\otimes k})=\bigoplus_{\substack{\mu\vdash k\\[0.2em] \ell(\mu)\leq N}} L_{\mu}^N\otimes E_{\bT}(S_{\mu})=L_{\lambda}^N\ .\]
\begin{conc}
If the parameters $\bc$ correspond to a standard Young tableau of shape $\lambda\vdash k$, then the subspace $W_{\bc}$ coincides with the irreducible representation $L^N_{\lambda}$ of $U_q(sl_N)$.
\vskip .2cm
\begin{center}
\!\!\!\!$\leadsto$\ \ The fused solution $\check R^{fus}(u)$ is thus a solution on $L_{\lambda}^N\otimes L_{\lambda}^N$.
\end{center}
\end{conc}

\paragraph{What is the algebra for $\check R^{fus}(u)$?} We still have to understand what could be the algebra controlling the fused solution $\check R^{fus}(u)$ on $L_{\lambda}^N\otimes L_{\lambda}^N$. So we need to find an algebra with a solution of the YB equation and a representation on $(L_{\lambda}^N)^{\otimes n}$ which recovers $\check R^{fus}(u)$. This is all implicit in what was said above, so let us make it explicit.

\vskip .2cm
Consider first the algebra $H_{kn}(q)$ in which the ``big'' solution $R^{(\bc)}(u)$ lives. In the representation, this solution was projected on an invariant subspace to obtain $R^{(fus)}(u)$. This projection corresponds to the following algebraic construction.

There is a subalgebra $H_k(q)\otimes\dots\otimes H_k(q)$ ($n$ factors) in $H_{kn}(q)$. The first factor is generated by $\si_1,\dots,\si_{k-1}$, the second factor is generated by $\si_{k+1},\dots,\si_{2k-1}$, and so on. So we use the idempotent $E_{\bT}$ of $H_k(q)$ to define the following element of $H_{kn}(q)$:
\[P_{\bT,n}=E_{\bT}\otimes\dots\otimes E_{\bT}\ .\]
This is clearly an idempotent of $H_{kn}(q)$. Then consider the subset 
$$H_{\bT,n}^{fus}(q)=P_{\bT,n}H_{kn}(q)P_{\bT,n}=\{P_{\bT,n}xP_{\bT,n}\ |\ x\in H_{kn}(q)\}$$ 
Let us say a few words about subsets of this form. First, this subset is obviously stable by linear combinations and by multiplication. Moreover, since $P_{\bT,n}$ is idempotent, then $P_{\bT,n}$ is the unit element for this multiplication. So $H_{\bT,n}^{fus}(q)$ is in fact an associative unital algebra; it is a subalgebra of $H_{kn}(q)$ but its unit element is the idempotent $P_{\bT,n}$.

Algebras like $P_{\bT,n}H_{kn}(q)P_{\bT,n}$ (a two-sided multiplication by an idempotent) are classical objects in representation theory, see for example in \cite[\S 6.2]{Gr}. Their main property for our story is the following. Consider any representation $\rho\ :\ H_{kn}(q)\mapsto\text{End}(M)$ and set $N:=\rho(P_{\bT,n})(M)$ the image of the idempotent in the representation $M$. The subspace $N$ is naturally a representation of the algebra $H^{fus}_{\bT,n}(q)$. Indeed $N$ is obviously invariant under the action of any element of the form $\rho(P_{\bT,n}xP_{\bT,n})$, and thus the action of $H^{fus}_{\bT,n}(q)$ on $N$ is given simply by restriction:
\begin{equation}\label{act-PAP}
\begin{array}{rcl}
H^{fus}_{\bT,n}(q)=P_{\bT,n}H_{kn}(q)P_{\bT,n} & \to & \text{End}\bigl(N\bigr) \\[0.5em]
P_{\bT,n}xP_{\bT,n} & \mapsto & {\rho(P_{\bT,n}xP_{\bT,n})}_{|_{W}}
\end{array}\ .
\end{equation}
It only remains to apply this to the representation $V^{\otimes kn}$ of $H_{kn}(q)$. Write the vector space like $V^{\otimes kn}=V^{\otimes k}\otimes\dots\otimes V^{\otimes k}$. From the definition of the idempotent $P_{\bT,n}$, it is clear that its image on $V^{\otimes kn}$ is the tensor product of the images of $E_{\bT}$ on $V^{\otimes k}$. But this is $L_{\lambda}^N$ as we already explained above. So we conclude that we have a representation:
\[H^{fus}_{\bT,n}(q)\ \to\ \text{End}\bigl((L_{\lambda}^N)^{\otimes n}\bigr)\ .\]

Now about a solution of the YB equation, recall that we have constructed a solution $\si_i^{(\bc)}(u)$ inside $H_{kn}(q)$ (this was the first step of the fusion procedure). If the parameters $\bc$ correspond to a standard Young tableau of shape $\lambda\vdash k$, then the second step (\ref{secstep-alg}) asserts that this solution $\si_i^{(\bc)}(u)$ commutes with the idempotent $P_{\bT,n}$. Thus and finally, we can send  $\si_i^{(\bc)}(u)$ inside $H^{fus}_{\bT,n}(q)$ and define:
\begin{equation}\label{fussol-alg}
\si_i^{fus}(u)=P_{\bT,n}\,\si_i^{(\bc)}(u)\,P_{\bT,n}\ \ \ \ \ \ \ \ \ \text{in $H^{fus}_{\bT,n}(q)$.}
\end{equation}
The function $\si_i^{fus}(u)$ satisfies the braided YB equation directly in $H^{fus}_{\bT,n}(q)$, since $\si_i^{(\bc)}(u)$ commutes with $P_{\bT,n}$ and already satisfies YB in $H_{kn}(q)$. We happily conclude this general study before plunging deep into an example.
\begin{conc} The fused Hecke algebra is:
\[H^{fus}_{\bT,n}(q):=P_{\bT,n}H_{kn}(q)P_{\bT,n}\ .\]
$\bullet$ It contains a solution $\si_i^{fus}(u)$ of the braided Yang--Baxter equation.\\
\noindent $\bullet$ It has a representation on $(L_{\lambda}^N)^{\otimes n}$ and the image of $\si_i^{fus}(u)$ recovers the fused solution $\check R_i^{fus}(u)$.
\end{conc}

\subsection{An example}\label{subsec-exa} 
Now we discuss the example which is behind the construction in the next section. Take the simplest partition $\lambda=(k)$ and the only standard Young tableau of shape $\lambda$ with its associated sequence of contents:
\[\bT=\fbox{1}\fbox{2}\dots \dots \fbox{k}\ \ \ \leadsto\ \ \bc=(1,q^2,\dots\,q^{2(k-1)})\ .\]
\paragraph{The idempotent $E_{\bT}$.} In this case, the idempotent $E_{\bT}$ obtained through the general Theorem \ref{thm-idem-Hn} is simply the $q$-symmetriser of $H_k(q)$:
\[
E_{\bT}=P_k=\frac{\sum_{w\in\mS_k}q^{\ell(w)}\si_w}{\sum_{w\in\mS_k}q^{2\ell(w)}}=\frac{q^{-k(k-1)/2}}{[k]_q!}\sum_{w\in\mS_k}q^{\ell(w)}\si_w\ .
\]
The two equalities are just different expressions for the normalisation factor. It might be profitable to see directly that $E_{\bT}=P_k$ in this case, since this is an easy particular case of Theorem \ref{thm-idem-Hn}.

Let us define the $q$-symmetriser $P_k$ directly by the property we actually want: it is the element of $H_k(q)$ projecting onto the representation $S_{(k)}$ in any representation of $H_k(q)$. In other words (from general considerations in representation theory) the $q$-symmetriser $P_k$ is characterised, up to a normalisation factor, by being a non-zero element satisfying either one of the two properties: 
\[\sigma_i P_k=q P_k\,,\ \ \forall i=1,\dots,k-1\,,\ \ \ \ \ \ \text{or}\ \ \ \ \ \ \ P_k\sigma_i=q P_k\,,\ \ \forall i=1,\dots,k-1\,.\]
It is rather easy to be convinced that these properties are satisfied with the defining formulas above for $P_k$ (and now it is clear that the normalisation is such that $P_k^2=P_k$). Let us see if this is also true for the following element $E_{\bT}$:
\[E_{\bT}=\Phi(1,q^2,\dots,q^{2(k-1)})=\prod_{i=1,\dots,k-1}^{\rightarrow} \si_i(q^{2i})\si_{i-1}(q^{2(i-1)})\dots\dots\si_1(q^2)\ ,\]
thereby verifying the particular case of Theorem \ref{thm-idem-Hn}. It is in fact not very difficult. For $k=1$ we have the easily verified fact: $\sigma_i(q^2)\sigma_1=q\sigma_1(q^2) $. Then we remark that $\sigma_i-q=\sigma_i(q^{-2})$ and, for $i>1$, in the product $E_{\bT}\sigma_i(q^{-2})$ we use the braided YB equation $\sigma_i(q^2u)\sigma_{i-1}(u)\sigma_i(q^{-2})=\sigma_{i-1}(q^{-2})\sigma_i(u)\sigma_{i-1}(q^2u)$ to move $\sigma_i(q^{-2})$ to the left. It is changed into $\sigma_{i-1}(q^{-2})$ and reached the idempotent corresponding to $k-1$ boxes, so we conclude using the induction hypothesis.

\paragraph{The idempotent $P_{\bT,n}$.} The idempotent $P_{\bT,n}$ of $H_{kn}(q)$ in this case is:
\[P_{\bT,n}=P_k\otimes\dots\otimes P_k\ ,\]
or said differently, it is the product of the $q$-symmetriser on $\si_1,\dots,\si_{k-1}$, with the $q$-symmetriser on $\si_{k+1},\dots,\si_{2k-1}$, and so on until the $q$-symmetriser on the $k-1$ last generators of $H_{kn}(q)$. An equivalent expression is:
\[P_{\bT,n}=\Bigl(\frac{q^{-k(k-1)/2}}{[k]_q!}\Bigr)^n\sum_{w\in\mS_k\times\dots\times \mS_k}q^{\ell(w)}\si_w\ .\]

\paragraph{The subspace $L_{(k)}^N$ (the $q$-symmetrised power $S_q^kV$).} We will call the irreducible representation $L_{(k)}^N$ of $U_q(sl_N)$ the $q$-symmetrised power of $V$ and denote it by $S_q^kV$, since from the Schur--Weyl duality, we have (this can be taken as the definition of $S_q^kV$):
\[S_q^kV=P_k(V^{\otimes k})\,.\]
To describe $S_q^kV$ explicitly, take $(e_1,\dots,e_N)$ the basis of $V$ that was used for defining the action of $H_k(q)$ in (\ref{rep-Hn}). The usual symmetrised power $S^k(V)$ (for $q=1$) has the following basis:
\begin{equation}\label{basis-sym}
e_{(i_1,\dots,i_k)}=\sum_{w\in\mS_k}e_{w(i_1)}\otimes\dots\otimes e_{w(i_k)}\,,\ \ \ \ \ \ \text{where $N\geq i_1\geq\dots\geq i_k\geq 1$\ .}
\end{equation}
For the $q$-symmetrised power $S^k_q(V)$, the basis is quite similar:
\begin{equation}\label{basis-qsym}
e^{(q)}_{(i_1,\dots,i_k)}=\sum_{w\in\mS_k}q^{\ell_{(i_1,\dots,i_k)}(w)}e_{w(i_1)}\otimes\dots\otimes e_{w(i_k)}\,,\ \ \ \ \ \ \text{where $N\geq i_1\geq\dots\geq i_k\geq 1$\ ,}
\end{equation}
where $\ell_{(i_1,\dots,i_n)}(w)$ is a modified length: it is the number of inversions of different indices $i_a,i_b$. More precisely, in formulas, it is given by:
\[\ell_{(i_1,\dots,i_n)}(w)=\text{Card}\{a<b\ \text{such that}\ w(a)>w(b)\ \text{and}\ i_a>i_b\}\ .\]
One can check that $e^{(q)}_{(i_1,\dots,i_k)}$ is simply $P_k(e_{i_1}\otimes\dots\otimes e_{i_k})$, and this is where the decreasing order on the indices $i_1,\dots,i_k$ is slightly more convenient.

In both formulas (\ref{basis-sym}) and (\ref{basis-qsym}), the sum does not have to be taken over the whole symmetric group $\mS_k$. It can be taken over the orbit of $(i_1,\dots,i_k)$ under the action of $\mS_k$ to avoid repetitions. That is to say, it can be taken over the symmetric group $\mS_k$ modulo the stabiliser of $(i_1,\dots,i_k)$.

As an example, take $k=3$ and $N=2$, then a basis of $S_q^k(V)$ consists of the vectors:
\[\begin{array}{l}
e_1\otimes e_1\otimes e_1\,,\\[0.4em]
e_2\otimes e_1\otimes e_1+qe_1\otimes e_2\otimes e_1+q^2e_1\otimes e_1\otimes e_2\,,\\[0.4em]
e_2\otimes e_2\otimes e_1+qe_2\otimes e_1\otimes e_2+q^2e_1\otimes e_2\otimes e_2\,,\\[0.4em]
e_2\otimes e_2\otimes e_2\ .
\end{array}\]

\newpage
\section{Fused Braids and the Fused Hecke Algebra}\label{sec-fusedH}

We fix an integer $k>0$. We took a long detour through quantum groups and the fusion procedure and we have reached a point where we know which algebras we are interested in. So we can forget for a while the preceding discussions, and only remember that we have an interest in the algebras:
\begin{equation}\label{def-Hfus}
H^{fus}_{k,n}(q):=P_{k,n}H_{kn}(q)P_{k,n}\,,
\end{equation}
that we call the fused Hecke algebras, and that are obtained from the construction of the preceding section applied to the $q$-symmetrised power $S_q^k(V)$ of the vector representation of $U_q(sl_N)$ (where $N=\dim(V)$). The idempotent $P_{k,n}$ and the space $S^k_q(V)$ were described explicitly in Section \ref{subsec-exa}.

The goal of this final section is to give a braid-like or diagrammatic description of $H^{fus}_{k,n}(q)$, which can be read independently of the preceding considerations. We leave to the readers to convince themselves that this is really the same algebra. All details can be found in \cite{CPdA1}. Finally, we will conclude that the fused Hecke algebra $H^{fus}_{k,n}(q)$ meets our expectations described in the introduction.

\subsection{Objects and multiplication}

\paragraph{Objects: fused braids.}  We consider the following objects, which are generalisations of usual braids. We keep the rectangular strip with a top line of $n$ fixed dots and a bottom line of $n$ fixed dots. And we put strands connecting dots but now we require the following: Each strand still connects a dot from the top line to a dot on the bottom line, but now, for each dot, there are $k$ strands attached to it. So in total there are $kn$ strands. 

To be more precise, it is better if we replace the dots by small ellipses, and we will do so from now on. Then the strands which are attached to the same ellipse are not really attached to the same point of the ellipse. Instead they are attached next to each other at the same ellipse. As before we require that at each point inside the strip at most two strands are crossing and we keep the same terminology of positive and negative crossings. Again as before we consider such diagrams up to isotopy, namely up to continuously moving the strands while leaving their end points fixed.

Such an object we call a fused braid. Needless to insist that for $k=1$, a fused braid is a usual braid. Some examples of fused braids can be found below.

\paragraph{The vector space.} We consider the vector space $Vect^{fus}_{k,n}$ of formal linear combinations of fused braids (in other words, a vector space with basis indexed by fused braids and we identify the basis vectors with their indices).
\begin{defi}\label{vector-fused-braids}
 The vector space $H^{fus}_{k,n}(q)$ is the quotient of $Vect^{fus}_{k,n}$ by the following relations:
 \begin{itemize}
   \item[(i)] The Hecke relation for all crossings:  
   \begin{center}
 \begin{tikzpicture}[scale=0.25]
\draw[thick] (0,2)..controls +(0,-2) and +(0,+2) .. (4,-2);
\fill[white] (2,0) circle (0.4);
\draw[thick] (4,2)..controls +(0,-2) and +(0,+2) .. (0,-2);
\node at (6,0) {$=$};
\draw[thick] (12,2)..controls +(0,-2) and +(0,+2) .. (8,-2);
\fill[white] (10,0) circle (0.4);
\draw[thick] (8,2)..controls +(0,-2) and +(0,+2) .. (12,-2);
\node at (17,0) {$-\,(q-q^{-1})$};
\draw[thick] (21,2) -- (21,-2);\draw[thick] (25,2) -- (25,-2);
\end{tikzpicture}
\end{center}
  \item[(ii)]  The idempotent relations for crossings near the ellipses:
 \begin{center}
 \begin{tikzpicture}[scale=0.4]
\fill (2,2) ellipse (0.8cm and 0.2cm);
%\draw[thick] (1.6,2)..controls +(0,-1) and +(1,1) .. (-1,0);
\draw[thick] (2.2,2)..controls +(0,-1.5) and +(1,1) .. (1.2,0);
\fill[white] (2,0.7) circle (0.2);
\draw[thick] (1.8,2)..controls +(0,-1.5) and +(-1,1) .. (2.8,0);
%\draw[thick] (2.4,2)..controls +(0,-1) and +(-1,1.4) .. (4.2,0);
\node at (4.5,1) {$=$};
\node at (7,1) {$q$};
\fill (9,2) ellipse (0.8cm and 0.2cm);
%\draw[thick] (8.6,2)..controls +(0,-1) and +(1,1) .. (6,0);
\draw[thick] (8.8,2)..controls +(0,-1.5) and +(0.5,0.5) .. (8.2,0);
\draw[thick] (9.2,2)..controls +(0,-1.5) and +(-0.5,0.5) .. (9.8,0);
%\draw[thick] (9.4,2)..controls +(0,-1) and +(-1,1.4) .. (11.2,0);

\node at (13,1) {and};

\fill (18,2) ellipse (0.8cm and 0.2cm);
%\draw[thick] (17.6,2)..controls +(0,-1) and +(1,1) .. (15,0);
\draw[thick] (17.8,2)..controls +(0,-1.5) and +(-1,1) .. (18.8,0);
\fill[white] (18,0.7) circle (0.2);
\draw[thick] (18.2,2)..controls +(0,-1.5) and +(1,1) .. (17.2,0);
%\draw[thick] (18.4,2)..controls +(0,-1) and +(-1,1.4) .. (20.2,0);
\node at (20.5,1) {$= $};
\node at (23,1) {$q^{-1}$};
\fill (25,2) ellipse (0.8cm and 0.2cm);
%\draw[thick] (24.6,2)..controls +(0,-1) and +(1,1) .. (22,0);
\draw[thick] (24.8,2)..controls +(0,-1.5) and +(0.5,0.5) .. (24.2,0);
\draw[thick] (25.2,2)..controls +(0,-1.5) and +(-0.5,0.5) .. (25.8,0);
%\draw[thick] (25.4,2)..controls +(0,-1) and +(-1,1.4) .. (27.2,0);
\end{tikzpicture}
\end{center}
\begin{center}
 \begin{tikzpicture}[scale=0.4]
\fill (2,0) ellipse (0.8cm and 0.2cm);
%\draw[thick] (1.6,0)..controls +(0,1) and +(1,-1.5) .. (-1,2);
\draw[thick] (1.8,0)..controls +(0,1.5) and +(-1,-1) .. (2.8,2);
\fill[white] (2,1.25) circle (0.2);
\draw[thick] (2.2,0)..controls +(0,1.5) and +(1,-1) .. (1.2,2);
%\draw[thick] (2.4,0)..controls +(0,1) and +(-1,-1.4) .. (4.2,2);
\node at (4.5,1) {$=$};
\node at (7,1) {$q$};
\fill (9,0) ellipse (0.8cm and 0.2cm);
%\draw[thick] (8.6,0)..controls +(0,1) and +(1,-1.5) .. (6,2);
\draw[thick] (9.2,0)..controls +(0,1.5) and +(-0.5,-0.5) .. (9.8,2);
\draw[thick] (8.8,0)..controls +(0,1.5) and +(0.5,-0.5) .. (8.2,2);
%\draw[thick] (9.4,0)..controls +(0,1) and +(-1,-1.4) .. (11.2,2);

\node at (13,1) {and};

\fill (18,0) ellipse (0.8cm and 0.2cm);
%\draw[thick] (17.6,0)..controls +(0,1) and +(1,-1.5) .. (15,2);
\draw[thick] (18.2,0)..controls +(0,1.5) and +(1,-1) .. (17.2,2);
\fill[white] (18,1.25) circle (0.2);
\draw[thick] (17.8,0)..controls +(0,1.5) and +(-1,-1) .. (18.8,2);
%\draw[thick] (18.4,0)..controls +(0,1) and +(-1,-1.4) .. (20.2,2);
\node at (20.5,1) {$=$};
\node at (23,1) {$q^{-1}$};
\fill (25,0) ellipse (0.8cm and 0.2cm);
%\draw[thick] (24.6,0)..controls +(0,1) and +(1,-1.5) .. (22,2);
\draw[thick] (25.2,0)..controls +(0,1.5) and +(-0.5,-0.5) .. (25.8,2);
\draw[thick] (24.8,0)..controls +(0,1.5) and +(0.5,-0.5) .. (24.2,2);
%\draw[thick] (25.4,0)..controls +(0,1) and +(-1,-1.4) .. (27.2,2);
\end{tikzpicture}
\end{center}
 \end{itemize}
\end{defi}
The first relation (the Hecke relation) is valid locally for all crossings as in the situation of classical braids and Hecke algebra. The idempotents relations are also local relations. In words, they impose the following: if two strands start from the same ellipse and their first crossing is crossing each other, then the original fused braid is equal to the fused braid obtained by removing this crossing and multiplying by $q^{\pm1}$ depending on the sign of the crossing; and similarly for two strands arriving at the same ellipse.

Note that the Hecke relation for a crossing near an ellipse is recovered as the difference of two idempotent relations, so that imposing these idepotent relations is quite compatible with the Hecke relation.

\paragraph{Multiplication.} Now we define a product on the vector space $H^{fus}_{k,n}(q)$, which makes it an associative unital algebra. Let $b,b'$ be two fused braids. We define $bb'$ as the result of the following procedure:
\begin{itemize}
\item \emph{(Concatenation)} We place the diagram of $b$ on top of the diagram of $b'$ by identifying the bottom ellipses of $b$ with the top ellipses of $b'$
\item \emph{(Removal of middle ellipses)} For each ellipse in the middle row, there are $k$ strands incoming and $k$ strands leaving. We remove this middle ellipse and replace it by the $q$-symmetriser $P_{k}$ of the Hecke algebra $H_k(q)$.
\end{itemize}
More explicitly, recall again the formula for the $q$-symmetriser in $H_k(q)$:
\[
P_k=\frac{q^{-k(k-1)/2}}{[k]_q!}\sum_{w\in\mS_k}q^{\ell(w)}\si_w\ .
\]
So in order to remove a middle ellipse, we take $w\in\mS_{k}$ and we first construct the diagram where this middle ellipse is replaced by the element $\si_w$ connecting the $k$ incoming strands with the $k$ outgoing ones. Then we make the sum over $w\in\mS_{k}$ of the resulting diagrams, each with the coefficient $q^{\ell(w)}$ (that is, multiplied by $q$ to the power the number of crossings we just added). Finally, we normalise as in the formula above. 

\vskip .2cm
The fused braid with only non-crossing vertical strands is the unit element of $H^{fus}_{k,n}(q)$. For example, for $k=2$, it is:
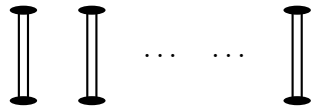
\begin{center}
\begin{tikzpicture}[scale=0.3]
\fill (1,2) ellipse (0.6cm and 0.2cm);\fill (1,-2) ellipse (0.6cm and 0.2cm);
\draw[thick] (0.8,2) -- (0.8,-2);\draw[thick] (1.2,2) -- (1.2,-2);
\fill (4,2) ellipse (0.6cm and 0.2cm);\fill (4,-2) ellipse (0.6cm and 0.2cm);
\draw[thick] (3.8,2) -- (3.8,-2);\draw[thick] (4.2,2) -- (4.2,-2);
\node at (7,0) {$\ldots$};
\node at (10,0) {$\ldots$};
\fill (13,2) ellipse (0.6cm and 0.2cm);\fill (13,-2) ellipse (0.6cm and 0.2cm);
\draw[thick] (12.8,2) -- (12.8,-2);\draw[thick] (13.2,2) -- (13.2,-2);
\end{tikzpicture}
\end{center}
The reason is the following. Say we plug this element on top of another fused braid, and we replace a middle ellipse by a certain $\si_w$ in $H_k(q)$. Then the whole lot of crossings added by $\si_w$ can be moved up along the parallel strands, and will thus hit an ellipse of the top row. Thus we can use the idempotent relation. Therefore, we see that replacing a middle ellipse by $\si_w$ will be the same as replacing the middle ellipse by the trivial braid, except for a factor $q^{\ell(w)}$. It remains only to perform the (normalised) sum and we conclude that multiplying by the above element does nothing.

\begin{exam}\label{ex:rem} We illustrate below the procedure to remove a middle ellipse:

$\bullet$ For $k=2$, two strands are incoming at each middle ellipse and two strands are leaving. There are two terms in the sum for the $q$-symmetriser: $\frac{1}{1+q^2}(1+q\si_1)$. So each ellipse is replaced by a sum of two terms:
 
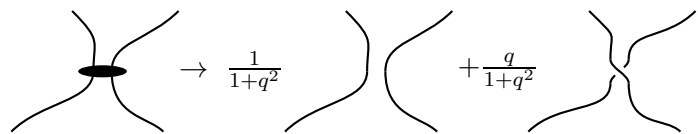
\begin{center}
 \begin{tikzpicture}[scale=0.4]
\fill (0,0) ellipse (0.8cm and 0.2cm);
%\node at (-1,1) {$\dots$};
\draw[thick] (-0.3,0)..controls +(0,1) and +(1,-1) .. (-1,2);
\draw[thick] (0.3,0)..controls +(0,1) and +(-1,-1) .. (2.5,2);
\draw[thick] (-0.3,0)..controls +(0,-1) and +(1,1) .. (-3,-2);
\draw[thick] (0.3,0)..controls +(0,-1.5) and +(-0.5,0.5) .. (2,-2);
\node at (3,0) {$\rightarrow$};
\node at (5,0) {$\frac{1}{1+q^2}$};
\draw[thick] (8.7,0)..controls +(0,1) and +(1,-1) .. (8,2);
\draw[thick] (9.3,0)..controls +(0,1) and +(-1,-1) .. (11.5,2);
\draw[thick] (8.7,0)..controls +(0,-1) and +(1,1) .. (6,-2);
\draw[thick] (9.3,0)..controls +(0,-1.5) and +(-0.5,0.5) .. (11,-2);
\node at (13,0) {$+\frac{q}{1+q^2}$};
\draw[thick] (16.7,0.4)..controls +(0,1) and +(1,-1) .. (16,2);
\draw[thick] (17.3,0.4)..controls +(0,1) and +(-1,-1) .. (19.5,2);
\draw[thick] (17.3,0.4)..controls +(0,-0.3) and +(0,0.3) .. (16.7,-0.4);
\fill[white] (17,0) circle (0.2);
\draw[thick] (16.7,0.4)..controls +(0,-0.3) and +(0,0.3) .. (17.3,-0.4);
\draw[thick] (16.7,-0.4)..controls +(0,-1) and +(1,1) .. (14,-2);
\draw[thick] (17.3,-0.4)..controls +(0,-1.5) and +(-0.5,0.5) .. (19,-2);
\end{tikzpicture}
\end{center}
Note then that if we are making a product in $H_{2,n}^{fus}(q)$, there are $n$ middle ellipses, and thus we have a sum of $2n$ terms.
\vskip .2cm
$\bullet$ For $k=3$, each middle ellipse has three strands arriving and leaving:  \begin{tikzpicture}[scale=0.3]
\fill (0,0) ellipse (0.8cm and 0.2cm);
\draw[thick] (-0.3,0)--(-0.5,1);
\draw[thick] (0,0)--(0,1);
\draw[thick] (0.3,0)--(0.5,1);
\draw[thick] (-0.3,0)--(-0.5,-1);
\draw[thick] (0,0)--(0,-1);
\draw[thick] (0.3,0)--(0.5,-1);
\end{tikzpicture}; each one is replaced by the following sum divided by $1+2q^2+2q^4+q^6=(1+q^2)(1+q^2+q^4)$:
\begin{center}
 \begin{tikzpicture}[scale=0.3]
\fill (1,12) circle (0.2cm);\fill (1,8) circle (0.2cm);
\draw[thick] (1,12) -- (1,8);
\fill (3,12) circle (0.2cm);\fill (3,8) circle (0.2cm);
\draw[thick] (3,12) -- (3,8);
\fill (5,12) circle (0.2cm);\fill (5,8) circle (0.2cm);
\draw[thick] (5,12) -- (5,8);
\node at (7,10) {$+$};
\node at (9,10) {$q$};
\fill (10,12) circle (0.2cm);\fill (10,8) circle (0.2cm);
\draw[thick] (12,12)..controls +(0,-2) and +(0,+2) .. (10,8);\fill[white] (11,10) circle (0.4);
\draw[thick] (10,12)..controls +(0,-2) and +(0,+2) .. (12,8);
\fill (12,12) circle (0.2cm);\fill (12,8) circle (0.2cm);
\fill (14,12) circle (0.2cm);\fill (14,8) circle (0.2cm);
\draw[thick] (14,12) -- (14,8);
\node at (16,10) {$+$};
\node at (18,10) {$q$};
\fill (19,12) circle (0.2cm);\fill (19,8) circle (0.2cm);
\draw[thick] (19,12) -- (19,8);
\fill (21,12) circle (0.2cm);\fill (21,8) circle (0.2cm);
\draw[thick] (23,12)..controls +(0,-2) and +(0,+2) .. (21,8);\fill[white] (22,10) circle (0.4);
\fill (23,12) circle (0.2cm);\fill (23,8) circle (0.2cm);
\draw[thick] (21,12)..controls +(0,-2) and +(0,+2) .. (23,8);
\node at (25,10) {$+$};
\node at (27,10) {$q^{2}$};
\draw[thick] (32,12)..controls +(0,-2) and +(0,+2) .. (28,8);\fill[white] (30.6,10.4) circle (0.4);\fill[white] (29.4,9.6) circle (0.4);
\fill (28,12) circle (0.2cm);\fill (28,8) circle (0.2cm);
\draw[thick] (28,12)..controls +(0,-2) and +(0,+2) .. (30,8);
\fill (30,12) circle (0.2cm);\fill (30,8) circle (0.2cm);
\draw[thick] (30,12)..controls +(0,-2) and +(0,+2) .. (32,8);
\fill (32,12) circle (0.2cm);\fill (32,8) circle (0.2cm);
\node at (34,10) {$+$};
\node at (36,10) {$q^{2}$};
\fill (37,12) circle (0.2cm);\fill (37,8) circle (0.2cm);
\draw[thick] (41,12)..controls +(0,-2) and +(0,+2) .. (39,8);
\draw[thick] (39,12)..controls +(0,-2) and +(0,+2) .. (37,8);
\fill[white] (38.4,10.3) circle (0.4);\fill[white] (39.5,9.7) circle (0.4);
\draw[thick] (37,12)..controls +(0,-2) and +(0,+2) .. (41,8);
\fill (39,12) circle (0.2cm);\fill (39,8) circle (0.2cm);
\fill (41,12) circle (0.2cm);\fill (41,8) circle (0.2cm);
\node at (43,10) {$+$};
\node at (45,10) {$q^{3}$};
\fill (46,12) circle (0.2cm);\fill (46,8) circle (0.2cm);
\fill (48,12) circle (0.2cm);\fill (48,8) circle (0.2cm);
\fill (50,12) circle (0.2cm);\fill (50,8) circle (0.2cm);
\draw[thick] (50,12)..controls +(0,-3) and +(0,+1) .. (46,8);
\fill[white] (48,9.2) circle (0.3);
\draw[thick] (48,12) -- (48,8);
\fill[white] (48,10.8) circle (0.3);\fill[white] (49.1,10) circle (0.3);
\draw[thick] (46,12)..controls +(0,-1) and +(0,+3) .. (50,8);
\end{tikzpicture}
\end{center}

$\bullet$ Here is an example of a product of two elements of $H^{fus}_{2,2}(q)$:
\begin{center}
\begin{tikzpicture}[scale=0.3]
\fill (1,2) ellipse (0.6cm and 0.2cm);\fill (1,-2) ellipse (0.6cm and 0.2cm);
\fill (4,2) ellipse (0.6cm and 0.2cm);\fill (4,-2) ellipse (0.6cm and 0.2cm);
\draw[thick] (3.8,2)..controls +(0,-2) and +(0,+2) .. (1.2,-2); \draw[thick] (4.2,2) -- (4.2,-2);
\fill[white] (2.5,0) circle (0.4);
\draw[thick] (0.8,2) -- (0.8,-2);\draw[thick] (1.2,2)..controls +(0,-2) and +(0,+2) .. (3.8,-2);  

\node at (6,0) {$.$};

\fill (8,2) ellipse (0.6cm and 0.2cm);\fill (8,-2) ellipse (0.6cm and 0.2cm);
\fill (11,2) ellipse (0.6cm and 0.2cm);\fill (11,-2) ellipse (0.6cm and 0.2cm);
\draw[thick] (10.8,2)..controls +(0,-2) and +(0,+2) .. (8.2,-2); \draw[thick] (11.2,2) -- (11.2,-2);
\fill[white] (9.5,0) circle (0.4);
\draw[thick] (7.8,2) -- (7.8,-2);\draw[thick] (8.2,2)..controls +(0,-2) and +(0,+2) .. (10.8,-2);  

\node at (17.5,0) {$=\displaystyle\frac{1}{(1+q^2)^2}\Bigg($};
\fill (22.5,4) ellipse (0.6cm and 0.2cm);
\fill (25.5,4) ellipse (0.6cm and 0.2cm);
\draw[thick] (25.3,4)..controls +(0,-2) and +(0,+2) .. (22.7,0); \draw[thick] (25.7,4) -- (25.7,0);
\fill[white] (24,2) circle (0.4);
\draw[thick] (22.3,4) -- (22.3,0);\draw[thick] (22.7,4)..controls +(0,-2) and +(0,+2) .. (25.3,0);  

\fill (22.5,-4) ellipse (0.6cm and 0.2cm);
\fill (25.5,-4) ellipse (0.6cm and 0.2cm);
\draw[thick] (25.3,0)..controls +(0,-2) and +(0,+2) .. (22.7,-4); \draw[thick] (25.7,0) -- (25.7,-4);
\fill[white] (24,-2) circle (0.4);
\draw[thick] (22.3,0) -- (22.3,-4);\draw[thick] (22.7,0)..controls +(0,-2) and +(0,+2) .. (25.3,-4);  

\node at (27.5,0) {$+ q$};
\fill (29.5,4) ellipse (0.6cm and 0.2cm);
\fill (32.5,4) ellipse (0.6cm and 0.2cm);
\draw[thick] (32.3,4)..controls +(0,-2) and +(0,+2) .. (29.7,0.5); \draw[thick] (32.7,4) -- (32.7,0);
\fill[white] (30.9,2) circle (0.4);
\draw[thick] (29.3,4) -- (29.3,0.5);\draw[thick] (29.7,4)..controls +(0,-2) and +(0,+2) .. (32.3,0);  

\draw[thick] (29.7,0.5)..controls +(0,-0.2) and +(0,+0.2) .. (29.3,-0.5);
\fill[white] (29.5,0) circle (0.2);
\draw[thick] (29.3,0.5)..controls +(0,-0.2) and +(0,+0.2) .. (29.7,-0.5);

\fill (29.5,-4) ellipse (0.6cm and 0.2cm);
\fill (32.5,-4) ellipse (0.6cm and 0.2cm);
\draw[thick] (32.3,0)..controls +(0,-2) and +(0,+2) .. (29.7,-4); \draw[thick] (32.7,0) -- (32.7,-4);
\fill[white] (30.8,-2.1) circle (0.4);
\draw[thick] (29.3,-0.5) -- (29.3,-4);\draw[thick] (29.7,-0.5)..controls +(0,-2) and +(0,+2) .. (32.3,-4);

\node at (34.5,0) {$+q$};
\fill (36.5,4) ellipse (0.6cm and 0.2cm);
\fill (39.5,4) ellipse (0.6cm and 0.2cm);
\draw[thick] (39.3,4)..controls +(0,-2) and +(0,+2) .. (36.7,0); \draw[thick] (39.7,4) -- (39.7,0.5);
\fill[white] (38.2,2.1) circle (0.4);
\draw[thick] (36.3,4) -- (36.3,0);\draw[thick] (36.7,4)..controls +(0,-2) and +(0,+2) .. (39.3,0.5);  

\draw[thick] (39.7,0.5)..controls +(0,-0.2) and +(0,+0.2) .. (39.3,-0.5);
\fill[white] (39.5,0) circle (0.2);
\draw[thick] (39.3,0.5)..controls +(0,-0.2) and +(0,+0.2) .. (39.7,-0.5);

\fill (36.5,-4) ellipse (0.6cm and 0.2cm);
\fill (39.5,-4) ellipse (0.6cm and 0.2cm);
\draw[thick] (39.3,-0.5)..controls +(0,-2) and +(0,+2) .. (36.7,-4); \draw[thick] (39.7,-0.5) -- (39.7,-4);
\fill[white] (38,-2) circle (0.4);
\draw[thick] (36.3,0) -- (36.3,-4);\draw[thick] (36.7,0)..controls +(0,-2) and +(0,+2) .. (39.3,-4);  

\node at (41.5,0) {$+q^2$};
\fill (43.5,4) ellipse (0.6cm and 0.2cm);
\fill (46.5,4) ellipse (0.6cm and 0.2cm);
\draw[thick] (46.3,4)..controls +(0,-2) and +(0,+2) .. (43.7,0.5); \draw[thick] (46.7,4) -- (46.7,0.5);
\fill[white] (45.1,2) circle (0.4);
\draw[thick] (43.3,4) -- (43.3,0.5);\draw[thick] (43.7,4)..controls +(0,-2) and +(0,+2) .. (46.3,0.5);  

\draw[thick] (43.7,0.5)..controls +(0,-0.2) and +(0,+0.2) .. (43.3,-0.5);
\fill[white] (43.5,0) circle (0.2);
\draw[thick] (43.3,0.5)..controls +(0,-0.2) and +(0,+0.2) .. (43.7,-0.5);

\draw[thick] (46.7,0.5)..controls +(0,-0.2) and +(0,+0.2) .. (46.3,-0.5);
\fill[white] (46.5,0) circle (0.2);
\draw[thick] (46.3,0.5)..controls +(0,-0.2) and +(0,+0.2) .. (46.7,-0.5);

\fill (43.5,-4) ellipse (0.6cm and 0.2cm);
\fill (46.5,-4) ellipse (0.6cm and 0.2cm);
\draw[thick] (46.3,-0.5)..controls +(0,-2) and +(0,+2) .. (43.7,-4); \draw[thick] (46.7,-0.5) -- (46.7,-4);
\fill[white] (44.9,-2) circle (0.4);
\draw[thick] (43.3,-0.5) -- (43.3,-4);\draw[thick] (43.7,-0.5)..controls +(0,-2) and +(0,+2) .. (46.3,-4);  

\node at (48.5,0) {$\Bigg)$};

\node at (17.5,-8) {$=\displaystyle\frac{1}{(1+q^2)^2}\Bigg($};
\fill (22,-6) ellipse (0.6cm and 0.2cm);\fill (22,-10) ellipse (0.6cm and 0.2cm);
\draw[thick] (21.8,-6) -- (21.8,-10);\draw[thick] (22.2,-6) -- (22.2,-10);
\fill (25,-6) ellipse (0.6cm and 0.2cm);\fill (25,-10) ellipse (0.6cm and 0.2cm);
\draw[thick] (24.8,-6) -- (24.8,-10);\draw[thick] (25.2,-6) -- (25.2,-10);

\node at (31,-8) {$+(q-q^{-1}+2q^3)$};

\fill (37,-6) ellipse (0.6cm and 0.2cm);\fill (37,-10) ellipse (0.6cm and 0.2cm);
\fill (40,-6) ellipse (0.6cm and 0.2cm);\fill (40,-10) ellipse (0.6cm and 0.2cm);
\draw[thick] (39.8,-6)..controls +(0,-2) and +(0,+2) .. (37.2,-10);  \draw[thick] (40.2,-6) -- (40.2,-10);
\fill[white] (38.5,-8) circle (0.4);
\draw[thick] (36.8,-6) -- (36.8,-10);\draw[thick] (37.2,-6)..controls +(0,-2) and +(0,+2) .. (39.8,-10);  

\node at (42,-8) {$+q^2$};

\fill (44,-6) ellipse (0.6cm and 0.2cm);\fill (44,-10) ellipse (0.6cm and 0.2cm);
\fill (47,-6) ellipse (0.6cm and 0.2cm);\fill (47,-10) ellipse (0.6cm and 0.2cm);
\draw[thick] (46.8,-6)..controls +(0,-2) and +(0,+2) .. (43.8,-10);  
\draw[thick] (47.2,-6)..controls +(0,-2) and +(0,+2) .. (44.2,-10);  
\fill[white] (45.5,-8) circle (0.4);
\draw[thick] (43.8,-6)..controls +(0,-2) and +(0,+2) .. (46.8,-10);  
\draw[thick] (44.2,-6)..controls +(0,-2) and +(0,+2) .. (47.2,-10);  

\node at (49,-8) {$\Bigg)$};

\end{tikzpicture}
\end{center}

For the four diagrams obtained after the first equality, we can proceed as follows. For the first diagram, we apply the Hecke relation and this results to the identity term and the term with coefficient $(q-q^{-1})$. For the second diagram, we take the strand connecting the first top ellipse to the first bottom ellipse, and we move it on the left of the diagram, then we apply the idempotent relations; this results with one term with coefficient $q^3$. We proceed similarly for the third diagram. We do almost nothing except moving the strands in the fourth diagram.
\end{exam}

\subsection{Fused permutations and standard basis}

The fused Hecke algebra involves ``topological objects'' having strands which can cross by passing over or below each other. For $q=1$, all topological information disappears and we obtain an algebra defined purely combinatorially. This is another example of a classical situation as shown in the following table (unfortunately, we have not discussed the Brauer algebra; it is a $q=1$ limit of the BMW algebra).
\begin{center}
\begin{tabular}{c|c}
Topological Algebras & Combinatorial Algebras\\
\hline & \\[-0.2em]
Hecke algebra & Symmetric group\\[0.5em]
BMW algebra & Brauer algebra\\[0.5em]
Fused Hecke algebra & Algebra of fused permutations
\end{tabular}
\end{center}
The connections between the two columns of the table is that the algebras on the left are flat deformations of the algebras on the right. Here ``flat'' means that the vector space remains the same while the multiplication is modified (with a parameter $q$). One interest is that we obtain a natural, or standard, basis for algebras on the left indexed by combinatorial objects (coming from the algebras on the right). Let us see how it works for the fused Hecke algebras.

\paragraph{Fused permutations.} We give a diagrammatic definition of fused permutations. On one hand, they generalise obviously the usual permutations of $\mS_n$ (which correspond to $k=1$). On the other hand, they are the combinatorial shadow of fused braids, meaning that they are obtained from fused braids if we forget all topological information (namely, the dots are not connected anymore by strands from a 3-dimensional world).

So, we place two horizontal rows of $n$ dots, one on top of another. And we connect the dots of the top row to the dots on the bottom row. We require the following: there are $k$ edges which start from each dot on the top row and there are $k$ edges which arrive at each dot on the bottom row. Again, as for fused braids, it is better to see the dots as small ellipses and edges are attached to nearby points on the ellipses.

The only information that matters is the following: which ellipse is connected to which other ellipse and by how many edges. More rigorously, take a diagram as above and let $a\in\{1,\dots,n\}$. There are $k$ edges starting from the $a$-th ellipse on top and we denote by $I_a$ the multiset $\{i_1,\dots,i_{k}\}$ indicating the bottom ellipses reached by these edges. This is indeed a multiset, meaning that repetitions are allowed since several of these edges can reach the same bottom ellipse. We consider two diagrams equivalent if their sequences of multisets $(I_1,\dots,I_n)$ coincide.

\begin{defi}
A fused permutation is an equivalence class of diagrams as explained above. We denote by $\mS^{fus}_{k,n}$ the set of fused permutations.
\end{defi}

\begin{rema} The set of fused permutations $\mS^{fus}_{k,n}$ is in bijection with:

$\bullet$ (by definition) the set of sequences $(I_1,\dots,I_n)$ of multisets such that:
\[|I_1|=\dots=|I_n|=k\ \ \ \quad\text{and}\quad\ \ \ I_1\cup\dots\cup I_n=\{\underbrace{1,\dots,1}_{k},\dots\dots,\underbrace{n,\dots,n}_{k}\}\,,\]
where the union is understood as concatenation of multisets;

$\bullet$ the set of $n$ by $n$ matrices with entries in $\mathbb{Z}_{\geq0}$ such that the sum of the entries in each row and each column is equal to $k$. The bijection is such that the entry in position $(a,b)$ indicates how many times the $a$-th top ellipse is connected to the $b$-th bottom ellipse;

$\bullet$ the set of double cosets in $\mS_{kn}$ modulo the subgroup $\mS_k\times\dots\times \mS_k$ ($n$ factors).
\end{rema}

\begin{exam}
$\bullet$ Let $n=2$ and  take $k=2$. We show below the three elements of $\mS_{2,2}^{fus}$. For each element, we give a diagram, the corresponding sequence of multisets, the corresponding matrix, and the corresponding double coset modulo the subgroup $H=\mS_2\times\mS_2$:

\begin{center}
\begin{tikzpicture}[scale=0.3]
\fill (1,2) ellipse (0.6cm and 0.2cm);\fill (1,-2) ellipse (0.6cm and 0.2cm);
\draw[thick] (0.8,2)..controls +(0,-2) and +(0,+2) .. (0.8,-2);\draw[thick] (1.2,2)..controls +(0,-2) and +(0,+2) .. (1.2,-2);  
\fill (4,2) ellipse (0.6cm and 0.2cm);\fill (4,-2) ellipse (0.6cm and 0.2cm);
\draw[thick] (3.8,2)..controls +(0,-2) and +(0,+2) .. (3.8,-2);\draw[thick] (4.2,2)..controls +(0,-2) and +(0,+2) .. (4.2,-2);
\node at (14,0) {$(\{1,1\},\{2,2\})$};
\node at (26,0) {$\left(\begin{array}{cc}2 &0\\0 & 2\end{array}\right)$};
\node at (38,0) {$H\phantom{(2,3)(1,4)H}$};
\end{tikzpicture}
\end{center}
\begin{center}
\begin{tikzpicture}[scale=0.3]
\fill (1,2) ellipse (0.6cm and 0.2cm);\fill (1,-2) ellipse (0.6cm and 0.2cm);
\draw[thick] (0.8,2) -- (0.8,-2);\draw[thick] (1.2,2)..controls +(0,-2) and +(0,+2) .. (3.8,-2);  
\fill (4,2) ellipse (0.6cm and 0.2cm);\fill (4,-2) ellipse (0.6cm and 0.2cm);
\draw[thick] (3.8,2)..controls +(0,-2) and +(0,+2) .. (1.2,-2); \draw[thick] (4.2,2)..controls +(0,-2) and +(0,+2) .. (4.2,-2);
\node at (14,0) {$(\{1,2\},\{1,2\})$};
\node at (26,0) {$\left(\begin{array}{cc}1 &1\\1 & 1\end{array}\right)$};
\node at (38,0) {$H(2,3)H\phantom{(1,4)}$};
\end{tikzpicture}
\end{center}
\begin{center}
\begin{tikzpicture}[scale=0.3]
\fill (1,2) ellipse (0.6cm and 0.2cm);\fill (1,-2) ellipse (0.6cm and 0.2cm);
\draw[thick] (0.8,2)..controls +(0,-2.1) and +(0,+1.9) .. (3.8,-2);\draw[thick] (1.2,2)..controls +(0,-1.9) and +(0,+2.1) .. (4.2,-2);
\fill (4,2) ellipse (0.6cm and 0.2cm);\fill (4,-2) ellipse (0.6cm and 0.2cm);
\draw[thick] (3.8,2)..controls +(0,-1.9) and +(0,+2.1) .. (0.8,-2);\draw[thick] (4.2,2)..controls +(0,-2) and +(0,+2) .. (1.2,-2); 
\node at (14,0) {$(\{2,2\},\{1,1\})$};
\node at (26,0) {$\left(\begin{array}{cc}0 &2\\2 & 0\end{array}\right)$};
\node at (38,0) {$H(1,3)(2,4)H$};
\end{tikzpicture}
\end{center}
\end{exam}

\paragraph{Algebra of fused permutations.} We will denote at once by $H^{fus}_{k,n}(1)$ the algebra of fused permutations, since it will become clear that this is the fused Hecke algebra for $q=1$. The algebra $H^{fus}_{k,n}(1)$ is the vector space with basis indexed by $\mS^{fus}_{k,n}$, and with the multiplication given as follows. Let $d,d'\in \mS^{fus}_{k,n}$.
\begin{itemize}
\item \emph{(Concatenation)} We place the diagram of $d$ on top of the diagram of $d'$ by identifying the bottom line of ellipses of $d$ with the top line of ellipses of $d'$.
\item \emph{(Removal of middle ellipses)} For each ellipse in the middle row, there are $k$ edges incoming and $k$ edges leaving. We delete this ellipse and sum over all possibilities of connecting the $k$ incoming edges with the $k$ leaving edges (there are $k!$ possibilities).
 \item  \emph{(Normalisation)} We divide the resulting sum by $(k!)^n$.
\end{itemize}
At the end of the procedure described above, we obtain a sum of diagrams representing a sum of fused permutations (with rational coefficients). This is what we define to be $dd'$ in $H^{fus}_{k,n}(1)$.

\begin{exam}$\ $

$\bullet$ Of course $H^{fus}_{1,n}(1)$ coincides with the group algebra of the symmetric group $\mS_n$.

$\bullet$ Here is an example of a product of two elements of $H^{fus}_{2,2}(1)$:
\begin{center}
\begin{tikzpicture}[scale=0.3]
\fill (1,2) ellipse (0.6cm and 0.2cm);\fill (1,-2) ellipse (0.6cm and 0.2cm);
\draw[thick] (0.8,2) -- (0.8,-2);\draw[thick] (1.2,2)..controls +(0,-2) and +(0,+2) .. (3.8,-2);  
\fill (4,2) ellipse (0.6cm and 0.2cm);\fill (4,-2) ellipse (0.6cm and 0.2cm);
\draw[thick] (3.8,2)..controls +(0,-2) and +(0,+2) .. (1.2,-2); \draw[thick] (4.2,2) -- (4.2,-2);
\node at (6,0) {$.$};
\fill (8,2) ellipse (0.6cm and 0.2cm);\fill (8,-2) ellipse (0.6cm and 0.2cm);
\draw[thick] (7.8,2) -- (7.8,-2);\draw[thick] (8.2,2)..controls +(0,-2) and +(0,+2) .. (10.8,-2);  
\fill (11,2) ellipse (0.6cm and 0.2cm);\fill (11,-2) ellipse (0.6cm and 0.2cm);
\draw[thick] (10.8,2)..controls +(0,-2) and +(0,+2) .. (8.2,-2); \draw[thick] (11.2,2) -- (11.2,-2);
\node at (13,0) {$=$};
\fill (15,4) ellipse (0.6cm and 0.2cm);\fill (15,0) ellipse (0.6cm and 0.2cm);
\draw[thick] (14.8,4) -- (14.8,0);\draw[thick] (15.2,4)..controls +(0,-2) and +(0,+2) .. (17.8,0);  
\fill (18,4) ellipse (0.6cm and 0.2cm);\fill (18,0) ellipse (0.6cm and 0.2cm);
\draw[thick] (17.8,4)..controls +(0,-2) and +(0,+2) .. (15.2,0); \draw[thick] (18.2,4) -- (18.2,0);

\fill (15,0) ellipse (0.6cm and 0.2cm);\fill (15,-4) ellipse (0.6cm and 0.2cm);
\draw[thick] (14.8,0) -- (14.8,-4);\draw[thick] (15.2,0)..controls +(0,-2) and +(0,+2) .. (17.8,-4);  
\fill (18,0) ellipse (0.6cm and 0.2cm);\fill (18,-4) ellipse (0.6cm and 0.2cm);
\draw[thick] (17.8,0)..controls +(0,-2) and +(0,+2) .. (15.2,-4); \draw[thick] (18.2,0) -- (18.2,-4);

\node at (20.5,0) {$=\frac{1}{4}\Bigl($};
\fill (22.5,4) ellipse (0.6cm and 0.2cm);
\draw[thick] (22.3,4) -- (22.3,0);\draw[thick] (22.7,4)..controls +(0,-2) and +(0,+2) .. (25.3,0);  
\fill (25.5,4) ellipse (0.6cm and 0.2cm);
\draw[thick] (25.3,4)..controls +(0,-2) and +(0,+2) .. (22.7,0); \draw[thick] (25.7,4) -- (25.7,0);

\fill (22.5,-4) ellipse (0.6cm and 0.2cm);
\draw[thick] (22.3,0) -- (22.3,-4);\draw[thick] (22.7,0)..controls +(0,-2) and +(0,+2) .. (25.3,-4);  
\fill (25.5,-4) ellipse (0.6cm and 0.2cm);
\draw[thick] (25.3,0)..controls +(0,-2) and +(0,+2) .. (22.7,-4); \draw[thick] (25.7,0) -- (25.7,-4);

\node at (27.5,0) {$+$};
\fill (29.5,4) ellipse (0.6cm and 0.2cm);
\draw[thick] (29.3,4) -- (29.3,0.5);\draw[thick] (29.7,4)..controls +(0,-2) and +(0,+2) .. (32.3,0);  
\fill (32.5,4) ellipse (0.6cm and 0.2cm);
\draw[thick] (32.3,4)..controls +(0,-2) and +(0,+2) .. (29.7,0.5); \draw[thick] (32.7,4) -- (32.7,0);

\draw[thick] (29.3,0.5)..controls +(0,-0.2) and +(0,+0.2) .. (29.7,-0.5);
\draw[thick] (29.7,0.5)..controls +(0,-0.2) and +(0,+0.2) .. (29.3,-0.5);

\fill (29.5,-4) ellipse (0.6cm and 0.2cm);
\draw[thick] (29.3,-0.5) -- (29.3,-4);\draw[thick] (29.7,-0.5)..controls +(0,-2) and +(0,+2) .. (32.3,-4);  
\fill (32.5,-4) ellipse (0.6cm and 0.2cm);
\draw[thick] (32.3,0)..controls +(0,-2) and +(0,+2) .. (29.7,-4); \draw[thick] (32.7,0) -- (32.7,-4);

\node at (34.5,0) {$+$};
\fill (36.5,4) ellipse (0.6cm and 0.2cm);
\draw[thick] (36.3,4) -- (36.3,0);\draw[thick] (36.7,4)..controls +(0,-2) and +(0,+2) .. (39.3,0.5);  
\fill (39.5,4) ellipse (0.6cm and 0.2cm);
\draw[thick] (39.3,4)..controls +(0,-2) and +(0,+2) .. (36.7,0); \draw[thick] (39.7,4) -- (39.7,0.5);

\draw[thick] (39.3,0.5)..controls +(0,-0.2) and +(0,+0.2) .. (39.7,-0.5);
\draw[thick] (39.7,0.5)..controls +(0,-0.2) and +(0,+0.2) .. (39.3,-0.5);

\fill (36.5,-4) ellipse (0.6cm and 0.2cm);
\draw[thick] (36.3,0) -- (36.3,-4);\draw[thick] (36.7,0)..controls +(0,-2) and +(0,+2) .. (39.3,-4);  
\fill (39.5,-4) ellipse (0.6cm and 0.2cm);
\draw[thick] (39.3,-0.5)..controls +(0,-2) and +(0,+2) .. (36.7,-4); \draw[thick] (39.7,-0.5) -- (39.7,-4);

\node at (41.5,0) {$+$};
\fill (43.5,4) ellipse (0.6cm and 0.2cm);
\draw[thick] (43.3,4) -- (43.3,0.5);\draw[thick] (43.7,4)..controls +(0,-2) and +(0,+2) .. (46.3,0.5);  
\fill (46.5,4) ellipse (0.6cm and 0.2cm);
\draw[thick] (46.3,4)..controls +(0,-2) and +(0,+2) .. (43.7,0.5); \draw[thick] (46.7,4) -- (46.7,0.5);

\draw[thick] (43.3,0.5)..controls +(0,-0.2) and +(0,+0.2) .. (43.7,-0.5);
\draw[thick] (43.7,0.5)..controls +(0,-0.2) and +(0,+0.2) .. (43.3,-0.5);
\draw[thick] (46.3,0.5)..controls +(0,-0.2) and +(0,+0.2) .. (46.7,-0.5);
\draw[thick] (46.7,0.5)..controls +(0,-0.2) and +(0,+0.2) .. (46.3,-0.5);

\fill (43.5,-4) ellipse (0.6cm and 0.2cm);
\draw[thick] (43.3,-0.5) -- (43.3,-4);\draw[thick] (43.7,-0.5)..controls +(0,-2) and +(0,+2) .. (46.3,-4);  
\fill (46.5,-4) ellipse (0.6cm and 0.2cm);
\draw[thick] (46.3,-0.5)..controls +(0,-2) and +(0,+2) .. (43.7,-4); \draw[thick] (46.7,-0.5) -- (46.7,-4);

\node at (48.5,0) {$\Bigr)$};

\node at (20,-8) {$=\frac{1}{4}$};
\fill (22,-6) ellipse (0.6cm and 0.2cm);\fill (22,-10) ellipse (0.6cm and 0.2cm);
\draw[thick] (21.8,-6) -- (21.8,-10);\draw[thick] (22.2,-6) -- (22.2,-10);
\fill (25,-6) ellipse (0.6cm and 0.2cm);\fill (25,-10) ellipse (0.6cm and 0.2cm);
\draw[thick] (24.8,-6) -- (24.8,-10);\draw[thick] (25.2,-6) -- (25.2,-10);

\node at (27,-8) {$+\frac{1}{2}$};

\fill (29,-6) ellipse (0.6cm and 0.2cm);\fill (29,-10) ellipse (0.6cm and 0.2cm);
\draw[thick] (28.8,-6) -- (28.8,-10);\draw[thick] (29.2,-6)..controls +(0,-2) and +(0,+2) .. (31.8,-10);  
\fill (32,-6) ellipse (0.6cm and 0.2cm);\fill (32,-10) ellipse (0.6cm and 0.2cm);
\draw[thick] (31.8,-6)..controls +(0,-2) and +(0,+2) .. (29.2,-10);  \draw[thick] (32.2,-6) -- (32.2,-10);

\node at (34,-8) {$+\frac{1}{4}$};

\fill (36,-6) ellipse (0.6cm and 0.2cm);\fill (36,-10) ellipse (0.6cm and 0.2cm);
\draw[thick] (35.8,-6)..controls +(0,-2) and +(0,+2) .. (38.8,-10);  
\draw[thick] (36.2,-6)..controls +(0,-2) and +(0,+2) .. (39.2,-10);  
\fill (39,-6) ellipse (0.6cm and 0.2cm);\fill (39,-10) ellipse (0.6cm and 0.2cm);
\draw[thick] (38.8,-6)..controls +(0,-2) and +(0,+2) .. (35.8,-10);  
\draw[thick] (39.2,-6)..controls +(0,-2) and +(0,+2) .. (36.2,-10);  
\end{tikzpicture}
\end{center}
It is the analogue for fused permutations of an example we gave earlier for the fused Hecke algebra. We can appreciate that the multiplication is much simpler for fused permutations.
\end{exam}

\paragraph{Standard basis of the fused Hecke algebra.} By construction, the algebra of fused permutations has a (diagrammatic) basis indexed by fused permutations, namely by the set $\mS_{k,n}^{fus}$. Let us see how to ``lift'' this basis as a basis of the fused Hecke algebra. There is an interpretation (or a more precise formulation) of this basis in terms of distinguished representatives of double cosets. We only give here a rough description in diagrammatic terms.

Take an element of $\mS_{k,n}^{fus}$. This is an equivalence class of diagrams so let us choose a canonical representative. Recall that the only information which defines the equivalence class is which ellipse is connected to which one, the way in which an edge connect an ellipse to another does not matter. To define a canonical representative, first we forbid that three different edges intersect at the same point. Second, we require that the number of intersections (we do not want to call them crossings at this point, because the word is reserved for braid-like pictures) between edges is minimal.

So for each $w\in\mS_{k,n}^{fus}$, we have fixed a canonical representative diagram. Now in this diagram, promote each intersection between edges as a crossing, and we decide that all these crossings are positive (the edges coming from the left passes over the other one). We obtain the diagram of a fused braid, and thus an element of the fused Hecke algebra $H^{fus}_{k,n}(q)$. Denote it $T_w$. All these elements together form the standard basis of the fused Hecke algebra:
\[\{T_w\,,\ \ w\in\mS_{k,n}^{fus}\}\ \text{ is the standard basis of $H^{fus}_{k,n}(q)$.}\]
So we can say that the fused Hecke algebra $H^{fus}_{k,n}(q)$ is a (flat) deformation of the algebra of fused permutations $H^{fus}_{k,n}(1)$ since they both have bases indexed by the same set, and the multiplication in $H^{fus}_{k,n}(1)$ is a particular case ($q=1$) of the one in $H^{fus}_{k,n}(q)$.

\begin{exam} $\bullet$ Here is the standard basis of $H^{fus}_{2,2}(q)$:
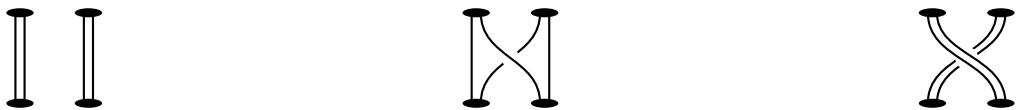
\begin{center}
\begin{tikzpicture}[scale=0.3]
\fill (1,2) ellipse (0.6cm and 0.2cm);\fill (1,-2) ellipse (0.6cm and 0.2cm);
\draw[thick] (0.8,2)..controls +(0,-2) and +(0,+2) .. (0.8,-2);\draw[thick] (1.2,2)..controls +(0,-2) and +(0,+2) .. (1.2,-2);  
\fill (4,2) ellipse (0.6cm and 0.2cm);\fill (4,-2) ellipse (0.6cm and 0.2cm);
\draw[thick] (3.8,2)..controls +(0,-2) and +(0,+2) .. (3.8,-2);\draw[thick] (4.2,2)..controls +(0,-2) and +(0,+2) .. (4.2,-2);

\fill (21,2) ellipse (0.6cm and 0.2cm);\fill (21,-2) ellipse (0.6cm and 0.2cm);
\fill (24,2) ellipse (0.6cm and 0.2cm);\fill (24,-2) ellipse (0.6cm and 0.2cm);
\draw[thick] (23.8,2)..controls +(0,-2) and +(0,+2) .. (21.2,-2); \draw[thick] (24.2,2)..controls +(0,-2) and +(0,+2) .. (24.2,-2);
\fill[white] (22.5,0) circle (0.4);
\draw[thick] (20.8,2) -- (20.8,-2);\draw[thick] (21.2,2)..controls +(0,-2) and +(0,+2) .. (23.8,-2);  

\fill (41,2) ellipse (0.6cm and 0.2cm);\fill (41,-2) ellipse (0.6cm and 0.2cm);
\fill (44,2) ellipse (0.6cm and 0.2cm);\fill (44,-2) ellipse (0.6cm and 0.2cm);
\draw[thick] (43.8,2)..controls +(0,-1.9) and +(0,+2.1) .. (40.8,-2);\draw[thick] (44.2,2)..controls +(0,-2) and +(0,+2) .. (41.2,-2); 
\fill[white] (42.5,0) circle (0.5);
\draw[thick] (40.8,2)..controls +(0,-2.1) and +(0,+1.9) .. (43.8,-2);\draw[thick] (41.2,2)..controls +(0,-1.9) and +(0,+2.1) .. (44.2,-2);
\end{tikzpicture}
\end{center}

$\bullet$ There are 21 distinct fused permutations in $\mS^{fus}_{2,3}$ and here are three examples. In each case, we draw a diagram in a canonical form, and the associated standard basis element of $H_{2,3}^{fus}(q)$ below it:
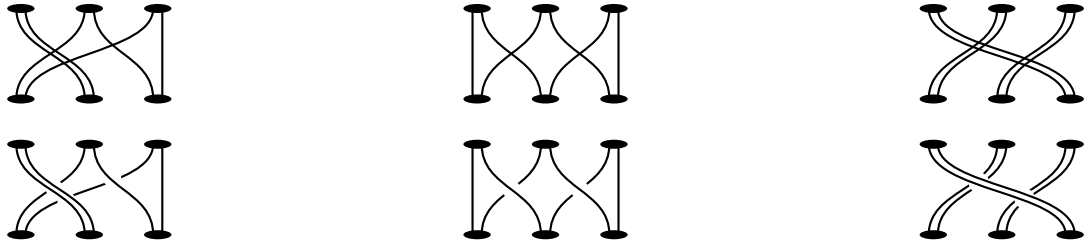
\begin{center}
\begin{tikzpicture}[scale=0.3]
\fill (1,2) ellipse (0.6cm and 0.2cm);\fill (1,-2) ellipse (0.6cm and 0.2cm);
\draw[thick] (0.8,2)..controls +(0,-2) and +(0,+2) .. (3.8,-2);\draw[thick] (1.2,2)..controls +(0,-2) and +(0,+2) .. (4.2,-2);  
\fill (4,2) ellipse (0.6cm and 0.2cm);\fill (4,-2) ellipse (0.6cm and 0.2cm);
\draw[thick] (3.8,2)..controls +(0,-2) and +(0,+2) .. (0.8,-2);\draw[thick] (4.2,2)..controls +(0,-2) and +(0,+2) .. (6.8,-2);
\fill (7,2) ellipse (0.6cm and 0.2cm);\fill (7,-2) ellipse (0.6cm and 0.2cm);
\draw[thick] (6.8,2)..controls +(0,-2) and +(0,+2) .. (1.2,-2);\draw[thick] (7.2,2) -- (7.2,-2);

\fill (1,-4) ellipse (0.6cm and 0.2cm);\fill (1,-8) ellipse (0.6cm and 0.2cm);
\fill (4,-4) ellipse (0.6cm and 0.2cm);\fill (4,-8) ellipse (0.6cm and 0.2cm);
\draw[thick] (3.8,-4)..controls +(0,-2) and +(0,+2) .. (0.8,-8);
\fill (7,-4) ellipse (0.6cm and 0.2cm);\fill (7,-8) ellipse (0.6cm and 0.2cm);
\draw[thick] (6.8,-4)..controls +(0,-2) and +(0,+2) .. (1.2,-8);\draw[thick] (7.2,-4) -- (7.2,-8);
\fill[white] (5,-5.5) circle (0.4);\fill[white] (2.5,-6) circle (0.4);\fill[white] (3,-6.5) circle (0.4);
\draw[thick] (4.2,-4)..controls +(0,-2) and +(0,+2) .. (6.8,-8);
\draw[thick] (0.8,-4)..controls +(0,-2) and +(0,+2) .. (3.8,-8);\draw[thick] (1.2,-4)..controls +(0,-2) and +(0,+2) .. (4.2,-8);  

\fill (21,2) ellipse (0.6cm and 0.2cm);\fill (21,-2) ellipse (0.6cm and 0.2cm);
\draw[thick] (20.8,2) -- (20.8,-2);\draw[thick] (21.2,2)..controls +(0,-2) and +(0,+2) .. (23.8,-2);  
\fill (24,2) ellipse (0.6cm and 0.2cm);\fill (24,-2) ellipse (0.6cm and 0.2cm);
\draw[thick] (23.8,2)..controls +(0,-2) and +(0,+2) .. (21.2,-2); \draw[thick] (26.8,2)..controls +(0,-2) and +(0,+2) .. (24.2,-2);
\fill (27,2) ellipse (0.6cm and 0.2cm);\fill (27,-2) ellipse (0.6cm and 0.2cm);
\draw[thick] (24.2,2)..controls +(0,-2) and +(0,+2) .. (26.8,-2);\draw[thick] (27.2,2) -- (27.2,-2);

\fill (21,-4) ellipse (0.6cm and 0.2cm);\fill (21,-8) ellipse (0.6cm and 0.2cm);
\draw[thick] (20.8,-4)-- (20.8,-8);
\fill (24,-4) ellipse (0.6cm and 0.2cm);\fill (24,-8) ellipse (0.6cm and 0.2cm);
\draw[thick] (23.8,-4)..controls +(0,-2) and +(0,+2) .. (21.2,-8); \draw[thick] (26.8,-4)..controls +(0,-2) and +(0,+2) .. (24.2,-8);
\fill (27,-4) ellipse (0.6cm and 0.2cm);\fill (27,-8) ellipse (0.6cm and 0.2cm);
\draw[thick] (27.2,-4) -- (27.2,-8);
\fill[white] (22.5,-6) circle (0.4);\fill[white] (25.5,-6) circle (0.4);
\draw[thick] (21.2,-4)..controls +(0,-2) and +(0,+2) .. (23.8,-8); 
\draw[thick] (24.2,-4)..controls +(0,-2) and +(0,+2) .. (26.8,-8);

\fill (41,2) ellipse (0.6cm and 0.2cm);\fill (41,-2) ellipse (0.6cm and 0.2cm);
\draw[thick] (40.8,2)..controls +(0,-2.1) and +(0,+1.9) .. (46.8,-2);\draw[thick] (41.2,2)..controls +(0,-1.9) and +(0,+2.1) .. (47.2,-2);
\fill (44,2) ellipse (0.6cm and 0.2cm);\fill (44,-2) ellipse (0.6cm and 0.2cm);
\draw[thick] (43.8,2)..controls +(0,-2) and +(0,+2) .. (40.8,-2);\draw[thick] (44.2,2)..controls +(0,-2) and +(0,+2) .. (41.2,-2); 
\fill (47,2) ellipse (0.6cm and 0.2cm);\fill (47,-2) ellipse (0.6cm and 0.2cm);
\draw[thick] (46.8,2)..controls +(0,-2) and +(0,+2) .. (43.8,-2);\draw[thick] (47.2,2)..controls +(0,-2) and +(0,+2) .. (44.2,-2);

\fill (41,-4) ellipse (0.6cm and 0.2cm);\fill (41,-8) ellipse (0.6cm and 0.2cm);
\fill (44,-4) ellipse (0.6cm and 0.2cm);\fill (44,-8) ellipse (0.6cm and 0.2cm);
\draw[thick] (43.8,-4)..controls +(0,-2) and +(0,+2) .. (40.8,-8);\draw[thick] (44.2,-4)..controls +(0,-2) and +(0,+2) .. (41.2,-8); 
\fill (47,-4) ellipse (0.6cm and 0.2cm);\fill (47,-8) ellipse (0.6cm and 0.2cm);
\draw[thick] (46.8,-4)..controls +(0,-2) and +(0,+2) .. (43.8,-8);\draw[thick] (47.2,-4)..controls +(0,-2) and +(0,+2) .. (44.2,-8);
\fill[white] (43,-5.7) circle (0.45);\fill[white] (45,-6.4) circle (0.45);
\draw[thick] (40.8,-4)..controls +(0,-2.1) and +(0,+1.9) .. (46.8,-8);\draw[thick] (41.2,-4)..controls +(0,-1.9) and +(0,+2.1) .. (47.2,-8);
\end{tikzpicture}
\end{center}
\end{exam}

\subsection{Conclusion}

\paragraph{Braid group.}

We denote by $\Sigma_i$ the fused braid in $H_{k,n}^{fus}(q)$ for which all strands starting from ellipse $i$ pass over the strands starting from ellipse $i+1$ and all other strands are vertical. A picture is better here, so for example, for $k=2$:
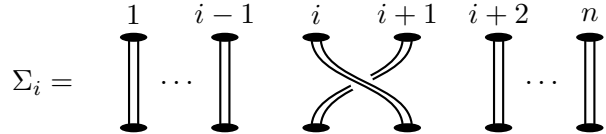
\begin{center}
 \begin{tikzpicture}[scale=0.3]
\node at (-2,0) {$\Sigma_i=$};
\node at (2,3) {$1$};\fill (2,2) ellipse (0.6cm and 0.2cm);\fill (2,-2) ellipse (0.6cm and 0.2cm);
\draw[thick] (1.8,2) -- (1.8,-2);\draw[thick] (2.2,2) -- (2.2,-2);
\node at (4,0) {$\dots$};
\draw[thick] (5.8,2) -- (5.8,-2);\draw[thick] (6.2,2) -- (6.2,-2);
\node at (6,3) {$i-1$};\fill (6,2) ellipse (0.6cm and 0.2cm);\fill (6,-2) ellipse (0.6cm and 0.2cm);
\node at (10,3) {$i$};\fill (10,2) ellipse (0.6cm and 0.2cm);\fill (10,-2) ellipse (0.6cm and 0.2cm);
\node at (14,3) {$i+1$};\fill (14,2) ellipse (0.6cm and 0.2cm);\fill (14,-2) ellipse (0.6cm and 0.2cm);

\draw[thick] (13.8,2)..controls +(0,-2) and +(0,+2) .. (9.8,-2);
\draw[thick] (14.2,2)..controls +(0,-2) and +(0,+2) .. (10.2,-2);
\fill[white] (12,0) circle (0.5);
\draw[thick] (10.2,2)..controls +(0,-2) and +(0,+2) .. (14.2,-2);
\draw[thick] (9.8,2)..controls +(0,-2) and +(0,+2) .. (13.8,-2);

\draw[thick] (17.8,2) -- (17.8,-2);\draw[thick] (18.2,2) -- (18.2,-2);
\node at (18,3) {$i+2$};\fill (18,2) ellipse (0.6cm and 0.2cm);\fill (18,-2) ellipse (0.6cm and 0.2cm);
\node at (20,0) {$\dots$};
\draw[thick] (21.8,2) -- (21.8,-2);\draw[thick] (22.2,2) -- (22.2,-2);\fill (22,2) ellipse (0.6cm and 0.2cm);\fill (22,-2) ellipse (0.6cm and 0.2cm);
\node at (22,3) {$n$};
\end{tikzpicture}
\end{center}
These elements satisfy the braid relation in $H_{k,n}^{fus}(q)$:
\[\Sigma_i\Sigma_{i+1}\Sigma_i=\Sigma_{i+1}\Sigma_i\Sigma_{i+1}\ .\]
Recall that the multiplication in the fused Hecke algebra  $H_{k,n}^{fus}(q)$ is quite intricate due to the appearances of the $q$-symmetriser replacing the middle ellipses each time we concatenate two fused braids. However, in the particular situation as for $\Sigma_i$, where the strands arriving at an ellipse are all coming from the same ellipse and are not intertwined with anything, then the $q$-symmetriser that will appear in place of this ellipse can be moved up along these strands and sent into an ellipse on top. Due to the idempotent relations in $H_{k,n}^{fus}(q)$, it will disappear. In other words, when we multiply various $\Sigma_i$ following the multiplication rule for the fused Hecke algebra, it reduces to simple concatenations. So the braid relations for the generators $\Sigma_i$ follows from the following kind of obvious relations in the usual braid group\footnote{The picture looks like some kind of tubes, but we really mean a relation for usual braids (here with 6 strands)} (here $k=2$, but simply add more parallel strands for arbitrary $k$):
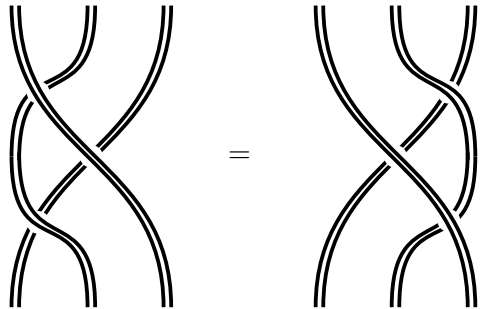
\begin{center}
\begin{tikzpicture}[scale=0.5]
\draw[line width=0.5mm] (4,2)..controls +(0,-4) and +(0,+4) .. (0,-6);
\draw[line width=0.5mm] (4.2,2)..controls +(0,-4) and +(0,+4) .. (0.2,-6);
\draw[line width=0.5mm] (2,2)..controls +(0,-3) and +(0,+3) .. (0,-2);
\draw[line width=0.5mm] (2.2,2)..controls +(0,-3) and +(0,+3) .. (0.2,-2.1);
\fill[white] (0.7,-0.3) circle (0.3);
\fill[white] (2.1,-2) circle (0.3);
\draw[line width=0.5mm] (0,2)..controls +(0,-4) and +(0,+4) .. (4,-6);
\draw[line width=0.5mm] (0.2,2)..controls +(0,-4) and +(0,+4) .. (4.2,-6);
\fill[white] (0.7,-3.75) circle (0.3);
\draw[line width=0.5mm] (0,-2)..controls +(0,-3) and +(0,+3) .. (2,-6);
\draw[line width=0.5mm] (0.2,-1.9)..controls +(0,-3) and +(0,+3) .. (2.2,-6);
\node at (6,-2) {$=$};
\draw[line width=0.5mm] (12,2)..controls +(0,-4) and +(0,+4) .. (8,-6);
\draw[line width=0.5mm] (12.2,2)..controls +(0,-4) and +(0,+4) .. (8.2,-6);
\fill[white] (11.5,-0.3) circle (0.3);
\fill[white] (10.1,-2) circle (0.3);
\draw[line width=0.5mm] (10,2)..controls +(0,-3) and +(0,+3) .. (12,-2.1);
\draw[line width=0.5mm] (10.2,2)..controls +(0,-3) and +(0,+3) .. (12.2,-2);
\draw[line width=0.5mm] (12,-1.9)..controls +(0,-3) and +(0,+3) .. (10,-6);
\draw[line width=0.5mm] (12.2,-2)..controls +(0,-3) and +(0,+3) .. (10.2,-6);
\fill[white] (11.5,-3.7) circle (0.3);
\draw[line width=0.5mm] (8,2)..controls +(0,-4) and +(0,+4) .. (12,-6);
\draw[line width=0.5mm] (8.2,2)..controls +(0,-4) and +(0,+4) .. (12.2,-6);
\end{tikzpicture}
\end{center}
In addition to the braid relations, the elements $\Sigma_1,\dots,\Sigma_{n-1}$ must also satisfy some characteristic equation since the algebra $H^{fus}_{k,n}(q)$ is finite-dimensional. In fact, for $n=2$, it is easy to see that the dimension is $k+1$ (we have $k+1$ different fused permutations) and that the element $\Sigma_1$ generates $H^{fus}_{k,2}(q)$. So the characteristic equation will be of degree $k+1$. To calculate the eigenvalues directly from the defining multiplication in  $H^{fus}_{k,2}(q)$ seems to be a little intricate, but we can use the representation of $H^{fus}_{k,2}(q)$ on the tensor product $S^kV\otimes S^kV$ of representations of $U_q(sl_N)$ (for any $N$). In this representation, the element $\Sigma_1$ corresponds to the $R$-matrix, and from this, one can obtain the characteristic equation. 

To summarize, the elements $\Sigma_1,\dots,\Sigma_{n-1}$ of $H^{fus}_{k,n}(q)$ satisfy the following relations:
\[\begin{array}{ll}
\Sigma_i\Sigma_{i+1}\Sigma_i=\Sigma_{i+1}\Sigma_i\Sigma_{i+1}\,,\ \ \  & \text{for $i\in\{1,\dots,n-2\}$}\,,\\[0.2em]
\Sigma_i\Sigma_j=\Sigma_j\Sigma_i\,,\ \ \  & \text{for $i,j\in\{1,\dots,n-1\}$ such that $|i-j|>1$}\,,\\[0.2em]
\prod_{l=0}^k\Bigl(\Sigma_i-(-1)^{k+l} q^{-k+l(l+1)}\Bigr)=0\,,\ \ \  & \text{for $i\in\{1,\dots,n-1\}$}\,.
\end{array}
\]
So inside the fused Hecke algebra $H^{fus}_{k,n}(q)$ the elements $\Sigma_1,\dots,\Sigma_{n-1}$ generate a subalgebra which is a quotient of the algebra of the braid group. In this quotient, the elementary braidings satisfy a characteristic equation of order $k+1$. 

It is important to note that these elements do not generate the whole algebra $H^{fus}_{k,n}(q)$ when $n>2$ and $k>1$. This is a striking difference compared to the usual Hecke algebra. We also note that other relations than the ones above (and not implied by the ones above) must be satisfied by the elements $\Sigma_1,\dots,\Sigma_{n-1}$ since the algebra they generate is finite-dimensional, and we know that relations above are not enough to define a finite-dimensional algebra. A description of the subalgebra generated by $\Sigma_1,\dots,\Sigma_{n-1}$, such as a set of defining relations, is an interesting open question.

\begin{rema}
The fact that $\Sigma_1,\dots,\Sigma_{n-1}$ satisfy the braid relations in $H^{fus}_{k,n}(q)$ is reminiscent of what happened in the fusion procedure (in fact it is really the fusion procedure without spectral parameters). Indeed, $\Sigma_1,\dots,\Sigma_{n-1}$, seen as elements of $H_{kn}(q)$, satisfy the braid relation (first step of the fusion procedure). And it turns out that they commute with the idempotent $P_{k,n}$ (second step of the fusion procedure) so they still satisfy the braid relation in $H^{fus}_{k,n}(q)=P_{k,n}H_{kn}(q)P_{k,n}$.
\end{rema}

\paragraph{YB equation.} Having identified the realisation of the braid group inside the fused Hecke algebra, now we go on in our wish list of the introduction, and ask for a solution of the YB equation. The main result is that the fused Hecke algebra admits a Baxterization formula. Of course, if one follows the fusion procedure of the preceding section, we already know that there must be a solution of the YB equation somewhere in $H^{fus}_{k,n}(q)$. The point here is to give this solution with an explicit Baxterization formula. 

First, we introduce some natural elements in the algebra $H^{fus}_{k,n}(q)$: the partial elementary braiding is denoted $\Sigma_i^{(p)}$ and corresponds to the fused braid for which the $p$ rightmost strands starting from ellipse $i$ pass over the $p$ leftmost strands starting from ellipse $i+1$, and all other strands are vertical. Again a picture is better (here $k=3$):
\begin{center}
 \begin{tikzpicture}[scale=0.25]
\node at (16,0) {$\Sigma_i^{(1)}=$};
\node at (20,0) {$\dots$};
\node at (22,3) {$i$};\fill (22,2) ellipse (0.8cm and 0.2cm);\fill (22,-2) ellipse (0.8cm and 0.2cm);
\node at (26,3) {$i+1$};\fill (26,2) ellipse (0.8cm and 0.2cm);\fill (26,-2) ellipse (0.8cm and 0.2cm);
\draw[thick] (21.7,2) -- (21.7,-2);
\draw[thick] (22,2) -- (22,-2);
\draw[thick] (25.7,2)..controls +(0,-2) and +(0,+2) .. (22.3,-2);
\fill[white] (24,0) circle (0.4);
\draw[thick] (22.3,2)..controls +(0,-2) and +(0,+2) .. (25.7,-2);
\draw[thick] (26,2) -- (26,-2);
\draw[thick] (26.3,2) -- (26.3,-2);
\node at (28,0) {$\dots$};
\node at (29.5,0) {$,$};

\node at (34,0) {$\Sigma_i^{(2)}=$};
\node at (38,0) {$\dots$};
\node at (40,3) {$i$};\fill (40,2) ellipse (0.8cm and 0.2cm);\fill (40,-2) ellipse (0.8cm and 0.2cm);
\node at (44,3) {$i+1$};\fill (44,2) ellipse (0.8cm and 0.2cm);\fill (44,-2) ellipse (0.8cm and 0.2cm);
\draw[thick] (39.7,2) -- (39.7,-2);
\draw[thick] (43.7,2)..controls +(0,-2) and +(0,+2) .. (40,-2);
\draw[thick] (44,2)..controls +(0,-2.5) and +(0,+2) .. (40.3,-2);
\fill[white] (42,-0.1) ellipse (0.6cm and 0.4cm);
\draw[thick] (40,2)..controls +(0,-2.5) and +(0,+2) .. (43.7,-2);
\draw[thick] (40.3,2)..controls +(0,-2) and +(0,+2) .. (44,-2);
\draw[thick] (44.3,2) -- (44.3,-2);
\node at (46,0) {$\dots$};
\node at (47.5,0) {$,$};

\node at (52,0) {$\Sigma_i^{(3)}=$};
\node at (56,0) {$\dots$};
\node at (58,3) {$i$};\fill (58,2) ellipse (0.8cm and 0.2cm);\fill (58,-2) ellipse (0.8cm and 0.2cm);
\node at (62,3) {$i+1$};\fill (62,2) ellipse (0.8cm and 0.2cm);\fill (62,-2) ellipse (0.8cm and 0.2cm);
\draw[thick] (61.7,2)..controls +(0,-2) and +(0,+2) .. (57.7,-2);
\draw[thick] (62,2)..controls +(0,-2.5) and +(0,+2) .. (58,-2);
\draw[thick] (62.3,2)..controls +(0,-3) and +(0,+2) .. (58.3,-2);
\fill[white] (60,-0.2) ellipse (1cm and 0.6cm);
\draw[thick] (57.7,2)..controls +(0,-3) and +(0,+2) .. (61.7,-2);
\draw[thick] (58,2)..controls +(0,-2.5) and +(0,+2) .. (62,-2);
\draw[thick] (58.3,2)..controls +(0,-2) and +(0,+2) .. (62.3,-2);
\node at (64,0) {$\dots$};
\node at (65.5,0) {$,$};

\end{tikzpicture}
\end{center}
where all strands starting from ellipses $1,\dots,i-1,i+2,\dots,n$ are vertical. The element $\Sigma_i^{(0)}$ is the identity element, while $\Sigma_i^{(k)}$ are the elements called $\Sigma_i$ satisfying the braid relations in the previous paragraph.

In addition to the $q$-numbers $[L]_q:=\frac{q^L-q^{-L}}{q-q^{-1}}$ and the $q$-factorial $[L]_q!:=[1]_q[2]_q\dots[L]_q$, we also define the $q$-binomials and the $q$-Pochhammer symbol:
\[
  \left[\begin{array}{c}L \\p\end{array}\right]_q:=\frac{[L]_q!}{[L-p]_q![p]_q!}\ ,\ \ \ \ \ \  (a\, ;\, q)_p=\prod_{r=0}^{p-1}(1-aq^r)\ .
\]
By convention, we have $[0]_q!=\left[\begin{array}{c}L \\0\end{array}\right]_q=(a\, ;\, q)_0=1$. 

\begin{theo}[\cite{CPdA2}]\label{th:bax}
The following function taking values in $H^{fus}_{k,n}(q)$
 \begin{equation}
 \check R_i( u )=\sum_{p=0}^k (-q)^{k-p}\left[\begin{array}{c}k \\p\end{array}\right]^2_q\ 
 \frac{ (q^{-2}\, ;\, q^{-2})_{k-p}  }{ (uq^{-2p}\,;\, q^{-2})_{k-p}}\ \Sigma_i^{(p)}\ , \label{eq:R}
 \end{equation}
satisfy the braided Yang--Baxter equations:
\[
\check  R_i(u)\check R_{i+1}(uv)\check R_i(v)=\check R_{i+1}(v) \check R_i(uv) \check R_{i+1}(v)\ .
\]
\end{theo}
Note that the (rational) dependence on the spectral parameter $u$ is clearly visible from this formula in $H_{k,n}^{fus}(q)$. 

\begin{exam}
For $k=1$, the formula above is the usual Baxterization formula of the Hecke algebra. For $k=2$, the formula is:
\[
\check R_i(u)=\Sigma_i^{(2)} - (q+q^{-1})\frac{(q^2-q^{-2}) }{1-uq^{-2}}\Sigma_i^{(1)} +q^2 \frac{(1-q^{-2})(1-q^{-4})}{(1-u)(1-uq^{-2})} \ .
\]
\end{exam}

\paragraph{Representation theory.} Let $W=S_q^k(V)$ be the $q$-symmetrised power of the vector representation of $U_q(sl_N)$. From the construction of $H^{fus}_{k,n}(q)$, it can be represented on $W^{\otimes n}$. So our situation again is that we have two algebras represented on the same vector space:
\[U_q(sl_N)\ \ \stackrel{\rho^{(n)}}{\longrightarrow}\ \ \ \ \text{End}(W^{\otimes n})\ \ \ \ \stackrel{\pi}{\longleftarrow}\ \ H^{fus}_{k,n}(q)\,,\ \]
With the formulation of $H^{fus}_{k,n}(q)$ using the idempotent as in (\ref{def-Hfus}), the next result follows almost by construction. 
\begin{theo}[Schur--Weyl duality for symmetrised power \cite{CPdA1}]
The centraliser $\text{End}_{U_q(sl_N)}(W^{\otimes n})$ is the image of the fused Hecke algebra $H^{fus}_{k,n}(q)$:
\[\text{End}_{U_q(sl_N)}(W^{\otimes n})=\pi\bigl(H_{k,n}^{fus}(q)\bigr)\ .\]
\end{theo}
It turns out that one may study directly (without assuming any knowledge on $U_q(sl_N)$-representations) the representation theory of $H^{fus}_{k,n}(q)$ and recover the well-known decomposition of the $U_q(sl_N)$ representation on $W^{\otimes n}$:
\[W^{\otimes n}=\bigoplus_{\substack{\lambda\vdash kn\\[0.3em]
\ell(\lambda)\leq N}} \bigl(L^N_{\lambda}\bigr)^{\oplus K_{\lambda,(k^n)}}\ ,
\]
where $K_{\lambda,(k^n)}$ is the Kostka number counting the number of semistandard Young tableaux of shape $\lambda$ containing $k$ times each number $1,\dots,n$. This is a particular case of the Littlewood--Richardson rule. These Kostka numbers are recovered as the dimensions of the irreducible representations of $H^{fus}_{k,n}(q)$. 

\vskip .2cm
To describe fully the centraliser $\text{End}_{U_q(sl_N)}(W^{\otimes n})$, we need to understand the kernel of the representation $\pi$ of $H^{fus}_{k,n}(q)$. More details are in \cite{CPdA1}.

\vskip 1.5cm
We conclude with the following summary:
\begin{conc}[Final]
We have constructed an algebra $H^{fus}_{k,n}(q)$ called the fused Hecke algebra:
\begin{itemize}
\item[$\bullet$] it contains elements $\Sigma_1,\dots,\Sigma_{n-1}$ satisfying the braid relations;
\item[$\bullet$] it admits an explicit Baxterization formula;
\item[$\bullet$] it has representations on vector spaces $W^{\otimes n}$, where $W=S_q^k(V)$, and through such representations it is in Schur--Weyl duality with $U_q(sl_N)$, where $N=\dim(V)$.
\end{itemize}
\end{conc}

\end{document}